\documentclass[11pt,preprint]{article}

\usepackage{jmlr2e}

\usepackage[ruled,vlined]{algorithm2e}
\usepackage{amsmath,amsfonts,amssymb}
\usepackage{graphicx,psfrag,epsf}
\usepackage{enumerate}
\usepackage{natbib}
\usepackage{url} 
\allowdisplaybreaks

\newcommand{\bbR}{\mathbb{R}}

\newcommand{\sD}{\mathcal{D}}
\newcommand{\sE}{\mathcal{E}}

\newcommand{\sN}{\mathcal{N}}

\newcommand{\sL}{\mathcal{L}}

\newcommand{\sS}{\mathcal{S}}
\newcommand{\sT}{\mathcal{T}}

\newcommand{\sP}{\mathcal{P}}
\newcommand{\sQ}{\mathcal{Q}}
\newcommand{\sU}{\mathcal{U}}
\newcommand{\sV}{\mathcal{V}}

\newcommand{\vb}{\mathbf{b}}

\newcommand{\ve}{\mathbf{e}}
\newcommand{\vf}{\mathbf{f}}

\newcommand{\vs}{\mathbf{s}}

\newcommand{\vu}{\mathbf{u}}

\newcommand{\vx}{\mathbf{x}}

\newcommand{\valpha}{\boldsymbol{\alpha}}

\newcommand{\vphi}{\boldsymbol{\phi}}
\newcommand{\vvarphi}{\boldsymbol{\varphi}}

\newcommand{\vzero}{\mathbf{0}}

\newcommand{\vbeta}{\boldsymbol{\beta}}
\newcommand{\veta}{\boldsymbol{\eta}}
\newcommand{\mA}{\mathbf{A}}
\newcommand{\mB}{\mathbf{B}}

\newcommand{\mD}{\mathbf{D}}

\newcommand{\mF}{\mathbf{F}}

\newcommand{\mI}{\mathbf{I}}

\newcommand{\mK}{\mathbf{K}}

\newcommand{\mM}{\mathbf{M}}
\newcommand{\mN}{\mathbf{N}}

\newcommand{\mQ}{\mathbf{Q}}

\newcommand{\mU}{\mathbf{U}}
\newcommand{\mV}{\mathbf{V}}
\newcommand{\mW}{\mathbf{W}}
\newcommand{\mX}{\mathbf{X}}

\newcommand{\mGamma}{\boldsymbol{\Gamma}}
\newcommand{\mOmega}{\boldsymbol{\Omega}}

\newcommand{\mDelta}{\boldsymbol{\Delta}}
\newcommand{\mSigma}{\boldsymbol{\Sigma}}
\newcommand{\mPi}{\boldsymbol{\Pi}}

\newcommand{\bvb}{\bar{\mathbf{b}}}
\newcommand{\bmA}{\bar{\mathbf{A}}}
\newcommand{\bmB}{\bar{\mathbf{B}}}

\newcommand{\bmD}{\bar{\mathbf{D}}}
\newcommand{\bmU}{\bar{\mathbf{U}}}

\newcommand{\bmV}{\bar{\mathbf{V}}}

\newcommand{\Expect}{\mathbb{E}}
\newcommand{\intd}{\,\mathrm{d}}
\newcommand{\diag}{\mathrm{diag}}

\newcommand{\tr}{\mathrm{tr}}
\newcommand{\Var}{\mathrm{Var}}
\newcommand{\rank}{\mathrm{rank}}

\def\trans{^\mathsf{T}}

\DeclareMathOperator*{\argmin}{arg\,min}
\newtheorem{condition}{Condition}

\usepackage{mathtools}
\DeclarePairedDelimiter\floor{\lfloor}{\rfloor}

\usepackage{color}

\usepackage{multirow}
\usepackage{color}
\usepackage{accents}
\newcommand{\ubar}[1]{\underaccent{\bar}{#1}}

\usepackage{comment}

\jmlrheading{1}{2022}{1--69}{4/00}{10/00}{}{}

\ShortHeadings{A Unified Analysis of 
	Multi-task Functional Linear Regression Models}{He, Ye and He}
\firstpageno{1}

\begin{document}

\title{A Unified Analysis of 
	Multi-task Functional Linear Regression Models with Manifold Constraint and Composite Quadratic Penalty}

\author{\name Shiyuan He \email  heshiyuan@ruc.edu.cn 
        \AND
       \name Hanxuan Ye \email ye86171958@ruc.edu.cn 
        \AND
    \name Kejun He \thanks{Corresponding author.}\email kejunhe@ruc.edu.cn\\
       \addr The Center for Applied Statistics, Institute of Statistics and Big Data\\
       Renmin University of China\\
       Beijing, 100872, China}
\editor{ }

\maketitle

\begin{abstract}
This work studies the multi-task functional linear regression models where both the covariates and the unknown regression coefficients (called slope functions) are curves.  For slope function estimation, we employ penalized splines to balance bias, variance, and computational complexity. The power of multi-task learning is brought in by imposing additional structures over the slope functions. We propose a general model with double regularization over the spline coefficient matrix: i) a matrix manifold constraint, and ii) a composite penalty as a summation of quadratic terms. 
Many multi-task learning approaches can be treated as special cases of this proposed model, such as a reduced-rank model and a graph Laplacian regularized model.
We show the composite penalty induces a specific norm, which helps quantify the manifold curvature and determine the corresponding proper subset in the manifold tangent space. The complexity of tangent space subset is then bridged to the complexity of geodesic neighbor via generic chaining. A unified upper bound of the convergence rate is obtained and specifically applied to the reduced-rank model and the graph Laplacian regularized model. The phase transition behaviors for the estimators are examined as we vary the configurations of model parameters.
\end{abstract}



\section{Introduction}\label{sec:intro}

Multi-task learning has been extensively adopted in various machine learning areas, including linear regression~\citep{solnon2012multi}, classification~\citep{cavallanti2010linear}, neural networks~\citep{crawshaw2020multi}, clustering~\citep{zhang2014convex}, and reinforcement learning \citep{teh2017distral}. 
By leveraging the shared information to learn multiple related tasks simultaneously, multi-task learning becomes an effective approach to improve the overall generalization performance of tasks. 
Its theoretical benefits were investigated in \citep{baxter2000model} under a class of probably approximately correct (PAC) models, showing the average estimation error of tasks can potentially decrease with the number of tasks. 
Multi-task learning can be achieved by different strategies, such as restricting model rank \citep{velu2013multivariate}, encouraging shared feature \citep{kolar2011union}, and learning common representation \citep{maurer2016benefit}. 
A comprehensive overview can be found in~\cite{thung2018brief} and~\cite{zhang2018overview}.

Existing literature of multi-tasking learning usually assumes the input space is a $d$-dimensional Euclidean space (though typically high-dimensional).
This work, however, considers the class of scalar-on-function regressions. The scalar-on-function regression has mostly been studied as a single task. 
The commonly-used model \citep{cardot1999functional}, known as functional linear regression, predicts a random variable $Y\in\mathbb{R}$ by a covariate curve $X(t)$, which is a random function over an interval $\sT$. 
The linear prediction is based on the integrated quantity $\alpha + \int_{\sT} X(t) \beta(t) \intd t$, where $\alpha$ and $\beta(\cdot)$ are the intercept and slope function, respectively. 
Abundant works for univariate $Y$ have been studied on different functional linear regression models, e.g., the least squares regression \citep{cardot03,yao2005functional, hall07,yuan2010reproducing}, generalized exponential family regression \citep{dou2012estimation}, and quantile regression \citep{Kato12}. 
Nevertheless, directly applying the above work to the multi-task applications will result in short of efficiency, since the intrinsic relatedness between tasks is ignored. 
Our work attempts to extend a broad class of functional linear regression models to the multi-task setting.  

Multi-task scalar-on-function regression models can provide wide applications in real world. 
In astronomy \citep{blanco2014determining}, researchers need to determine multiple atmospheric parameters (e.g., effective temperature, surface gravity, metallicity) and individual chemical abundances from high-resolution stellar spectrum. 
Each spectrum can be treated as a functional covariate (viewed as a function of wavelength), and the stellar parameters and chemical abundances correspond to multiple scalar responses. 
Besides, the study of Alzheimer's disease shows that multiple cognitive and memory scores \citep{li2016sparse} can be potentially predicted from neural imaging by functional regression technique \citep[e.g.,][]{wang2014regularized}. 
There are also some applications where the datasets are collected from several locations. 
For example, \cite{ramsay05} predicted total annual precipitation for some Canadian weather stations from yearly temperature variation; and \cite{jiang2020functional} predicted the death rate caused by cardiovascular disease by the annual curves of air pollutant for several cities. 
Such kind of problems can be naturally formulated to be multi-task when a functional linear regression model is conducted for each location (viewed as a task). 
The slope functions of all locations can be expected to share similarity depending on the spatial proximity.

Suppose there are $M$ tasks with output variables $Y_1, \dots, Y_M \in\bbR$ and their associated functional covariates $X_1(t),\cdots, X_M(t)$ over the common domain $\sT$.
The first step towards modeling the multi-task functional linear regression is to represent the slope functions $\beta_1, \dots, \beta_M$ in an appropriate space. 
One potential choice is to model the slope functions in a reproducing kernel Hilbert space \citep[RKHS,][]{yuan2010reproducing} or by natural cubic splines \citep{crambes09}. 
This approach, known as smoothing splines, has computational cost as high as the cubic order of sample size (or observation grid size). 
Alternatively, one can also use a regression spline (like B-spline) space $\mathbb{S}_K$ with degrees of freedom $K$ and order $\mathfrak{o}+1$. 
In this approach, the degree of freedom $K$ is usually set to be a relatively small value, compared with the sample size. 
A small $K$ significantly reduces the computational cost, but it also increases the approximation bias. 
As a compromising solution, penalized splines \citep{cardot03} use a moderately large $K$ to reduce bias and employ a roughness penalty  to control the model complexity. 
Penalized splines can properly balance the computational cost, model bias and variance. Though penalized splines have many appealing practical characteristics, analyzing its rates of convergence and  phase transition behavior  is much more challenging. 
In the context of non-parametric regression \citep[e.g.,][]{claeskens2009asymptotic, kauermann2009some, xiao2019asymptotic, huang2021asymptotic} and covariance function estimation \citep{xiao2020asymptotic}, the corresponding penalized spline estimator is known to exhibit distinct rates of convergence when we vary the spline order, degrees of freedom, the penalty derivative order, and the penalty tuning parameter, 
as the sample size goes to infinity.  However, little is known about the theoretical properties of the penalized spline estimator in a functional linear regression model.

After using the penalized splines, the second challenge   is  to impose additional structures to obtain an improved estimation for multi-task functional linear regression. One possible approach is to assume the slope functions $\beta_1,\cdots, \beta_M$ come from an unknown subspace. 
Let $\vphi = (\phi_1,\cdots, \phi_K)\trans$ denote the vector of basis functions in the spline space $\mathbb{S}_K$. 
As each slope function rewritten as $\beta_m = \phi\trans\vb_m$ with coefficient vector $\vb_m \in \bbR^{K}$, $m=1,\dots, M$, the subspace assumption amounts to forcing $\mB=(\vb_1,\cdots, \vb_M) \in \bbR^{K\times M}$ to reside on a fixed-rank matrix manifold. This is called the reduced multi-task model in our work. 
Other structures are also easy to be imposed in our multi-task setting with penalized splines.
For example, sometimes, an external graph structure is available where the relation between tasks is encoded.  In the graph, each task is treated as a vertex, and the closeness between tasks is represented by edge weight.  The idea of graph Laplacian regularization \citep{evgeniou2005learning,zhu2015subspace,yousefi2018local} can be adopted to encourage similarity of slope functions between contiguous tasks. This is referred to as the graph regularized multi-task model.

Motivated by the above two special models, we propose a general model with double regularization over the spline coefficient matrix $\mB$. 
The first regularization is constraining $\mB$ to an embedded matrix submanifold $\mathbb{M}$ ($\subseteq\bbR^{K\times M}$). 
The second regularization is through a composite quadratic penalization $\sP_{\veta}(\mB) = \sum_{j=1}^{P} \eta_j\tr\big(\mB\trans\mPi_{j1}\mB\mPi_{j2}\big)$, where $\mPi_{j1}$'s and $\mPi_{j2}$'s are positive semi-definite matrices and $\eta_j$'s are penalty parameters. 
We will show that both the reduced multi-task model and the graph regularized multi-task model can be treated as special cases of the proposed model with double regularization.  Our proposed method is a more general model in the sense that the manifold constraint $\mathbb{M}$ in the first regularization, and the number $P$ and positive semi-definite matrices $\mPi_{j1}$'s and $\mPi_{j2}$'s in the second regularization are not specified. Overall, in this work, we aim to develop a unified treatment to the general model, and  provide a set of analysis tools that allows for the easy uncovering of the estimator's asymptotic properties and phase transition behaviors.

\subsection{Contributions of This Work}

Understanding the theoretical properties of the proposed model is far from being straightforward because the model consists of several non-trivial components. 
First, the general manifold constraint makes an explicit solution unavailable. 
Second, studying the estimator with a general composite quadratic penalty has also been known to be challenging. 
Special cases include the roughness penalty for penalized splines \citep{huang2021asymptotic} and the graph Laplacian penalty \citep{green2021minimax} for spatial similarity. 
To our best knowledge, little is known about the phase transition behavior of the penalized spline estimator in the setting of functional linear regression. 
As for the graph Laplacian regularization, although it has been widely adopted in applications, its theoretical study is mostly limited to non-parametric regression models of single task \citep{kirichenko2017estimating,green2021minimax,garcia2020maximum} and multiple tasks \citep{yousefi2018local}. 
The graph Laplacian regularized estimator in the context of  functional linear regression has been barely explored, especially for its phase transition behavior.

This work attempts to overcome the above challenges and includes several contributions to the literature. The first contribution of this work is integrating the penalized spline regularization and the Laplacian regularization into a general framework, namely the composite quadratic regularization. We develop a unified solution and a comprehensive analysis method for this general regularization. 
In the literature, it is well known that the estimation error is closely connected to local model complexity \citep{bartlett2005local,yousefi2018local}. However, the question remains for specifying an appropriate local set for the estimator with the composite quadratic penalty. 
We find that properly characterizing an ellipsoid-like neighbor (see~\eqref{eqn:defineQnorm} and~\eqref{eqn:generalNeighbor}) in the parameter space is a key element for understanding the composite quadratic penalty, where the ellipsoid-like neighbor is induced by the penalty itself. 
Via utilizing the generic chaining technique \citep{talagrand2014upper}, we are able to effectively evaluate model complexity and reveal the phase transition, which would be difficult to access via other existing analysis routines.

The second important contribution of our work is to allow the model parameter to be simultaneously constrained  by a \textit{general} submanifold (including but not limited to a low-rank manifold). A unified treatment is  also provided for the empirical process on the manifold constraint through the generic chaining technique.  With the penalty induced norm, we quantify the manifold curvature via its second fundamental form. When the manifold curvature is restricted, we show the complexity of the ellipsoid-like neighbor in manifold tangent space can be bridged to  that of its geodesic neighbor. The empirical processes can therefore be controlled by quantifying the complexity of local tangent space (see Section~\ref{sec:manifold}).  Restricting manifold curvature also allows us to bound the estimator perturbation by the magnitude of the corresponding tangent vector, which is also measured by the penalty induced norm (see Lemma~\ref{lemma:metriccompatible}).

This work also includes several novel contributions to the penalized spline literature. Our theoretical analysis tools lead to a comprehensive analysis of penalized splines as sample size $N$ goes to infinity, under extensive settings of the spline order $\mathfrak{o}+1$, the spline degrees of freedom $K$, the penalty derivative order $d$, the penalty parameter $\eta_1$, and the smoothness order $\nu$ of the slope function.  Our analysis is more involved than the existing works on non-parametric regression \citep[][]{huang2021asymptotic} and covariance function estimation \citep{xiao2020asymptotic}, since the smoothness and eigenvalue decay rate of the covariance function of $X_m$ also affect the spline approximation error and estimation error. 
The simultaneous diagonlization technique in the scope of penalized spline functional linear regression models is established in Proposition~\ref{prop:splineStruct}, and the corresponding approximation error is quantified in Proposition~\ref{proposition:general:splineerror}.  Coupled with penalized spline penalty, the related empirical norm is  shown to converge under weaker condition in Proposition~\ref{proposition:empiricalnorm}. See Section~\ref{sec:prelim} for more discussions.


%
%
%

\begin{table}[h!]
	\centering
	\caption{Rates of convergence for the reduced model (see \eqref{eqn:method:reduced}) under different configurations of $d$, $\eta_1$, and $K$ as $N\to\infty$. In the table, 
		$\tau = \nu \wedge (\mathfrak{o}+1) + \{ q \wedge (\mathfrak{o}+1)\}/2$, 
		$\iota = q+d$, and $\kappa = \tau+d-\nu$. \label{tbl:reduceRates}}
	\def\arraystretch{1.2}
	\begin{tabular}{|c|c|l|l|l|}
		\hline
		&		$d$  & $\eta_1$ & $K$ & Rate \\
		\hline
		(i) &		\multirow{2}{*}{$d\le \nu$ }	& $\;\lesssim (MN/R)^{-2(\iota\vee \tau)/(2\tau + 1)} $ &  $\; \asymp (MN/R)^{1/(2\tau + 1)}$ &
		$(MN/R)^{-\tau/(2\tau + 1) }$  \\
		\cline{1-1}		\cline{3-5}
		(ii) &				& $\;\asymp (MN/R)^{-2\iota/(2\iota+1)}$ & $ \;\gtrsim (MN/R)^{\iota / [(2\iota + 1)(\iota\wedge\tau)]}$ &
		$	(MN/R)^{-\iota/(2\iota+1)} $  \\
		\hline
		(iii) &		\multirow{2}{*}{$d > \nu$ }& $\; \lesssim (MN/R)^{-2\iota/(2\tau + 1)}$ & $\; \asymp (MN/R)^{1/(2\tau + 1)}$ & $(MN/R)^{-\tau/(2\tau + 1)} $  \\
		\cline{1-1}		\cline{3-5}
		(iv) &		&  $\; \asymp (MN/R)^{-2\iota \kappa/(\kappa+2\iota\tau)}$ & $\; \asymp  (MN/R)^{\iota/(\kappa+2\iota\tau)} $ & 
		$	(MN/R)^{-\iota\tau/(\kappa +2\iota\tau)}$ \\
		\hline
	\end{tabular}
\end{table}

Integrating all the above technical tools, a unified upper bound of the convergence rate for the proposed model is reached under a general class of loss functions and multi-task relationships (in Section~\ref{sec:maintheory}). 
The obtained unified result is then applied to the reduced rank-$R$ model in Section~\ref{sec:reducedrate} and  to the graph regularized model in Section~\ref{sec:graphrate}. 

For the the reduced rank-$R$ model (see Section~\ref{sec:reducedrate}), the rates of convergence are summarized in Table~\ref{tbl:reduceRates} for various parameter settings.  Table~\ref{tbl:reduceRates} answers how penalized splines behave in the classical single-task functional linear regression by plugging in $M=R=1$.  The optimal rate for estimating a single-task slope function is known as $N^{-(q+\nu) /(2q+2\nu+1)}$ \citep{yuan2010reproducing}.  Generally, in Settings (i) and (iii) of Table~\ref{tbl:reduceRates} where the effect of the roughness penalty $\eta_1$ is weak and $K$ is tuned to be optimal, penalized splines in the reduced model behaves like regression splines. 
In Setting (ii) where $\eta_1$ is tuned to be optimal and $K$ is relatively large, the behavior is like smoothing splines. 
With fixed $M$ and $R$, the optimal rate $N^{-(q+\nu) /(2q+2\nu+1)}$ can be obtained in~Setting (ii) with $d=\nu$. 

For the graph regularized model (see Section~\ref{sec:graphrate}), the obtained rates of convergence are summarized in Table~\ref{tbl:graphRates} for $d\le \nu$ and various parameters. Table~\ref{tbl:graphRates} reveals a more interesting phenomenon according to the strength of the graph regularization parameter $\eta_2$.  When the graph regularization is weak (Settings (i) and (ii) in Table~\ref{tbl:graphRates}), the estimator behaves as if we conduct independently estimation for each task.  When the graph regularization is strong and the number of tasks is large enough (Settings (iii) and (v) in Table~\ref{tbl:graphRates}), the rates of convergence can achieve faster than the optimal rate $N^{-(q+\nu) /(2(q+\nu)+1)}$ of the single-task case. In Settings (iv) and (vi), where graph regularization is strong but the number of tasks is small, the estimator can exhibit a much slower rate of convergence because the penalty bias dominates.

\begin{table}[h!]
	\centering 
	\caption{Rates of convergence for the graph regularized model (see \eqref{eqn:method:graph}) under $d\le \nu$ and different configurations of $\eta_1$, $\eta_2$, and $K$. In the table, 	$\tau = \nu \wedge (\mathfrak{o}+1) + \{ q \wedge (\mathfrak{o}+1)\}/2$,
	$\iota = q+d$, $r_1= \frac{\tau}{\tau(2+\mu) + 1}$, and $r_2= \frac{\iota}{\iota(2+\mu) + 1}$. The rows are divided into two groups: weak graph regularization ($\eta_2\lesssim M^{-2/\mu}$) and strong graph regularization ($\eta_2\gtrsim M^{-2/\mu}$). \label{tbl:graphRates}}
	\def\arraystretch{1.65}
	\resizebox{0.98\textwidth}{!}{	\begin{tabular}{|c|c|  l | l | l | l|}
			\hline
 &				Graph Reg. &		$\eta_1$  & $K$ & $\eta_2$  & Rate \\
			\hline
(i) &				\multirow{2}{*}{Weak} &	 $\lesssim N^{-\frac{2(\iota \vee \tau)}{2\tau + 1} }$	
			& $\; \asymp N^{1/(2\tau + 1)}$ &
			$\;\lesssim M^{-\frac{2}{\mu} } \wedge N^{-\frac{2\tau}{2\tau + 1} }$ &  
			$N^{- \frac{\tau}{2\tau + 1}}$  \\
\cline{1-1}			\cline{3-6}
(ii) &				&$\asymp N^{-\frac{2\iota}{2\iota +1} }$	& $ \;\gtrsim N^{\frac{\iota}{(2\iota+1)( \iota \vee \tau)} }$  	& $\;\lesssim M^{-\frac{2}{\mu} } \wedge N^{-\frac{2\iota}{2\iota+1} } $ & 
			$N^{- \frac{ \iota}{2\iota+1} }$  \\
			\hline
(iii) &				\multirow{4}{*}{Strong} &	 \multirow{2}{*}{$\lesssim (MN)^{-\frac{2(\iota\vee\tau)r_1 }{\tau} }$	}
			& \multirow{2}{*}{ $\; \asymp (MN)^{r_1 /\tau}$ } 
			& $\; \asymp (MN)^{-2r_1}$ if 	$M\gtrsim N^{\frac{\mu r_1 }{1-\mu r_1} }$ &
			$(MN)^{-r_1}$  \\
\cline{1-1}			\cline{5-6}		
(iv) &				&			& & $\; \asymp M^{-2/\mu}$ if  $M\ll N^{\frac{\mu r_1}{1-\mu r_1} }$ &  $M^{-\frac{1}{\mu}}$ \\
\cline{1-1}			\cline{3-6}		
(v) &				&
			\multirow{2}{*}{$\asymp (MN)^{-2r_2}$}		& 
			\multirow{2}{*}{$ \;\gtrsim (MN)^{\frac{r_2}{\iota\wedge \tau} }$} &
			$\; \asymp (MN)^{-2r_2} $  if $M\gtrsim N^{\frac{\mu r_2}{1-\mu r_2} }$& 
			$(MN)^{-r_2}$  \\
\cline{1-1}			\cline{5-6}		
(vi) &				& 	 & &  $ \asymp M^{-2/\mu}$ if $M\ll N^{\frac{\mu r_2}{1-\mu r_2} }$ &  $M^{-\frac{1}{\mu}}$ \\
			\hline
	\end{tabular}}
\end{table}

\subsection{Organization and Notations}


The rest of this paper is organized as follows.  Section~\ref{sec:singleregression} reviews single-task functional linear regression with penalized splines. The proposed multi-task functional linear regression with double regularization is presented in Section~\ref{sec:method}. 
For the analysis of penalized splines in the context of functional linear regression, the fundamental tools are established in Section~\ref{sec:prelim}. 
We start to formally examine the proposed multi-task model in Section~\ref{sec:approError}, where the approximation error is defined. 
The estimation error is quantified in Section~\ref{sec:manifold} by controlling the empirical processes over manifold. 
A unified upper bound of the convergence rate is developed for the proposed model and presented in Section~\ref{sec:maintheory}. 
We finally apply the upper bound to the reduced multi-task model and the graph regularized multi-task model in Sections~\ref{sec:reducedrate} and \ref{sec:graphrate}, respectively.
The main conclusions of this paper are summarized in Section \ref{sec:discussion} with some remarks on future work. 
Technical proofs are all provided in Appendix.

Throughout the manuscript, for two sequences of numbers $\{a_n\}$ and $\{b_n\}$, we write $a_n\lesssim b_n$ if $a_n\le C\cdot b_n$ for some positive constant $C$. When $a_n\lesssim b_n$ and $b_n\lesssim a_n$, their relation is denoted as $a_n \asymp b_n$.
We write $a_n \ll b_n$ if $a_n/b_n \to 0$ as $n\to \infty$. 
For two numbers $a,b\in \bbR$, we denote $a\wedge b = \min\{a,b\}$ and $a\vee b  = \max\{a,b\}$. 
Let $\mathbb{L}_2(\sT)$ denote the set of square-integrable functions on domain $\sT$. 
For $f_1,f_2 \in \mathbb{L}_2(\sT)$, their inner product is denoted as $\langle f_1,\, f_2\rangle = \int_{\sT} f_1(t) \cdot f_2(t) \intd t$.
The $L_2$ norm $\|f\|_{L_2}$ is determined as $\|f\|_{L_2}^2 =\int_{\sT} f^2(t)\intd t$.
Table~\ref{tbl:notation} lists the frequently used notations in this work.

\begin{table}[t]
	\centering
	\caption{The frequently-used notations throughout the work. \label{tbl:notation}}
		\def\arraystretch{1.08}
	\resizebox{0.87\textwidth}{!}{
		\begin{tabular}{|r|l|l|}
			\hline
			Notation & Meaning &  \\
			\hline
			$\mB$ & The spline coefficient matrix. & See~\eqref{eqn:lossfun:mat}. \\
			$\bmB_0$ & The optimal parameter   \textit{without} the  constraint  $\mathbb{M}$. & See~\eqref{def:popuestimate}. \\
			$\bmB$ & The optimal parameter  \textit{with} the  constraint  $\mathbb{M}$. &
			See~\eqref{def:popuestimate:constraint}. \\
			$d$ & The penalty derivative order. & See~\eqref{eqn:method:penaltyP}.\\
			$\widehat{\delta}_N$ &  The critical radius. & See~\eqref{eqn:criticalradius}. \\
			$\mGamma$ & The roughness penalty matrix. & See~\eqref{eqn:penaltyGamma}. \\
			$h$ & The kernel bandwidth parameter. &  See~\eqref{eqn:edgeweights}. \\
			$\iota$ & $\iota = q+d$. & \\
			$K$ & The degree of freedom of the spline basis $\vphi(\cdot )$.  & \\
			$\mathbb{M}$ & The constraint matrix manifold for $\mB$. & See~\eqref{eqn:method:mixed}. \\	
			$M$ & The number of tasks. & \\
			$N$ & The number of samples for each task. & \\
			$\mathfrak{o}+1$ & The order of the spline basis $\vphi(\cdot )$. & See Proposition \ref{proposition:general:splineerror}. \\
			$\mOmega$ & The graph Laplacian matrix. &  	See~\eqref{eqn:method:penaltyGraph:p1}. \\
			$q$ & The smoothness of the covariance function. & See Condition~\ref{ass:Kx}.\\		
			$\mathcal{S}$ & The manifold for the auxiliary variables. &  See Section~\ref{sec:graphrate:review}.\\
			$\mathcal{T}$ & The domain of the slope function $\beta(t)$. & \\
			$\tau$ & $\tau = \nu \wedge (\mathfrak{o}+1) +[ q \wedge (\mathfrak{o}+1)]/2$. & \\
			$\nu$ & The smoothness of the true slope function. & See Condition~\ref{ass:beta}. \\
			$\mu$ & The intrinsic dimension of $\sS$. & \\
			$\|\cdot\|_X$, $\|\cdot\|_{\Gamma}$ & Two norms for the slope function $\beta$. &  See~\eqref{eqn:twobasicnorms}.\\
			\hline
		\end{tabular}
	}
\end{table}

\section{Single-task Functional Linear Regression with Penalized Splines}\label{sec:singleregression}

The classical functional linear linear regression \citep{cardot03}  models a single response variable $Y\in\bbR$ and a random functional covariate $X(t)$ on a compact domain $\sT$. 
The regression model predicts $Y$ via the integrated quantity $U=\alpha + \int_{\sT} X(t) \beta(t) \intd t$,
where $\alpha$ and $\beta(\cdot)$ are the intercept and slope function, respectively.
When $Y | X$ follows a distribution of the exponential family,
we can consider the generalized functional linear model 
\begin{equation}\label{eqn:singleGflm}
	g\{E (Y | X)\} = \alpha + \int_{\sT} X(t) \beta(t) \intd t ,
\end{equation}
for some link function $g$. 
For model \eqref{eqn:singleGflm} with the canonical link, the conditional distribution $Y$ given the canonical parameter $U$ takes the form of
$$
    P(Y| U) \propto \exp\Big\{\frac{YU  - \psi(U)}{c(\sigma)}\Big\},
$$
where $\psi^{\prime -1}=g$. In this case, the corresponding loss function is the negative log-likelihood $\ell(y, u) = -yu+\psi(u)$ for the estimation of $\alpha$ and $\beta(\cdot)$.

Generally, the conditional mean of $Y$ is just one way to summarize the conditional distribution of $Y$.
To characterize more aspects of the conditional distribution, we can instead focus on the conditional quantile of $Y$  \citep{cardot2005quantile}
\begin{equation}\label{eqn:singleQuanflm}
	Q_{Y |X}(w) = \alpha + \int_{\sT} X(t) \beta(t) \intd t ,
\end{equation}
where $Q_{Y |X}(w) $ is the $w$-quantile ($w \in (0,1)$) for the conditional distribution $F_{Y|  X}(y)=P(Y \leq y|  X)$, i.e.,
$Q_{Y|X}(w)=F_{Y| X}^{-1}(w)=\inf\big\{ y: F_{Y| X}(y)\geq w\big\}.$
As for the quantile regression model~\eqref{eqn:singleQuanflm}, the loss function is usually chosen as $\ell(y, u) = (y-u)\times \{w - I(y<u) \}$.

Suppose the pair of random elements $(X,Y)$ follows some model, like~\eqref{eqn:singleGflm} or~\eqref{eqn:singleQuanflm}. 
We have $N$ independent realizations $\{(x_n, y_n)\}_{n=1}^N$ of the pair $(X,Y)$, and aim to estimate the true intercept $\alpha_0$ and slope function $\beta_0$ from these samples.
The method of penalized splines approximately represents  $\beta_0$ by a function $\beta$ in a spline space $\mathbb{S}_K$ with $K$ degrees of freedom. 
Let $\vphi = (\phi_1,\cdots, \phi_K)\trans$ denote a vector of B-spline basis functions of $\mathbb{S}_K$, then $\beta(\cdot)= \vphi{\trans}(\cdot) \vb$ where $\vb\in\mathbb{R}^K$ is the spline coefficient vector to be estimated.
Using an appropriate loss function $\ell(\cdot, \cdot)$ and a roughness penalty to avoid overfitting, 
we obtain an estimate $(\widehat{\alpha}, \widehat{\beta})$  by solving
\begin{equation}\label{eqn:singleObj}
(\widehat{\alpha}, \widehat{\beta}) = \argmin_{\alpha\in\bbR, \beta\in\mathbb{S}_K} 
	\sum_{n=1}^{N}\ell \Big(y_{n},\  \alpha + \int_{\sT} x_{n}(t) \beta(t) \intd t\Big) +
	\eta_1 \int_{\sT} \big\{\beta^{(d)}(t)\big\}^2\intd t,
\end{equation}
where $\eta_1$ is the penalty parameter for the roughness penalty and the superscript $d$ represents the $d$-th order of derivative. Given $\beta(\cdot) = \vphi{\trans} (\cdot )\vb$, the roughness penalty in \eqref{eqn:singleObj} has an explicit form in terms of $\vb$, 
\begin{equation} \label{eqn:method:penaltyP:single}
    \sP_{\eta_1}(\vb) := \eta_1 
	\int_{\sT} \big\{\beta^{(d)}(t)  \big\}^2 \intd t
	= \eta_1  \vb{\trans}\mGamma\vb,
\end{equation}
with 
\begin{equation} \label{eqn:penaltyGamma}
	\mGamma = \int_{\sT} \vphi^{(d)}(t) \big\{\vphi^{(d)}(t)\big\}{\trans} \intd t.
\end{equation} 
The essence of the above penalized spline technique is to use a moderately large $K$ to balance computational complexity and approximation bias. Meanwhile, it exploits the penalty $\sP_{\eta_1}(\vb)$ to prevent overfitting.

\section{Multi-task Functional Linear Regression with Double Regularization}
\label{sec:method}

Beyond the single-task regression model in Section~\ref{sec:singleregression}, we are interested in simultaneously estimating the intercepts and slope functions for $M$ regression tasks.
For the $m$-th task, $m=1, \dots, M$, it follows some functional linear regression model, such as models \eqref{eqn:singleGflm} or \eqref{eqn:singleQuanflm}, with unknown true intercept $\alpha_{0m}$ and slope function $\beta_{0m}$. 
For simplicity, we assume the observation numbers are the same for all tasks, i.e., there are $N$ pairs of observations $\{(x_{nm},y_{nm})\}_{n=1}^N$ for each task, and the associated loss function is $\ell_m(y,u)$, $m=1, \dots, M$. 
Based on the samples, the aggregated loss for estimation is
\begin{equation}  \label{eqn:lossfun:L2}
	\sL(\valpha,\vbeta) = \frac{1}{NM}\sum_{m=1}^M\sum_{n=1}^N \ell_m\Big(y_{nm},\, \alpha_m+ \int x_{nm}(t)\beta_m(t) \intd t\Big),
\end{equation}
where $\valpha= (\alpha_1, \cdots, \alpha_{M}){\trans}$ and $\vbeta= (\beta_1(t),\cdots, \beta_M(t))\trans$.
The above loss is viewed as a function of $(\valpha,\vbeta) \in \bbR^M\times [\mathbb{L}_2(\mathcal{T}) ]^M$, where
$\mathbb{L}_2(\mathcal{T})$ is the set of square-integrable functions over the domain $\mathcal{T}$.

Assisted by the penalized spline technique, we represent  $\beta_m(\cdot)$ in the spline space $\mathbb{S}_K$ via $\beta_m(\cdot) = \vphi{\trans} (\cdot )\vb_m$. 
All spline coefficients can be stacked into a matrix $\mB = (\vb_1, \cdots, \vb_M) \in \bbR^{K\times M}$. The aggregated loss function~\eqref{eqn:lossfun:L2} can be written as a function with respect to $\valpha$ and $\mB$: 
\begin{equation}\label{eqn:lossfun:mat}
	\sL(\valpha,\mB) = \frac{1}{NM}\sum_{m=1}^M\sum_{n=1}^N \ell_m(y_{nm},\, \alpha_m+ \vx_{nm}{\trans}\vb_m),
\end{equation}
where $\vx_{nm} = \int_{\sT} x_{nm}(t) \vphi(t) \intd t$
is the vector obtained through the integration of the covariate $x_{nm}$ with the spline basis $\vphi$.
Combining  the roughness penalties for all slope functions in the same form of \eqref{eqn:singleObj},  we get  a penalty in terms of $\mB$, i.e., 
\begin{equation} \label{eqn:method:penaltyP}
	\sP_{\eta_1}(\mB) := \eta_1 \sum_{m=1}^M
	\int_{\sT} \big\{\beta_m^{(d)}(t)  \big\}^2 \intd t
	= \eta_1 \sum_{m=1}^M \vb_m{\trans}\mGamma\vb_m
	= \eta_1 \tr\big(\mB{\trans}\mGamma\mB\big).
\end{equation}
A naive penalized spline estimator for multi-task problem can be obtained via solving
\begin{equation} \label{eqn:method:simpleAggregation}
(\widehat{\valpha}, \widehat{\mB})	=\argmin_{\valpha,\mB}\, \sL(\valpha,\mB) + \sP_{\eta_1} (\mB).
\end{equation}
The estimated slope function for the $m$-th task is
$\widehat{\beta}_m(\cdot) = \vphi{\trans} (\cdot ) \widehat{\vb}_m$
where $\widehat{\vb}_m$ is the $m$-th column of $\widehat{\mB}$.
However, it is evident that the estimator in~\eqref{eqn:method:simpleAggregation} does not enjoy any improvement over the single-task setting, under which \eqref{eqn:singleObj} is applied to each task independently. 

\subsection{The Reduced Multi-task Model}
\label{sec:method:reduced}
One  remedy for improving the estimator in~\eqref{eqn:method:simpleAggregation} from the setting of single-task regression models is to impose low-rank structure among the slope functions. 
Specifically, it is assumed that each slope function can be well approximated by a combination of $R$ representation functions where $R$ is much smaller than $M$. 
Denote the $R$ representation functions by $\psi_{1}, \cdots, \psi_{R}$. Each slope function can be approximated by $\beta_{0m}(\cdot) \approx \sum_{r=1}^R A_{mr} \psi_{r}(\cdot)$ for some coefficients $A_{mr}$, $r=1,\dots, R$. 
Employing spline expansion with basis functions $\vphi = (\phi_1,\cdots, \phi_K)\trans$ and ignoring the approximation errors, we further write $\psi_{r}(\cdot) = \sum_{k=1}^{K} D_{kr} \phi_k(\cdot)$, where $\mD = \big(D_{kr}\big) \in\bbR^{K\times R}$ is an unknown spline coefficient matrix to be estimated. 
This approach induces an approximated low-rank structure since the coefficient matrix $\mB$ satisfies the decomposition $\mB = \mD\mA{\trans}$ with $\mA = \big(A_{mr}\big) \in\bbR^{M\times R}$ and $\mD\in\bbR^{K\times R}$.  
In other words, each slope functions $\beta_{m}$ has the expression $\beta_m(\cdot) = \vphi\trans(\cdot )\vb_m= \sum_{r=1}^{R} A_{mr} \sum_{k=1}^{K} D_{kr} \phi_k(\cdot)$.
This leads us to formulate the reduced (rank) multi-task model
\begin{equation} \label{eqn:method:reduced}
    (\widehat{\valpha}, \widehat{\mB})	=\	\argmin_{\valpha,\rank(\mB) = R}\, \sL(\valpha,\mB) + \sP_{\eta_1} (\mB),
\end{equation}
where the loss function $\sL(\valpha,\mB)$ and penalty term $\sP_{\eta_1} (\mB)$ remain the same as~\eqref{eqn:lossfun:mat} and \eqref{eqn:method:penaltyP}, respectively.
In contrast to~\eqref{eqn:method:simpleAggregation}, an additional constraint $\rank(\mB) = R$ is imposed over $\mB$ for the spline coefficient matrix.
	
\subsection{The Graph Regularized Multi-task Model}\label{sec:method:graph}

In some applications, the relationships between tasks can be determined by some external covariates. 
For example, we may be interested in predicting the average level of air pollutants (as the response) from the wind speed curve (as the functional covariate) at different locations \citep{he2022varying}. 
Each location corresponds to a regression task and we can expect nearby spatial locations have similar slope functions. 
In this case, the spatial coordinates can be treated as external covariate, and two tasks are similar if their spatial coordinates are close to each other. 

In these examples, external covariates provide extra information of measuring the similarity between different tasks. To be specific, suppose the $m$-th task is associated with an external covariate $\vs_m\in\bbR^{s}$, for $m=1, \dots, M$. 
The similarity $w_{vv'}$ between the $v$-th and $v'$-th tasks can be determined by $\vs_v$, $\vs_{v'}$, and a decreasing function $G(\cdot):[0,\infty)\to[0,\infty)$ via
\begin{equation} \label{eqn:edgeweights}
 	w_{vv'}=\frac{2}{\sigma_G h^{\mu+2}M} G\big(-\|\vs_v-\vs_{v'}\|_2 / h\big).
\end{equation}
In the above, $h$ is a bandwidth parameter and $\sigma_G = \int_{\bbR^s} s_1^2 G(\|\vs\|_2)\intd \vs $ with $s_1$ being the first coordinate of $\vs$. 
As in \cite{trillos2020error}, we let $G(\cdot):[0,\infty)\to[0,\infty)$ have support $[0, 1]$ and be Lipschitz continuous. Normalizing $G$ allows us to assume $\int_{\bbR^s} G(\|\vs\|_2)\intd \vs = 1$. 

Given the weights $w_{vv'}$ measuring the similarity between each pair of tasks, we introduce a penalty for the slope functions
$\vbeta = (\beta_1,\cdots, \beta_M)$ as
\begin{align} 
	\sP_{\veta}(\vbeta) &= \eta_1 \sum_{m=1}^M \int_{\sT} \big\{\beta_m^{(d)}(t)  \big\}^2 \intd t 	
	+\eta_2	\sum_{v,v'=1}^M
	\frac{w_{vv'}}{NM} \sum_{n,m}	\langle x_{nm},\beta_v -  \beta_{v'} \rangle^2 \nonumber \\
	& \qquad \qquad \qquad + \eta_{1}\eta_2
	\sum_{v,v'=1}^M w_{vv'}\int \big\{\beta_v^{(d)}(t) -  \beta_{v'}^{(d)}(t) \big\}^2 \intd t. \label{eqn:method:penaltyGraph:fun}
\end{align}
The first term on the right hand side of~\eqref{eqn:method:penaltyGraph:fun} is exactly the roughness penalty~\eqref{eqn:method:penaltyP}. The last two terms on the right hand side of~\eqref{eqn:method:penaltyGraph:fun} encourage between-task similarity of the slope functions.  In particular, the second term encourages the predicted values by similar tasks to be close.  Meanwhile, the third term encourages the adjacent tasks to share similar $d$-th order derivative values of their slope functions.

The similarity weights~\eqref{eqn:edgeweights} induces a weighted graph $\mathcal{G} = (\sV,\sE)$. 
Each element in the vertex set $\sV$ represents a task. There exists an edge $e_{v v'}\in \sE$ connecting the $v$-th and $v'$-th tasks if  $w_{vv'}>0$. We can define a weighted adjacency matrix $\mW = \big(w_{vv'}\big)\in\bbR^{M\times M}$. 
The penalization term~\eqref{eqn:method:penaltyGraph:fun} encodes the intrinsic structure of the graph $\mathcal{G}$ via its graph Laplacian~\citep{chung1997spectral}. We let the degree of the $v$-th vertex be $d_v = \sum_{v'=1}^{M} w_{v v'}$.	
The diagonal matrix with degrees $d_1, \dots, d_M$ in the diagonal is called the degree matrix and denoted by $\mD\in\bbR^{M\times M}$.  The unnormalized graph Laplacian matrix is defined as $\mOmega = \mD - \mW$. 
When each slope function is expressed as $\beta_m(\cdot)=\vphi(\cdot)\trans \vb_m$ in the spline space $\mathbb{S}_K$, the second term on the right hand side of~\eqref{eqn:method:penaltyGraph:fun} has equivalent expression
\begin{equation}  \label{eqn:method:penaltyGraph:p1}
	\sum_{v,v'=1}^M
	\frac{w_{vv'}}{NM} \sum_{n,m}
	\langle x_{nm},\beta_v -  \beta_{v'} \rangle^2  
	= \sum_{v,v'=1}^M w_{vv'} 
	\|\hat{\mSigma}^{1/2}(\vb_v - \vb_{v'} ) \|^2
	=  \tr\big(\mB\mOmega \mB\trans \widehat{\mSigma}\big),
\end{equation}
where $\widehat{\mSigma} = \frac{1}{MN} \sum_{n=1}^{N} \sum_{m=1}^{M} \vx_{nm}\vx_{nm}\trans$ is the pooled covariance matrix. The last term in~\eqref{eqn:method:penaltyGraph:fun} then becomes
\begin{equation} \label{eqn:method:penaltyGraph:p2}
	\sum_{v,v'=1}^M
	w_{vv'}\int \big\{\beta_v^{(d)}(t) -  \beta_{v'}^{(d)}(t) \big\}^2 \intd t 
	=  \sum_{v,v'=1}^M w_{vv'} \|\mGamma^{1/2}(\vb_v - \vb_{v'} )\|^2
	=  \tr\big(\mB\trans\mGamma\mB\mOmega \big).
\end{equation}
Therefore, the penalization term~\eqref{eqn:method:penaltyGraph:fun} can be rewritten as
\begin{equation} \label{eqn:method:penaltyGraph}
	\sP_{\veta}(\mB) = \eta_1 \tr\big(\mB\trans\mGamma\mB\big) +
	\eta_2 \tr\big(\mB\mOmega\mB\trans\widehat{\mSigma}\big) + \eta_1\eta_2
	\tr\big(\mB\trans\mGamma\mB\mOmega\big),
\end{equation}
where $\veta = (\eta_1,\eta_2)$ is the set of penalty parameters. In summary, we have the following objective function for graph regularized multi-task learning
\begin{equation} \label{eqn:method:graph}
(\widehat{\valpha}, \widehat{\mB})	=	\argmin_{\valpha,\mB\in\bbR^{K\times M}}\, \sL(\valpha,\mB) + \sP_{\veta} (\mB),
\end{equation}
where the loss function is the same as~\eqref{eqn:lossfun:mat} but \eqref{eqn:method:penaltyGraph} is employed in the penalization term. 

\subsection{The General Model with Double Regularization}\label{sec:method:general}
The above two multi-task models~\eqref{eqn:method:reduced} and~\eqref{eqn:method:graph} can be unified through a general model. 
With the loss function~\eqref{eqn:lossfun:mat}, we propose the following penalized estimator
\begin{equation}\label{eqn:method:mixed}
    (\widehat{\valpha}, \widehat{\mB})	=	\argmin_{\valpha, \mB \in \mathbb{M}}\, \sL(\valpha, \mB) + \sP_{\veta} (\mB),
\end{equation}
with double regularization on the spline coefficient matrix $\mB$.
The first regularization over $\mB$ is the constraint set $\mathbb{M} \subset \bbR^{K\times M}$. We consider the setting where $\mathbb{M}$ is a Riemannian embedded submanifold of $\bbR^{K\times M}$ without boundary. 
The second regularization over $\mB$ is the penalty $\sP_{\veta}(\cdot)$. 
It is a general composite quadratic penalty with  parameter(s) $\veta = (\eta_1,\cdots, \eta_P)$ and takes the form of
\begin{equation} \label{eqn:method:generalPenalty}
	\sP_{\veta} (\mB) = \sum_{j=1}^{P} \eta_j\tr\big(\mB\trans\mPi_{j1}\mB\mPi_{j2}\big),
\end{equation}
where $\mPi_{j1}$'s and $\mPi_{j2}$'s are symmetric positive semi-definite matrices. 
Because the penalized splines are employed for function estimation in this work, the first term in the summation of~\eqref{eqn:method:generalPenalty} is assumed to be the roughness penalty~\eqref{eqn:method:penaltyP}, i.e., $\mPi_{11} = \mGamma$ and $\mPi_{12} = \mI$. 

It is evident the model~\eqref{eqn:method:mixed} includes~\eqref{eqn:method:reduced} and~\eqref{eqn:method:graph} as special cases. 
The model~\eqref{eqn:method:mixed} becomes the reduced (rank) model~\eqref{eqn:method:reduced} when $\mathbb{M}$ is the rank-$R$ matrix manifold, i.e., $\{\mB\in\bbR^{K\times M}:\ \rank(\mB) = R\}$, and the penalty~\eqref{eqn:method:generalPenalty} is specified as $P =1$ with $\mPi_{11} = \mGamma$ and $\mPi_{12} = \mI$.
On the other hand, model \eqref{eqn:method:graph} corresponds to the case where $\mathbb{M} =\bbR^{K\times M}$, consisting of all matrices of size $K\times M$, and the penalty function \eqref{eqn:method:penaltyGraph} satisfies $P=3$ with $\mPi_{11} = \mGamma$, $\mPi_{12} = \mI$, $\mPi_{21} = \mOmega$, $\mPi_{22} = \widehat{\mSigma}$, $\mPi_{31} = \mGamma$, and $\mPi_{32} = \mOmega$, respectively.

In the following sections, we will first develop a unified upper bound of the convergence rate for the general multi-task functional linear regression model~\eqref{eqn:method:mixed}. 
The general result will then be applied to two special structures: the reduced model \eqref{eqn:method:reduced} and the graph regularized model \eqref{eqn:method:penaltyGraph}.
Note when each component of $\vbeta(\cdot) =(\beta_{1}(\cdot) ,\cdots, \beta_{M}(\cdot) )\trans$   is expressed by  splines with $\beta_{m}(\cdot) =\vphi(\cdot)\trans \vb_{m}$, we set $\vb_{m}$ in the $m$-th column of $\mB$, and thus $\vbeta(\cdot)$ and $\mB$ present the same object in essence. To simplify the presentation, we will also use $\sP_{\veta} (\vbeta)\equiv \sP_{\veta}(\mB)$ for the   composite quadratic penalty~\eqref{eqn:method:generalPenalty} in the following.

\section{Preliminaries on the Penalized Spline Technique} \label{sec:prelim}

In the doubly regularized multi-task model~\eqref{eqn:method:mixed}, we use the penalized splines to estimate the slope functions. 
This section develops the technical tools for analyzing penalized spline in the context of functional linear regression. Our results extend those for non-parametric regression in \cite{huang2021asymptotic}. 
The resulted tool will further facilitate to derive the upper bound of the convergence rate of~\eqref{eqn:method:mixed}.
In particular, Section~\ref{sec:prelim:error} provides the approximation error of true slope function using the spline space. Simultaneously diagonalization is also developed for two quadratic forms based on the covariance function and penalty. 
In Section~\ref{sec:prelim:ellipsoid}, we argue that the estimation error of the penalized spline estimator is connected to the complexity of an ellipsoid formed by the two quadratic forms. More precisely, the phase-transition behavior of penalized spline either like a regression spline estimator or like a smoothing spline estimator is determined by the complexity of the ellipsoid. 
After that, in Section~\ref{sec:prelim:empirical}, a novel result on the convergence in terms of empirical norm follows. To our best knowledge, the technical results in this section serve as novel contributions to the literature of penalized spline estimator and provide insights into our unified multi-task model~\eqref{eqn:method:mixed}.

\subsection{Spline Approximation and Simultaneous Diagonalization}\label{sec:prelim:error}
As our model resides in the spline space, we begin with investigating the approximation power of the spline space to the true slope function $\beta_{0m}$ in terms of prediction error. 
For this purpose, some regularity assumptions are required on the true slope function and the covariance function of the functional predictor. The first condition assumes $\beta_{0m}$ is smooth and belongs to the Sobolev space of order $\nu$.

\begin{condition}\label{ass:beta}
	The true slope function $\beta_{0m}$ belongs to the Sobolev space of order $\nu$, i.e.,
	$ \beta_{0m} \in \mathbb{L}_2^{\nu}(\sT) := \{\beta:\; \beta^{(k)} \in \mathbb{L}_2(\sT) \text{ for } k\le \nu\}\,$, where $\beta^{(k)}$ represent the weak derivative of $\beta$ of order $k$.
\end{condition}

Suppose each functional predictor $x_{nm}$ has zero mean, and covariance function $\mathcal{C}_m(t,t') = \Expect \{x_{nm}(t) x_{nm}(t')\}$ for the $m$-th task. 
For succinct presentation, we assume the covariance functions are the same across different tasks, i.e., $\mathcal{C}_m\equiv \mathcal{C}$. 
It is important to note that the conclusions derived in this work can be generalized to a general setting with diverse covariance functions for various tasks. 
See Remark~\ref{remark:disinctC} at the end of this section for the detailed discussion on this general setting with a proof outline.

\begin{condition}\label{ass:Kx}
	The covariance function $\mathcal{C}$ satisfies the following properties  for some positive integer $q$ and non-negative integer $p$:
	\begin{enumerate}
        \item[(i)] Denote $\{\lambda_{0j} \}$ as the non-increasing sequence of the eigenvalues of $\mathcal{C}$. 
        The eigenvalues decay with the order $\lambda_{0j}\asymp j^{-2q}$.
        \item[(ii)] Denote $\mathcal{C}^{(k,l)}(t, t') = \frac{\partial^{k+l}}{\partial u^k \partial v^l  } \mathcal{C}(t, t')$. 
        For $i,j=0,\cdots, q-1$, the (weak) derivatives  $\mathcal{C}^{(i,j)}$, $\mathcal{C}^{(q,q-1)}$, and $\mathcal{C}^{(q-1,q)}$ exist and are square integrable.
        \item[(iii)] $\mathbb{L}_2^q(\sT) = H(\mathcal{C}) \oplus \mathbb{P}_{p}$, where $H(\mathcal{C})$ is the reproducing kernel Hilbert space with kernel $\mathcal{C}$, $\mathbb{P}_{p}$ is a   subspace with  dimension $p$, and $\oplus$ represents direct sum of subspace.
	\end{enumerate}
\end{condition}

Condition~\ref{ass:Kx} is related to the smoothness requirement of the covariance function $\mathcal{C}$.  In particular, Point (iii) of Condition~\ref{ass:Kx} means the support of the probability measure of the random $x_{nm}$ can be a proper subspace of the full Sobolev space $\mathbb{L}_2^q(\sT)$, such that there is a null space of dimension $p$ over which the random $x_{mn}$ has no variability. Note the value of $p$ can be flexible. Its value can be $0$ (i.e., $\mathbb{P}_{p}$ is an empty set) or a large number depending on the particular functional data of interest. 

Condition~\ref{ass:Kx} is mild and one of its sufficient conditions is the \textit{Sacks-Ylvisaker condition} \citep{ritter1995}. The Sacks-Ylvisaker condition was discussed in the literature of functional linear regression with smoothing splines \citep{yuan2010reproducing,du2014penalized} to justify the eigenvalue decay and the sample path smoothness of $x_{nm}$.  
In this work, we use this simplified version of the Sacks-Ylvisaker condition, because Condition~\ref{ass:Kx} highlights the essential properties of the covariance function for the analysis of the penalized spline estimator. 
These essential properties include the eigenvalue decay rate, the smoothness of the covariance function, and the possible existence of the null space $\mathbb{P}_p$. 
It can be seen that the covariance functions of many stochastic processes satisfy  Condition~\ref{ass:Kx}. As an example, the Brownian motion covariance function  $\mathcal{C}(s,t) = \min (s,t)$ satisfies Condition~\ref{ass:Kx} with $q=1$ and $p=1$, and $\mathbb{P}_p$ is the space of constant functions. In this case, the Brownian motion has no variability in the subspace $\mathbb{P}_p$. More examples of covariance functions satisfying Condition~\ref{ass:Kx} can be found in Appendix~\ref{sec:covexamples}.

The estimation performance for functional linear regression is intimately connected to the covariance function $\mathcal{C}$ of the predictors. In this work, of particular importance is the decay rate of the eigenvalues of the covariance function. A covariance with faster eigenvalue decay rate will lead to a faster rate of convergence. 
On the other hand, the eigenfunctions of the covariance function and the subspace $\mathbb{P}_p$ do \textit{not} play a significant role in our analysis and will not affect the rates of convergence.

The spline approximation error is measured based on two (semi-)norms, which will play a fundamental role throughout this work.
For a sufficiently smooth $\beta\in \mathbb{L}_2(\sT)$, we define two (semi-)norms $\|\cdot\|_{X} $ and $\|\cdot\|_{\Gamma} $ as follows
\begin{equation} \label{eqn:twobasicnorms}
	\|\beta\|_{X} = 
	\big( \Expect 	\langle x_{nm}, \beta\rangle^2\big)^{1/2}\quad  \text{ and }\quad
	\|\beta\|_\Gamma =  \Big[ \int
	\big\{\beta^{(d)}(t)\big\}^2\intd t \Big]^{1/2}. 
\end{equation}
Note that $\|\cdot\|_{X}$ is the same for all tasks as we have assumed their functional predictors $x_{nm}$'s share a common covariance function. Meanwhile, $\|\cdot\|_\Gamma$ is a semi-norm related to the roughness penalty. The next proposition characterizes the spline approximation error together with the penalty term.

\begin{proposition}\label{proposition:general:splineerror}
Under Conditions~\ref{ass:beta} and~\ref{ass:Kx}, the spline approximation satisfies
\begin{equation} \label{eqn:prop:approError}
	\inf_{\beta\in\mathbb{S}_{K}} 
	\big\{ \| \beta - \beta_{0m}\|_{X}  +\eta_1^{1/2} \| \beta\|_{\Gamma}\big\}
	\asymp  K^{-\tau} + \eta^{1/2}_1 K^{(d-\nu)_{+}},
\end{equation}
with $\tau = \nu \wedge (\mathfrak{o}+1) + \{ q \wedge (\mathfrak{o}+1)\}/2$, where $\mathfrak{o}+1$ is the order of the spline basis. 
\end{proposition}

On the left hand side of~\eqref{eqn:prop:approError}, $\| \beta - \beta_{0m}\|_{X}$ can be interpreted as the expected prediction error when we use a spline approximation $\beta\in\mathbb{S}_{K}$ in place of the true slope function $\beta_{0m}$. 
The second term $\eta_1^{1/2} \| \beta\|_{\Gamma}$ is the amount of incurred penalty for $\beta$ with penalized spline estimation.
The right hand side of~\eqref{eqn:prop:approError} states the approximation error and the penalty term is of order $K^{-\tau} $ and $ \eta^{1/2}_1 K^{(d-\nu)_{+}}$, respectively. 
The penalty order $d$ for the norm $\|\cdot\|_{\Gamma}$ is allowed to be larger than the actual smoothness order $\nu$ of $\beta_{0m}$. When the penalty order $d$ is strictly larger than the actual smoothness order $\nu$ (i.e. $d>v$), increasing the knot number $K$ will increase the penalty bias as well.

Proposition~\ref{proposition:general:splineerror} parallels Theorem~3.1 of~\cite{huang2021asymptotic}, but Proposition~\ref{proposition:general:splineerror} is established in the setting of functional linear regression. In \cite{huang2021asymptotic}, the approximation error of $ f \in\mathbb{S}_{K}$ to a regression function $f_{0}$ is measured in the $L_2$ sense, and they concluded $\| f - f_{0}\|_{L_2} \asymp K^{-\nu\wedge (\mathfrak{o}+1)}$. 
On the other hand, the order of approximation error $K^{-\tau}$ in  Proposition~\ref{proposition:general:splineerror} is smaller than $K^{-\nu\wedge (\mathfrak{o}+1)}$. This is because we have employed a different (semi-)norm $\|\cdot\|_X$ and taken the smoothness of the covariance function into account.  

A key technique of this work is to simultaneously diagonalize the two norms $\|\cdot\|_{X}$ and $\|\cdot\|_{\Gamma}$ defined above.
Simultaneous diagonalization facilitates establishing convergence rate for both smoothing splines \citep{wahba90,gu2013smoothing,yuan2010reproducing} and penalized splines \citep{claeskens2009asymptotic,huang2021asymptotic}. 
Suppose $\widetilde{\vphi}(\cdot)$ is an arbitrary vector of basis functions in $\mathbb{S}_{K}$ (such as the normalized B-spline basis in Section~4.3 of~\citealp{schumaker2007spline}).
For $\beta(\cdot) = \widetilde{\vphi}\trans(\cdot)\widetilde{\vb}\in \mathbb{S}_{K}$ with some spline coefficient vector $\widetilde{\vb}$, it is not difficult to find  
$$\|\beta\|_{X}^2  =
	\widetilde{\vb}\trans \widetilde{\mSigma} \widetilde{\vb} 
	\quad	\text{ with } \quad
	\widetilde{\mSigma} = \mathrm{Var}\Big\{
	\int x_{nm}(t) \widetilde{\vphi}(t) \intd t \Big\}, $$ 
and	
$$\|\beta\|_\Gamma^2  =
	\widetilde{\vb}\trans \widetilde{\mGamma} \widetilde{\vb}
	\quad	\text{ with } \quad
	\widetilde{\mGamma} = \int \widetilde{\vphi}^{(d)}(u) 
	\big\{\widetilde{\vphi}^{(d)}(u)\big\}\trans \intd u.$$
In words, the squares of the two (semi-)norms are simply quadratic forms of the spline coefficient vector $\widetilde{\vb}$. 
In the following, we construct another basis $\vphi(\cdot)$ from the original $\widetilde{\vphi}$ in an appropriate way, such that the $\widetilde{\mSigma}$ and $\widetilde{\mGamma}$ simultaneously become diagonal matrices.
	
\begin{proposition} \label{prop:splineStruct}
Under Condition~\ref{ass:Kx}, there exists an invertible matrix $\mQ$ with which we can define $\vb = \mQ^{-1} \widetilde{\vb}$ and $\vphi(\cdot) = \mQ\widetilde{\vphi}(\cdot)$. 
It follows $\beta(\cdot) = \widetilde{\vphi}\trans(\cdot)\widetilde{\vb} = \vphi\trans(\cdot)\vb$. 
Further, it holds for some $\bar{p} (\le p)$ that
\begin{equation} \label{prop:splineStruct:resultEqn}
    \|\beta \|_{X}^2 = \widetilde{\vb}\trans\widetilde{\mSigma}\widetilde{\vb} = \vb\trans (\mI_{K-\bar{p} }\oplus \vzero_{\bar{p}})\vb, \qquad
    \|\beta \|_\Gamma^2 = \widetilde{\vb}\trans\widetilde{\mGamma}\widetilde{\vb} = \vb\trans\mGamma\vb,
\end{equation}
where  $\mI_{K-\bar{p}}$ is the identity matrix of size $K-\bar{p}$, and $\vzero_{\bar{p}}$ is a square matrix of size $\bar{p} \times \bar{p}$ filled up with zeros.  Besides, $\mGamma = \mathrm{diag}(\gamma_1,\gamma_2,\cdots, \gamma_K)$ is a diagonal matrix whose diagonal elements $\gamma_k$'s are non-negative and  monotone increasing. They satisfy $\gamma_k \gtrsim k^{(2q + 2d)}$ for $k>2d$ and $\gamma_k\ge 0$ for $k\le 2d$.
\end{proposition}

Recall that Condition~\ref{ass:Kx} assumes the support of the probability measure of the random $x_{nm}$ may not be the full Sobolev space $\mathbb{L}_2^q(\sT)$, but up to an additional finite-dimensional null subspace $\mathbb{P}_p$. 
This implies the (semi-)norm $\|\beta\|_{X}$ could possibly be zero for a non-zero $\beta\in\mathbb{S}_K$ in the spline space. Equivalently, the quadratic term $\|\beta\|_{X}^2$ can have finite zero eigenvalues with respect to $\|\beta\|_{\mathbb{L}_2}^2$ for $\beta \in\mathbb{S}_K$. In Proposition~\ref{prop:splineStruct}, $\bar{p}$ represents the replicate number of the zero eigenvalues. On the other hand, the finite-dimensional subspace $\mathbb{P}_p$ in Condition~\ref{ass:Kx} will not have influence on the prediction error, and therefore will not affect the upper bound of the convergence rate. 
As a consequence, this null space is usually directly ignored in the literature \citep[e.g.,][]{yuan2010reproducing}.
Following the same strategy, we simply set $p=0$ in~(iii) of Condition~\ref{ass:Kx} for presentation convenience (i.e., $\bar{p}=0$) in the rest of this work. Meanwhile,  we will assume the employed spline basis $\vphi$ has already been constructed as in Proposition~\ref{prop:splineStruct}, such that both (semi-)norms  $\|\cdot \|_{X}^2$ and $\|\cdot \|_{\Gamma}^2$ have been diagonalized. 

\subsection{Ellipsoid and the Transition Behavior of Penalized Spline}\label{sec:prelim:ellipsoid}
Suppose $p=0$ and the (semi-)norms $\|\cdot\|_X$ and $\|\cdot\|_\Gamma$ have been diagonalized as discussed at the end of the previous subsection. The two (semi-)norms together with the penalty parameter $\eta_1$ determine an ellipsoid $E$ for the spline coefficient vector in $\bbR^K$, where
\begin{align} 	
	E &= \big\{ \vb\in\bbR^K:\; \beta
	(\cdot ) = \vphi\trans(\cdot) \vb
    \text{ and } \|\beta\|_X^2+\eta_1\|\beta\|_{\Gamma}^2\le 1\big\}  \nonumber \\
    &= \Big\{ \vb\in\bbR^K:\;	\sum_{k=1}^{K} (1+\eta_1 \gamma_{k}) b_k^2 \le 1 \Big\}. \label{eqn:complexity0:ellipsoid}
\end{align}
The half lengths of its principal axes are $1/\sqrt{1+\eta_1 \gamma_{k}}$ for $k=1,\cdots, K$.  In our analysis, we find the complexity of $E$ plays a crucial role in determining the transition behavior of penalized splines, i.e., either like smoothing splines or regression splines. It is evident the complexity of $E$ has intricate dependence over $K$ and $\eta_1$.  Intuitively, when $K$ diverges to infinity slowly and $\eta_1$ decreases to zero fast, the complexity of $E$ is close to that of standard Euclidean $(K-1)$-sphere $\mathbb{S}^{K-1}= \big\{(b_k)_{k=1}^K:\, \sum_{k=1}^K  b_k^2 \le 1\big\}$. In this case, the estimation error of the penalized spline estimator is close to that of the regression spline estimator.  On the other hand, when $K$ diverges to infinity fast and $\eta_1$ decreases to zero slowly, the complexity of $E$ is close to that of infinite dimension Sobolev ellipsoid $\mathcal{E}_\infty = \big\{(b_k)_{k=1}^\infty:\, \sum_{k=1}^\infty k^{2(d+q)} b_k^2 \le 1\big\}$. The corresponding estimation error of the penalized spline estimator approaches that of the smoothing spline estimator. Based on this intuition, the goal is to find the breakpoint of $K$ and $\eta_1$ at which this complexity transition occurs for the above ellipsoid $E$.

In this work, the complexity measurement is assisted by the generic chaining technique \citep{talagrand2014upper}.  The generic chaining technique can provide sharper upper and lower bounds compared to the classical Dudley's integral entropy bounds in certain cases. For instance, as discussed in Section 2.5 of \cite{talagrand2014upper}, Dudley's integral entropy bound may fail to accurately describe the behavior of the empirical process over an ellipsoid. The generic chaining, however, can be both accurate and easy to calculate for characterizing the complexity of an ellipsoid $E$. The tuition discussed in this section will be further extended to the unified model~\eqref{eqn:method:mixed}, where a general ellipsoid-like neighbor is studied. See Section~\ref{sec:manifold:epbound} for more discussions.
	

Given a set $T$ and a metric $d(\cdot,\cdot)$ defined on it, the generic chaining characterizes the complexity of $T$ via the
\textit{$\gamma_\alpha$-functional}
\begin{equation} \label{eqn:gammafunctional}
	\gamma_{\alpha}(T, d(\cdot,\cdot)) = \inf_{\{T_n\}}\sup_{t \in T} \sum_{n=0}^{\infty} 2^{n/\alpha}d(t, T_n),
\end{equation}
where $\alpha\ge 0$ and $\{T_n\}_{n \ge 0}$ is a sequence of subsets of $T$. The subset sequence $\{T_n\}_{n \ge 0}$ should be  \textit{admissible}, which means the cardinality of each $T_n$ is limited by $| T_0| = 1$ and $|T_n| \le 2^{2^n}$.

As a direct consequence of Equation (2.115) and Theorem 4.1.11 of \citep{talagrand2014upper}, the $\gamma_2$-functional of the ellipsoid $E$ is related to the summation of the squared half lengths of its principal axes
$$
    \gamma_{2}(E,d) \asymp \Big( \sum_{k=1}^{K} \frac{1}{1+\eta_1 \gamma_{k}} \Big)^{1/2}.
$$
Because $1+\eta_1 \gamma_k \ge 1$, it is obvious that
$\sum_{k=1}^K {1}/{(1+ \eta_1  \gamma_k)} \le K$; meanwhile, based on Proposition~\ref{prop:splineStruct}, it holds that
$$
    \sum_{k=1}^K \frac{1}{1+ \eta_1  \gamma_k} 
    \lesssim \int_{0}^\infty \frac{1}{1+ \eta_1  y^{2q+2d}}  \intd y \asymp \eta_1^{-\frac{1}{2q+2d}}.
$$
In summary, as $\eta_1\to 0$ and $K\to \infty$, the $\gamma_2$ functional of the ellipsoid $E$ is
\begin{equation}\label{eqn:complexity0:keybound}
	\gamma_{2}(E,d)  \asymp (K^{\frac{1}{2}}\wedge \eta_1^{-\frac{1}{4(q+d)}}).
\end{equation}
Result~\eqref{eqn:complexity0:keybound} indicates the  complexity of the ellipsoid $E$ is determined by the relative magnitude of $K^{1/2}$ and $\eta_1^{-1/(4q+4d)}$. 
Such magnitude comparison underlies the  penalized spline analysis in Theorem~3.2 of~\cite{huang2021asymptotic}.  We can similarly expect that for functional linear regression: the penalized spline estimator imitates regression spline estimator when $K^{1/2} \gtrsim \eta_1^{-1/(4q+4d)}$; otherwise, it will behave like the smoothing spline estimator. 

\begin{remark}
The presented results have several distinctions from those in~\cite{huang2021asymptotic}.  As a study of non-parametric regression, \cite{huang2021asymptotic} considers two norms $\|\cdot\|_{L_2}$ and $\|\cdot\|_\Gamma$ (in our notations) instead of  $\|\cdot\|_X$ and $\|\cdot\|_\Gamma$.  In the context of non-parametric regression,  the convergence rates of the penalized spline  estimator are in fact  determined by a different ellipsoid as $\big\{\vb\in\bbR^K:\; f(\cdot ) = \vphi\trans(\cdot) \vb \text{ and } \; \|f\|_{L_2}^2+\eta_1\|f\|_{\Gamma}^2\le 1\big\}$.  After simultaneous diagonalization of $\|\cdot\|_{L_2}$ and $\|\cdot\|_\Gamma$, the diagonal elements of the matrix $\mGamma$ scale as the rate of $\gamma_k \asymp k^{2d}$, instead of $k^{2(d+q)}$ in our context. 
\end{remark}

\begin{remark}
	To extend the above analysis tool to the general   model~\eqref{eqn:method:mixed}, we note that the roughness penalty is the first summand of the general  penalty~\eqref{eqn:method:generalPenalty}. Inspired by~\eqref{eqn:complexity0:ellipsoid}, we can further consider a set of the form
$$\Big\{\mB=(\vb_1,\cdots, \vb_M): \;\; \beta_m(\cdot) = \vphi\trans(\cdot)\vb_m\text{ and }
\sum_{m=1}^{M}\|\beta_m\|_X^2 + \sP_{\veta} (\vbeta)\le 1 \Big\}$$
for the unified model~\eqref{eqn:method:mixed} with the composite quadratic penalty $\sP_{\veta}(\cdot) $ in~\eqref{eqn:method:generalPenalty}. 
At the same time, the constraint structure of matrix manifold $\mathbb{M}$ should be taken into account. 
The details of these two aspects will be clearly presented in Sections~\ref{sec:manifold} and~\ref{sec:maintheory}.
\end{remark}

\subsection{Convergence of Empirical Norm}\label{sec:prelim:empirical}

As a direct application of the above result~\eqref{eqn:complexity0:keybound}, we develop a convergence result of the empirical norm $\|\beta\|_{Nm}^2: = ({1}/{N}) \sum_{n=1}^{N} \langle x_{nm}, \beta\rangle^2$ to its expected counterpart $\|\beta\|_{X}^2$ in the setting of penalized splines. 
For this purpose, we further assume the random covariate $x_{nm}$ is sub-Gaussian as follows. 

\begin{condition} \label{ass:subGaussian}
	There exists some positive constant $C_g$, such that for any $\beta\in \mathbb{S}_K$, the inner product $\langle\beta, x_{nm}\rangle$ is sub-Gaussian with $\|\langle\beta, x_{nm}\rangle\|_{\psi_2} \le C_g \| \beta \|_{X}$.
\end{condition}

The  convergence of the empirical norm  $\|\cdot \|_{Nm}$ to $\|\cdot \|_X$  requires  $K/N\to 0$, as $K$ and $N$ diverge to infinity. This requirement can be interpreted from the perspective of random matrix theory. 
With  $\vphi$ constructed in Proposition~\ref{prop:splineStruct}, the convergence of $\|\beta\|_{Nm}$ to $\|\beta\|_{X}$ for  any $\beta(\cdot) = \vphi\trans(\cdot)\vb$ is equivalent to the convergence of $\widehat{\mSigma}_m = (1/N) \sum_{n=1}^{N} \vx_{nm}\vx_{nm}\trans$ to  $\mI_K$. This convergence in operator norm entails $K/N\to 0$. See also \citep[][]{huang1998projection} for empirical norm convergence in the context of  a non-parametric function  fitting.

On the other hand, for penalized spline models, two norms $\|\cdot\|_{X}$ and 
$\|\cdot\|_\Gamma$ frequently appear together. 
Considering the  summation of the form
$\|\beta\|_{Nm}^2  +\eta_1 \| \beta\|_{\Gamma}^2$ (or $\|\beta\|_{X}^2  +\eta_1 \| \beta\|_{\Gamma}^2$),
we can turn the condition $K/N\to 0$ into a milder one. 	Essentially, the convergence of $\widehat{\mSigma}_m +\eta_1\mGamma$	to the matrix $\mI+\eta_1\mGamma$ allows a larger budget for relative error, because the diagonal elements of $\mGamma$ increase to infinity as $K\to\infty$. 

\begin{proposition} \label{proposition:empiricalnorm}
	Under Conditions~\ref{ass:Kx} and \ref{ass:subGaussian}, suppose  $\{K\wedge \eta_1^{-1/ (2q+2d) } \}/ N\to 0$ as $N\to\infty$. Then, with probability  at least  $1-\exp\big\{ - K\wedge \eta_1^{-1/ (2q+2d) } \big\}$,	it holds that
	\begin{equation} \label{eqn:event:empLower}
		(1-\epsilon)\big\{\|\beta\|_{X}^2   +\eta_1 \| \beta\|_{\Gamma}^2\big\}
		\le \|\beta\|_{Nm}^2 + \eta_1 \| \beta\|^2_{\Gamma}  \le (1+\epsilon)\big\{\|\beta\|_{X}^2   +\eta_1 \| \beta\|_{\Gamma}^2\big\},
	\end{equation}
		for all $\beta\in\mathbb{S}_K$ in the spline space, and for $\epsilon = C_e \big\{K\wedge \eta_1^{-1/ (2q+2d) } / N\big\}^{1/2}$ with some constant $C_e>0$.
\end{proposition}

In the setting of penalized splines, as $N\to \infty$, we usually have $K\to \infty$ to reduce the spline approximation bias and $\eta_1\to 0$ to reduce the penalty bias. 
The above proposition states that, under the weaker condition $\{K\wedge \eta_1^{-1/ (2q+2d) } \}/ N\to 0$, the empirical norm of a function coupled with the corresponding roughness penalty converges to its expected counterpart in terms of relative error. 
This result is valuable for our theoretical analysis, because it allows $K$ to grows faster than $N$ as long as the penalty parameter $\eta_1$ does not decrease to $0$ too fast.

\begin{remark} \label{remark:disinctC}
We assume the covariance functions are the same across various tasks, i.e., $\mathcal{C}_m\equiv \mathcal{C}$. When these functions differ, the developed theoretical tools remain applicable.  To see this,  we first note that eigenfunctions do \textit{not} play any significant role in our analysis of convergence rates.  It is thus absolutely fine for the eigenfunctions to be  different among multiple tasks.

As for the eigenvalues, we now consider the case that, for $m=1,\ldots, M$, each covariance function $\mathcal{C}_m$ satisfies Condition~\ref{ass:Kx} with task-specific parameters $(p_m, q_m)$, instead of the common constants $(p,q)$ as in the current work. In other words, the eigenvalue decay rates $q_m$'s (recall that $\lambda_{m,j}\asymp j^{-2q_m}$ for the $j$-th eigenvalue of $\mathcal{C}_m$ according to Condition~\ref{ass:Kx}) are distinct across various tasks.  In this case, instead of a common norm $\|\beta\|_{X}$  shared by all tasks and studied in~\eqref{eqn:twobasicnorms}, we define $\|\beta\|_{X,m} = \big( \Expect \langle x_{nm}, \beta\rangle^2\big)^{1/2} $	for each task ($m=1,\ldots, M$).  Then, Propositions~\ref{proposition:general:splineerror}, \ref{prop:splineStruct}, and \ref{proposition:empiricalnorm} can be applied to each $\|\beta\|_{X,m} $ separately.  Afterward, \eqref{eqn:complexity0:keybound} suggests a complexity measurement $K^{1/2}\wedge \eta_1^{-1/\{4({q_m}+d)\}}$ for the local neighborhood of each task. It further implies that the phase transition behaviour would be different for each task because the values of their corresponding $q_m$'s are different.  Moreover, the tuning parameter associated with the penalized spline (i.e., $\eta_1$ in the current manuscript) needs to be assigned distinct values for various tasks to recover the optimal rate of convergence. 
\end{remark}

\section{Optimal Model Parameter and Approximation Error}
\label{sec:approError}

We now start to address the theoretical properties of the doubly regularized estimator~\eqref{eqn:method:mixed}. 
Developing the upper bound of the convergence rate of the slope functions $\beta_m$ is of primary interest of this work. 
For simplicity, we assume the intercepts $\alpha_m$'s are zero and focus on analyzing the estimator of $\beta_m$'s. 
Taking $\alpha_m$ into consideration will not affect the rate of convergence but only make the technical proofs more complicated. 
Given $\valpha=\vzero$, we write $\sL(\vbeta)\equiv \sL(\valpha,\vbeta)$ for the loss function with respect to slope functions in~\eqref{eqn:lossfun:L2}
and $\sL( \mB)\equiv \sL(\valpha, \mB)$ for the corresponding loss function with respect to the spline coefficient matrix in~\eqref{eqn:lossfun:mat}.
Taking expectation with respect to both the responses $y_{nm}$'s and functional covariates $x_{nm}$'s,
we denote the expected loss functions $\bar{\sL}(\vbeta) = \Expect \sL(\vbeta)$ and $\bar{\sL}(\mB) = \Expect \sL(\mB)$.

To quantify the approximation error for the general model~\eqref{eqn:method:mixed}, 	
we define two versions of optimal spline coefficient matrix $\mB$ associated with the expected loss $\bar{\sL}(\mB)$.
The first one is the \textit{unconstrained optimal parameter} $\bmB_0$. It is computed with the expected loss $\bar{\sL}$ and the original penalty $\sP_{\veta}(\mB)$, but
without the manifold constraint $\mathbb{M}$, i.e.,
\begin{equation} \label{def:popuestimate}
	\bmB_0:=\argmin_{\mB\in\bbR^{K\times M}} 
	\bar{\sL} (\mB)
	+ \sP_{\veta}(\mB).
\end{equation}
In addition to~\eqref{def:popuestimate}, we define the \textit{constrained  optimal parameter} $\bmB$, which is computed under the constraint $\mathbb{M}$ imposed upon the spline coefficient matrix, i.e.,
\begin{equation} \label{def:popuestimate:constraint}
	\bmB:=\argmin_{\mB\in\mathbb{M}} \bar{\sL} (\mB)
	+ \sP_{\veta}(\mB).
\end{equation}
To avoid the intricacy of multiple optimal solutions, we assume the objective function in~\eqref{def:popuestimate:constraint} is strictly convex in a local neighbor of $\bmB$ over $\mathbb{M}$. Equivalently, the intersection between the level set
$ \big\{\mB\in\bbR^{K\times M}:\;\bar{\sL} (\mB) + \sP_{\veta}(\mB)
= \bar{\sL} (\bmB)+ \sP_{\veta}(\bmB)\big\}$
and a small neighbor of $\bmB$ over $\mathbb{M}$ is trivially the single point $\{\bmB\}$.

The optimal parameters help us quantify the overall model approximation and penalty biases. Given $\bmB_0 = (\bar{\vb}_{01},\cdots, \bar{\vb}_{0M})$ from~\eqref{def:popuestimate}, we set $\bar{\vbeta}_0=(\bar{\beta}_{01},\cdots, \bar{\beta}_{0M})\trans$ with $\bar{\beta}_{0m}(\cdot)=\vphi\trans(\cdot)\bar{\vb}_{0m}$, $m =1,\dots, M$.
The \textit{spline approximation error} $\sE(\mathbb{S}_K)$ for the general model~\eqref{eqn:method:mixed} is defined as
\begin{equation} \label{eqn:perror:spline}
	\sE(\mathbb{S}_K):= \bigg\{\sum_{m=1}^{M} \| \bar{\beta}_{0m} - \beta_{0m}\|_{X}^2 + \sP_{\veta} (\bar{\vbeta}_0)\bigg\}^{1/2}.
\end{equation}
The quantity $\sE(\mathbb{S}_K)$ can be interpreted as the bias due to modeling the slope function $\beta_{0m}$ in the spline space $\mathbb{S}_K$ with a penalization term in our model. 
Similarly, using $\bmB = (\bar{\vb}_{1},\cdots, \bar{\vb}_{M})$ defined in~\eqref{def:popuestimate:constraint}, we set $\bar{\vbeta}=(\bar{\beta}_{1},\cdots, \bar{\beta}_{M})\trans$ with $\bar{\beta}_{m}(\cdot)=\vphi\trans(\cdot)\bar{\vb}_{m}$, $m=1,\dots, M$.
The additional \textit{manifold constraint error} $\sE(\mathbb{M})$ is quantified as 
\begin{equation} \label{eqn:perror:manifold}
	\sE(\mathbb{M}) := 	\bigg\{\sum_{m=1}^{M} \| \bar{\beta}_{0m} - \bar{\beta}_{m}\|_{X}^2 + \sP_{\veta} (\bar{\vbeta}_{0} - \bar{\vbeta})\bigg\}^{1/2}.
\end{equation}
The above $\sE(\mathbb{M})$ compares the difference between  $\bar{\vbeta}$ and $\bar{\vbeta}_0$. 
We will use both $\sE(\mathbb{S}_K)$ and $\sE(\mathbb{M})$ to describe the overall model bias in our analysis of convergence rate for the penalized spline estimator under the manifold constraint $\mathbb{M}$ over the spline coefficient matrix.

\section{Manifold Local Complexity}\label{sec:manifold}

Quantifying the estimation error of~\eqref{eqn:method:mixed} amounts to examining a loss-related empirical process $\sV(\mB)$ indexed by $\mB\in\mathbb{M}$ as
\begin{equation}  \label{eqn:localEP}
	\sV(\mB) :=  \sL (\mB) - 
	\bar{\sL} (\mB)= \frac{1}{N}\sum_{m=1}^M\sum_{n=1}^N
	\big\{\ell_m(y_{nm}, \vx_{nm}\trans\vb_m) - \Expect\,\ell_m(y_{nm}, \vx_{nm}\trans\vb_m)\big\}.
\end{equation}
We consider to control the magnitude of $\sV(\mB)-\sV(\bmB)$
for $\mB$ in a local neighbor of  $\bmB$ over the manifold $\mathbb{M}$.  
To develop the upper bound, we first review a few concepts and  notations for submanifold in Section~\ref{sec:manifold:concept}. 
Interested readers are referred to \cite{Lee2018} for a detailed description of manifolds.
After that, in Section~\ref{sec:manifold:epbound}, we study the complexity of a manifold local neighbor induced by the penalty~\eqref{eqn:method:generalPenalty} and obtain the upper bound of the uniform magnitude of $\sV(\mB)-\sV(\bmB)$ in the  local neighbor.

\subsection{Review of Riemannian Submanifold}
\label{sec:manifold:concept}

In model~\eqref{eqn:method:mixed}, we consider the constraint set $\mathbb{M}$ ($\subset \bbR^{K\times M}$) as a Riemannian embedded submanifold without boundary. 
The manifold $\mathbb{M}$ is a subset of matrices that is locally homeomorphic to the Euclidean space \citep{Lee2018}. 
At any $\mB\in\mathbb{M}$, the manifold $\mathbb{M}$ is approximated by a tangent space $T_{\mB}\mathbb{M}$ to the first order.  
We set the metric of $\mathbb{M}$ as being induced from the ambient space $\bbR^{K\times M}$. 
This means, at any $\mB\in\mathbb{M}$, the metric value of two tangent vectors is simply the value of their Euclidean inner product.
The second order structure is induced by connection over manifold.
Suppose $\mDelta,\mDelta'$ are two tangent vector fields and $\widetilde{\nabla}$ is the Euclidean connection of $\bbR^{K\times M}$. 
For the ambient connection, $\widetilde{\nabla}_{\mDelta} \mDelta'$ can be viewed as the direction derivative of $\mDelta'$ in the direction of $\mDelta$ in the Euclidean space $\mathbb{R}^{K\times M}$. 
The Levi-Civita connection $\nabla$ for $\mathbb{M}$ can then be determined via $\nabla_{\mDelta} \mDelta' =  \mathrm{P}_{\mB}(\widetilde{\nabla}_{\mDelta} \mDelta')$, where $\mathrm{P}_{\mB}$ is the orthonormal projection onto the tangent space $T_B\mathbb{M}$.

A \textit{geodesic} $\gamma(t, \mDelta)$ is a smooth curve over $\mathbb{M}$ indexed by $t$ in an interval including $0$. 
The geodesic  starts at $\mB= \gamma(0, \mDelta)$ with initial velocity $\dot{\gamma}(0, \mDelta) = \mDelta \in T_{\mB}\mathbb{M}$, and has zero acceleration (i.e., $\nabla_{\dot \gamma(t, \mDelta)} \dot \gamma(t, \mDelta) = \vzero$) in the tangent space. 
The geodesic defines the \textit{exponential mapping} $\exp_{\mB}(\cdot)$ which maps a tangent vector $\mDelta\in T_{\mB}\mathbb{M}$ to $\exp_{\mB}( \mDelta) = \gamma(1, \mDelta)$. In particular, it maps the zero tangent vector $\vzero\in  T_{\mB}\mathbb{M}$ to the point $\mB$ itself, i.e., $\exp_{\mB}(\vzero) = \mB$. 
The domain $\sD_{\mB}$ of $\exp_{\mB}(\cdot)$ is a star-shaped subset of $T_{\mB}\mathbb{M}$ containing $\vzero$ \citep[Proposition~5.19 of][]{Lee2018}.
When $\mathbb{M}$ is complete, the domain $\sD_{\mB} = T_{\mB} \mathbb{M}$ is the full tangent space. 
Let $B(\vzero,r) = \{\mDelta\in T_{\mB} \mathbb{M}:\ \|\mDelta\|_F\le r\}$ be the ball with radius $r$ in the tangent space. 
The \textit{injective radius} ($\mathrm{inj}(\mB)$) at $\mB$ is the supermum of $r$ such that the exponential mapping $\exp_{\mB}(\cdot)$ is a diffeomorphism over $B(\vzero,r) \subseteq T_{\mB} \mathbb{M}$. 

Our theory will restrict the curvature of the submanifold, where the curvature is quantified via second fundamental form. 
The \textit{second fundamental form} ${\rm II}(\cdot,\cdot)$ is a mapping from the product of two tangent vector fields onto the normal vector bundle \citep[see Chapter~8 of][]{Lee2018}.
It holds that ${\rm II}(\mDelta,\mDelta') = \mathrm{P}^{\perp}_{\mB}\big(\widetilde{\nabla}_{\mDelta} \mDelta'\big)$, where $\mathrm{P}^{\perp}_{\mB}$ at $\mB$ is the projection onto the normal space $N_{\mB}\mathbb{M}=\big( T_{\mB}\mathbb{M}\big)^{\perp}$. 
Given a geodesic $\gamma(\mDelta, t)$ (which can also be viewed as a curve of $\mathbb{R}^{K\times M}$), its acceleration vector $\ddot\gamma(\mDelta, t)$ in the ambient space $\mathbb{R}^{K\times M}$ can be computed from the second fundamental form via $\ddot\gamma(\mDelta, t) = {\rm II}\big(\dot\gamma(\mDelta, t),\dot\gamma(\mDelta, t)\big)$. 
Meanwhile, as $\mathbb{M}$ is a submanifold of the Euclidean space $\bbR^{M\times K}$, its curvature tensor is determined by its second fundamental form $\rm II(\cdot,\cdot)$ due to the Gaussian Equation \citep[see Theorem~8.5 of][]{Lee2018}.

\subsection{Local Empirical Process over Manifold} \label{sec:manifold:epbound}

Controlling the magnitude of $\sV(\mB)-\sV(\bmB)$ in a proper geodesic neighbor of $\bmB$ will assist to derive the estimation error for the doubly regularized estimator~\eqref{eqn:method:mixed}. 
Inspired by the ellipsoid in~\eqref{eqn:complexity0:ellipsoid}, 
we introduce a norm $\sQ_{\veta}(\cdot)$ to determine the size of a neighbor set around $\vzero$ in the tangent space  $T_{\bmB}\mathbb{M}$.
In particular, for a matrix $\mB\in\bbR^{K\times M}$, we can define $\sQ_{\veta}(\mB)$ via
\begin{equation}\label{eqn:defineQnorm}
	\sQ^2_{\veta}(\mB) := 
	\|\mB\|_F^2 + \sP_{\veta} (\mB)= \sum_{m=1}^M \|\beta_m\|_{X}^2 + \sP_{\veta}(\vbeta),
\end{equation}
which combines the Frobenius norm and the general composite quadratic penalty~\eqref{eqn:method:generalPenalty}. 
Now, let $\mathbb{N}(\bmB,\delta)$ denote the local neighbor of $\vzero$ in the tangent space $T_{\bmB}\mathbb{M}$ determined by $\sQ_{\veta}(\cdot)$ via
\begin{equation} \label{eqn:generalNeighbor}
	\mathbb{N}(\bmB, \delta ) = \big\{\mDelta\in \sD_{\bmB}:\ \sQ_{\veta}(\mDelta) \le \delta \big\}.  
\end{equation}
The exponential mapping maps $\mathbb{N}(\bmB, \delta)$ back to the manifold via 
$$\exp_{\bmB}\big(\mathbb{N}(\bmB, \delta )\big)
=\{\exp_{\bmB}(\mDelta):\;  \mDelta \in \mathbb{N}(\bmB,\delta)\},$$
The set $\exp_{\bmB}\big(\mathbb{N}(\bmB, \delta)\big)$ is a local geodesic neighbor of $\bmB$ over $\mathbb{M}$. 
We will apply the generic chaining to both $\mathbb{N}(\bmB, \delta)$ and $\exp_{\bmB}\big(\mathbb{N}(\bmB, \delta )\big)$. Generic chaining will help properly characterize their complexities and control the magnitude of $\sV(\mB)-\sV(\bmB)$ over the geodesic neighbor $\exp_{\bmB}\big(\mathbb{N}(\bmB, \delta )\big)$.

\begin{remark}
The neighbor $\mathbb{N}(\bmB, \delta )$ in~\eqref{eqn:generalNeighbor} is ellipsoid-like. To see this, we only need to rewrite the squared norm in~\eqref{eqn:defineQnorm} as
\begin{equation}\label{eqn:ellNormMain}
\sQ^2_{\veta}(\mB) =  \|\mB\|_F^2 +  \sum_{j=1}^{P} \eta_j\tr\big(\mB\trans\mPi_{j1}\mB\mPi_{j2}\big)
= \vb\trans\Bigg( \mI + \sum_{j=1}^{P} \eta_j  \mPi_{j2}\otimes \mPi_{j1}\Bigg) \vb,
\end{equation}
where $\vb = \mathrm{vec}(\mB)$ is the vectorization of the matrix $\mB\in\mathbb{R}^{K\times M}$. Without a manifold constraint (i.e., $\mathbb{M} =\bbR^{K\times M}$), this norm exactly induces an ellipsoid neighbor in the Euclidean space $\mathbb{R}^{K\times M}$ due to the penalty associate term $\sum_{j=1}^{P} \eta_j  \mPi_{j2}\otimes \mPi_{j1}$ in \eqref{eqn:ellNormMain}. When a proper
submanifold constraint (i.e., $\mathbb{M} \subset \bbR^{K\times M}$) is considered and the neighbor is small enough, this norm also induces an ellipsoid neighbor in the tangent space of the submanifold. 
Generic chaining is known to provide sharp characterization of the ellipsoid complexity, while Dudley's bound may fail to do so. As noted in Section~2.5 of \cite{talagrand2014upper}, the complexity characterization of a general ellipsoid via the Dudley's bound can be worse by a factor of $\sqrt{\log(\xi +1)}$, where $\xi$ is the intrinsic dimension of the constraint manifold $\mathbb{M}$ and it can be as large as $\xi = MK$. 
\end{remark}

Because the manifold $\mathbb{M}$ at $\bmB$ is approximated by its tangent space $T_{\bmB}\mathbb{M}$ on the first order, it can be expected that, as long as the manifold $\mathbb{M}$ has a bounded curvature, the two local sets, $\mathbb{N}(\bmB, \delta )$ and $\exp_{\bmB}\big(\mathbb{N}(\bmB, \delta )\big)$, should have the same $\gamma_{\alpha}$ complexity level. 
We restrict the curvature of $\mathbb{M}$ by the following condition. 

\begin{condition} \label{ass:manifold}
	The manifold $\mathbb{M}$ at $\bmB$ has strictly positive injective radius ($\mathrm{inj}(\bmB)>0$). 
	In addition, there exists a constant $C_{\mathrm{II}}>0$, such that the second fundamental form $\mathrm{II}(\cdot, \cdot)$ of the manifold $\mathbb{M}$ is bounded with respect to $\sQ_{\cdot} (\mDelta)$:
	\begin{equation} \label{eqn:manifoldass}
		\sQ_{\veta} ( \mathrm{II}(\mDelta,\mDelta) ) \le C_{\mathrm{II}}
		\sQ_{\veta}^2 (\mDelta),
	\end{equation}
	for all $\mDelta\in T_{\mB}\mathbb{M}$ and all $\mB$ in a local neighbor of $\bmB$.
\end{condition}

Since the second fundamental form is bilinear, \eqref{eqn:manifoldass} is equivalent to the requirement that $	\sQ_{\veta} ( \mathrm{II}(\mDelta,\mDelta) ) \le C_{\mathrm{II}}$ for all $\mDelta$ satisfying $\sQ_{\veta} (\mDelta)\le 1$.
Based on this condition, we can show the manifold is locally close to the tangent space as measured by both the Frobenius norm and the norm of $\sQ_{\veta}$.

\begin{lemma} \label{lemma:metriccompatible}
	Under Condition~\ref{ass:manifold}, there exists a radius $R_M>0$ (depending on $C_{\rm II}$), such that for all $\mDelta\in  \mathbb{N}(\bmB, R_M)$, it holds that
	\begin{equation}  \label{lemma:metriccompatible1}
	(1/2)	\sQ_{\veta} (  \mDelta)\le
	    \sQ_{\veta} (\exp_{\bmB}(\mDelta) - \bmB) \le
	    2\sQ_{\veta} (  \mDelta),
    \end{equation}
    and that
	\begin{equation} \label{eqn:manifold:secondorderappr}
	    \sQ_{\veta} (	\exp_{\bmB}(\mDelta) - \bmB - \mDelta) 
	    \le  2  \sQ_{\veta}^2 \big( \mDelta\big).
    \end{equation}
    Besides, for any $\mDelta_1,\mDelta_2\in  \mathbb{N}(\bmB, R_M)$, we have
	\begin{equation}  \label{lemma:metriccompatible2}
		(1/4) \|  \mDelta_1 - \mDelta_2 \|_F\le
		\big\|\exp_{\bmB}(\mDelta_1) - \exp_{\bmB}(\mDelta_2)\big\|_F \le
		4 \|  \mDelta_1 - \mDelta_2 \|_F.
	\end{equation}
\end{lemma}

In the above, \eqref{lemma:metriccompatible1} indicates that, in terms of the norm $\sQ_{\veta}$, the magnitude of the deviation between $\exp_{\bmB}(\mDelta)$ and $\bmB$ has the same order of that of $\mDelta$.  The bound for~\eqref{lemma:metriccompatible1} in the special case with $\veta=\vzero$ (i.e., $\sQ_{\veta} ( \cdot)$ is simply the Frobenius norm) has been  used in the literature \citep[e.g.,][]{trillos2020error,berenfeld2021density}. For the second result~\eqref{eqn:manifold:secondorderappr}, observe that $\bmB + \mDelta$ is a first-order approximation to $\exp_{\bmB}(\mDelta)$, and thus the term $\exp_{\bmB}(\mDelta) - \bmB - \mDelta$ can be viewed as the error from high orders.  The result~\eqref{eqn:manifold:secondorderappr} implies that the magnitude of the high-order error can by controlled by $2\sQ_{\veta}^2 \big( \mDelta\big)$.  The third result \eqref{lemma:metriccompatible2} states the length of the difference of two tangent vectors is compatible with the ambient distance between their images under the exponential mapping.  To our best knowledge, the above results are novel in the literature with a general norm of $\sQ_{\veta}(\cdot)$ and a pair of $\mDelta_1,\mDelta_2$.

\begin{figure}
	\centering
	\includegraphics[width = 0.8\textwidth]{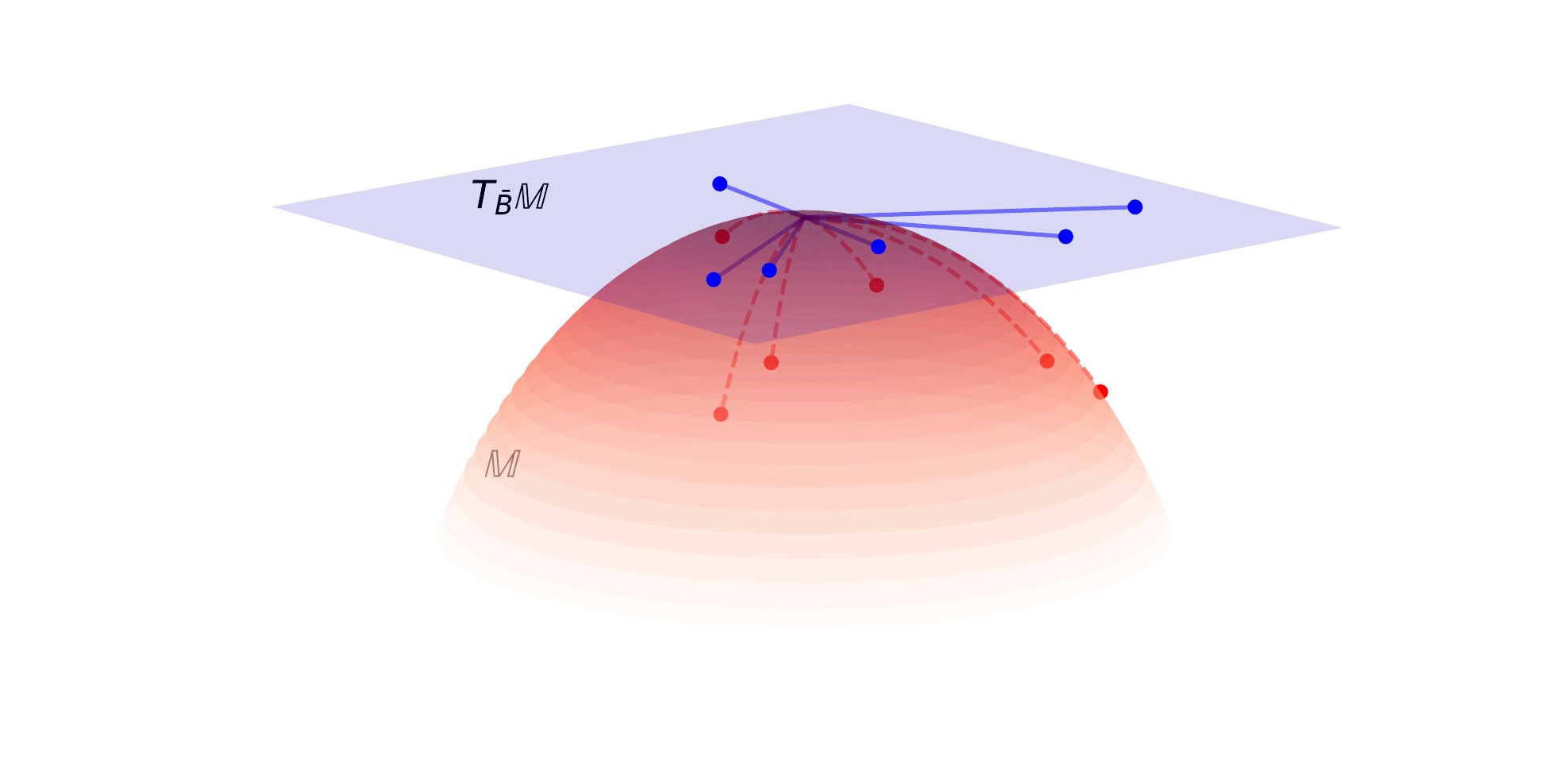}
	\caption{An illustration for the sets $\bar{T}_n$ and $\bar{S}_n$. The red curved surface is the manifold $\mathbb{M}$. The blue hyperplane is the tangent space $T_{\bmB}\mathbb{M}$. The blue points represent the set $\bar{T}_n$ in the tangent space. 
	Through the exponential mapping (the red dashed curves), the blues points are mapped onto the manifold. 
	The resulting red points represent the set $\bar{S}_n$. \label{fig:tangentfig}}
\end{figure}

Lemma~\ref{lemma:metriccompatible} allows us to transfer the $\gamma_2$ functional of the local set $T=\mathbb{N}(\bmB,\delta)$  in the tangent space to the $\gamma_2$ functional of  $S=\exp_{\bmB}(\mathbb{N}(\bmB,\delta))$ over the manifold $\mathbb{M}$. 
The $\gamma_2$ functional of these local sets are computed with the Frobenius norm, e.g., $\mathbb{N}(\bmB,\delta)$ is endowed with the metric $d(\mDelta_1, \mDelta_2) = \|\mDelta_1 - \mDelta_2\|_F$. 
Suppose $\{\bar{T}_n\}$ is an admissible sequence of subsets of $T = \mathbb{N}(\bmB,\delta)$ (with $\delta < R_M$) satisfying
\begin{align}
	\sup_{\mDelta \in \mathbb{N}(\bmB,\delta)} \sum_{n=0}^{\infty} 2^{n/2}d(\mDelta, \bar{T}_n) \le 2 \gamma_{2}(\mathbb{N}(\bmB,\delta), \|\cdot\|_F).\label{eqn:adimissibleT}
\end{align}
We map each set $\bar{T}_n$ from the tangent space to the manifold via
$$
\bar{S}_n = \exp_{\bmB}(\bar{T}_n) := \{\mB:\;\mB = \exp_{\bmB}(\mDelta)\text{ for } \mDelta\in \bar{T}_n\}.
$$
See Figure~\ref{fig:tangentfig} for an illustration of the sets $\bar{T}_n$ and $\bar{S}_n$.
It is readily checked that $\{\bar{S}_n\}$ is an admissible sequence of subsets of the geodesic neighbor $S=\exp_{\bmB}(\mathbb{N}(\bmB,\delta))$, since the cardinality of $\bar{S}_n$ is the same as that of $\bar{T}_n$. 
Then, it holds that
\begin{align}
	\gamma_{2}(S, \|\cdot\|_F) & =  \inf_{\{S_n\}}\sup_{\mB \in S} 
	\sum_{n=0}^{\infty} 2^{n/2}d(\mB, S_n)   \nonumber  \\ 
	&\le \sup_{\mB \in S} \sum_{n=0}^{\infty} 2^{n/2}d(\mB, \bar{S}_n)  \nonumber  \\
	& =  \sup_{\mDelta \in \mathbb{N}(\bmB,\delta)} \sum_{n=0}^{\infty} 2^{n/2}d(\exp_{\bmB}(\mDelta), \exp_{\bmB}(\bar{T}_n) )  \nonumber  \\
	& \stackrel{(i)}{\le} 4	\sup_{\mDelta \in \mathbb{N}(\bmB,\delta)} \sum_{n=0}^{\infty} 2^{n/2}d(\mDelta, \bar{T}_n) \nonumber  \\
    & \stackrel{(ii)}{\le} 8 \, \gamma_{2}(T, \|\cdot\|_F) . \label{eqn:gammaalpha:mani2tan}
\end{align}
In the above, the inequalities~(i) and (ii) are due to \eqref{lemma:metriccompatible2} of Lemma~\ref{lemma:metriccompatible} and \eqref{eqn:adimissibleT}, respectively. 
By a similar argument, when $\delta < \min\{\mathrm{inj}(\bmB), R_M\}$, we can also show it holds
$\gamma_{2}(S, \|\cdot\|_F) \ge \gamma_{2}(T, \|\cdot\|_F)/8$, and in this case, we can conclude that
$	\gamma_{2}(S, \|\cdot\|_F) \asymp  \gamma_{2}(T, \|\cdot\|_F) $.
This means that the $\gamma_2$ functional of the sets $T=\mathbb{N}(\bmB,\delta)$ and $S=\exp_{\bmB}(\mathbb{N}(\bmB,\delta))$ are of the same order, as long as the radius $\delta$ is sufficiently small.

The above discussion reveals that, to control $\sV(\mB)-\sV(\bmB)$	for $\mB$ in a local neighbor of $\bmB$, it suffices to consider the complexity of the local set $\mathbb{N}(\bmB,\delta)$ in the tangent space.
Meanwhile, bounding $\sV(\mB) - \sV(\bmB)$ also requires quantifying the random oscillation of each summand in~\eqref{eqn:localEP}. We impose a Lipschitz continuous assumption on the loss function.

\begin{condition} \label{ass:lossLip}
	For $m=1,\cdots, M$, the loss functions $\ell_m(y, u)$'s are uniformly Lipschitz continuous with respect to $u$, i.e.,
	\begin{equation} \label{ass:lossfun:lip}
		|\ell_m(y, u)-\ell_m(y, u')| \le C_{L} |u - u'|
	\end{equation}
	for some constant $C_L$.
\end{condition}

The Lipschitz continuous assumption is commonly found in the literature \citep[e.g.,][]{van2008high,geoffrey2020robust}.  
Examples of such loss function include the logistic regression loss, huber loss and quantile regression loss, among many others. 
The least squares loss also meets the requirement when the functional covariate $x_{nm}$ and the response $y_{nm}$ are bounded. 
Condition~\ref{ass:lossLip} simply provides a convenient way for analyzing the empirical process $\sV(\mB)- \sV(\bmB)$.  The proposition below relates the magnitude of the empirical process to  
the complexity $ \gamma_2( \mathbb{N}(\bmB,\delta),\; \|\cdot\|_F) $.

\begin{proposition} \label{lemma:epupper}
	Under Conditions \ref{ass:subGaussian}--\ref{ass:lossLip} and for $\delta< R_M$, with probability at least $1 - 2 \exp(-u^2)$, it holds that
	\begin{align*}
		\sup_{\mDelta \in \mathbb{N}(\bmB, \delta)}	\Big|\sV(\exp_{\bmB}(\mDelta)) -	\sV(\bmB) \Big|  
		\le      \frac{C_{V}C_L(C_g+2)}{N^{1/2}} \big\{ \gamma_2( \mathbb{N}(\bmB,\delta),\; \|\cdot\|_F)   + u \delta\big\},
	\end{align*}
	where $C_{V}$ is a positive absolute constant.
\end{proposition}

\begin{remark}
In the special case that $\mathbb{M} = \mathbb{R}^{K\times M}$ (i.e., without a proper manifold constraint), it follows from Eqn. (2.115) and Theorem 4.1.11 of \cite{talagrand2014upper} that we readily have the general expression
$$
\gamma_2( \mathbb{N}(\bmB,\delta),\; \|\cdot\|_F) \asymp \Big(\sum_{j=1}^{MK} \frac{\delta^2}{1+\lambda_j(\mathbf{P} )}\Big)^{1/2},
$$
where $\mathbf{P}:=\sum_{j=1}^{P} \eta_j  \mPi_{j2}\otimes \mPi_{j1}$ is the quadratic penalty associate term in~\eqref{eqn:ellNormMain} and $\lambda_j(\mathbf{P} )$ is its $j$-th largest eigenvalue.
\end{remark}

\section{The Main Result}
\label{sec:maintheory}

We are now ready to derive the unified upper bound for the estimator of the general model~\eqref{eqn:method:mixed} with double regularization.
We consider the finite sample loss for the $m$-th task  $$\sL_{m}(\beta_m) :=  (1/N) \sum_{n=1}^{N}\ell_m\Big(y_{nm}, \int x_{nm}(t)\beta_m(t)\intd t\Big),$$ 
and its expected counterpart $\bar{\sL}_{m}(\beta_m) :=  \Expect \sL_{m}(\beta_m)$, where the expectation is taken with respect to both the response $y_{nm}$ and functional covariate $x_{nm}$.   We view the losses $\sL_{m}(\beta_m)$ and $\bar{\sL}_{m}(\beta_m)$  as functions of $\beta_m \in \mathbb{L}_2(\sT)$. The next condition states the expected loss $\bar{\sL}_{m}(\beta)$ is strongly convex and smooth, for $\beta$ in a local neighbor of the true $\beta_{0m}$. The locality means the norm $\|\beta - \beta_{0m}\|_X$ is small.

\begin{condition} \label{ass:lossfunConvex}
Denote $E_m(\beta,\beta'):=\bar{\sL}_{m}(\beta) - \bar{\sL}_{m}(\beta') - D \bar{\sL}_{m}(\beta')[\beta - \beta']$, where $D \bar{\sL}_{m}(\beta)[\cdot]$ is the Fr\'echet derivative of $ \bar{\sL}_{m}$ at $\beta$ in $\mathbb{L}_2(\sT)$.   
For $m=1,\cdots, M$, there exist constants $C_c > c_c > 0$ such that 
\begin{equation} \label{ass:lossfun:convex}
	c_c \| \beta - \beta'\|_{X}^2 \le E_m(\beta,\beta') \le C_c \| \beta - \beta'\|_{X}^2,
\end{equation}
with $\beta$ and $\beta'$ in a local neighbor of the true $\beta_{0m}$. 
\end{condition}

Examples of loss function satisfying  Condition~\ref{ass:lossfunConvex} include the least squares loss, logistic regression loss, etc. 
For the least squares loss, we can see $\bar{\sL}_{m}(\beta)  = \Expect (y_{1m} - \langle x_{1m}, \beta\rangle )^2$ and $D \bar{\sL}_{m}(\beta')[\beta - \beta'] = -2\Expect ( y_{1m} \langle x_{1m}, \beta-\beta'\rangle)$. It follows $E_m(\beta,\beta') = \Expect  \langle x_{1m}, \beta-\beta'\rangle^2 = \|\beta - \beta'\|_X^2$, and Condition~\ref{ass:lossfunConvex} holds with $C_c = c_c = 1$. In Section~\ref{sec:appendix:mainproof:quantile} of the Appendix, we show the loss of quantile regression also satisfies  Condition~\ref{ass:lossfunConvex} under mild regularity conditions.

According to Lemma~\ref{lemma:epupper}, it shows that the local empirical process can be controlled by
$\gamma_2( \mathbb{N}(\bmB,\delta),\; \|\cdot\|_F)$. Based on the $\gamma_2$ functional, we can determine the estimation error by the \textit{critical radius} $\widehat{\delta}_N$ of
\begin{equation} \label{eqn:criticalradius}
	\widehat{\delta}_N:=\inf\Big\{\delta>0:\ 
	C_{V}C_L(C_g+2)\gamma_2( \mathbb{N}(\bmB,\delta),\; \|\cdot\|_F)  \le  \delta^2N^{1/2} \Big\}.
\end{equation}
The above definition indicates the two functions
$\gamma_2( \mathbb{N}(\bmB,\delta),\; \|\cdot\|_F) /\delta$ and $N^{1/2}\delta $ intersect at the critical radius $\widehat{\delta}_N$.
Note $N^{1/2}\delta $ as a function of $\delta$ is linear with slope $N^{1/2}$. 
It is also easy to check $\gamma_2( \mathbb{N}(\bmB,\delta),\; \|\cdot\|_F) /\delta$ is a constant for $\delta\in(0,\mathrm{inj}(\bmB))$. 
The value of $\widehat{\delta}_N$ is therefore well-defined for large enough $N$. 

\begin{remark}
Similar definition of critical radius can be found in  \cite{wainwright2019high} and \cite{yang2017randomized}.   In Section~13.2 of \cite{wainwright2019high}, the critical radius $\delta^*$ for non-parametric regression is defined as 
\begin{equation} \label{eqn:criticalradius:wainwright}
    \delta^*:=\inf\Big\{ \delta>0:\;
    2\sigma \mathcal{G}_n(\delta, \mathcal{F}^*)   \le \delta^2
\Big\},
\end{equation}
where $\mathcal{F}^*$ is a localized function class, $\sigma$ is the standard deviation of additive noise. In addition,   $\mathcal{G}_n(\delta, \mathcal{F}^*)$ is the local Gaussian complexity
$$
\mathcal{G}_n(\delta, \mathcal{F}^*)  = \Expect \bigg\{\sup_{g\in\mathcal{F}^*,\ \|g\|_n\le \delta }
\frac{1}{N} \Big| \sum_{n=1}^{N} e_i g(x_i) \Big| \bigg\},
$$
where $e_i$ i.i.d follows the standard Gaussian distribution and $\|g\|_n^2 = (1/N)\sum_{n=1}^{N} g(x_i)^2$ is the related empirical norm.
We remark that, although we study a very different model  and use generic chaining to quantify complexity, 
the key difference between~\eqref{eqn:criticalradius} and~\eqref{eqn:criticalradius:wainwright} is that we use the proposed norm $\sQ_{\veta}(\cdot)$ to define the local set in~\eqref{eqn:generalNeighbor}. This key difference helps us to reveal the phase transition behavior of the estimators in Sections~\ref{sec:reducedrate} and \ref{sec:graphrate}.
\end{remark}

Combining the critical radius $\widehat{\delta}_N$ with the spline approximation error $\sE(\mathbb{S}_K)$ and the manifold constraint error $\sE(\mathbb{M})$ (in Section~\ref{sec:approError}), we establish a unified upper bound of the convergence rate for the estimator of the general model~\eqref{eqn:method:mixed} with double regularization.

\begin{theorem} \label{theorem:main}
	Suppose Conditions~\ref{ass:beta}--\ref{ass:lossfunConvex} hold, and define $\ubar{c}_c = \min\{c_c,1\}$ and $\bar{C}_c = \max\{C_c, 1\}$. Assume the manifold constraint error is sufficiently small such that
	$	\sE(\mathbb{M})\le	 \ubar{c}_c / (32\bar{C}_c) $.
    Then, for a sufficiently large $N$ and a given $u$, it holds with probability at least $1- 2\exp(-u^2)$ that
	there exists a local optimal estimate $\hat{\vbeta} = \big(\hat{\beta}_1,\cdots, \hat{\beta}_M \big)\trans$ of the model~\eqref{eqn:method:mixed} satisfying
	\begin{align}
		(1/M) \Big\{\sum_{m=1}^{M} & \|\hat{\beta}_m-\beta_{0m}\|_{X} +  \sP_{\veta}^{1/2}(\hat{\vbeta})\Big\} \nonumber \\
		& \le
		\frac{C_U}{\ubar{c}_cM^{1/2}}
		\Big [\ubar{c}_c \big\{ \sE(\mathbb{M}) + \sE(\mathbb{S}_K) \big\}
		+  \widehat{\delta}_N + C_L(C_g+2)u/\sqrt{N}\Big], \label{eqn:mainbound}
	\end{align} 
	where $C_U$ is an absolute constant.
\end{theorem}

In the above theorem, it is required the manifold constraint error $\sE(\mathbb{M})$ is sufficiently small. 
This is satisfied for the reduced model~\eqref{eqn:method:reduced}  if we set $R$ to be moderately large. Generally, for large enough $N$, we can relax the manifold constraint $\mathbb{M}$ to a larger subset of $\bbR^{K\times M}$ such that $\sE(\mathbb{M})$ is small.  Alternatively, this condition on $\sE(\mathbb{M})$ can be removed if 
the expected loss $\bar{\sL}(\mB)$ is geodesically strongly convex in a neighbor of $\bmB$ over the manifold $\mathbb{M}$. 

In the following sections, we apply Theorem~\ref{theorem:main} to specific models, such as the reduced model~\eqref{eqn:method:reduced} or the graph regularized model~\eqref{eqn:method:graph}. 
We will derive the quantities $\sE(\mathbb{S}_K)$, $\sE(\mathbb{M})$, and $\widehat{\delta}_N$ for these models. 
Some interesting phase transition behaviors can be concluded after plugging these obtained values into~\eqref{eqn:mainbound}.

\section{Application I: The Reduced Multi-task Model} 
\label{sec:reducedrate}

In this section, we consider the reduced model~\eqref{eqn:method:reduced}, where the constraint set $\mathbb{M}$ is the set of rank-$R$ matrices, which forms a fixed-rank manifold.
We apply Theorem~\ref{theorem:main} to derive the rate of convergence when $N$ diverges to infinity and the number of task $M$ is fixed. For simplicity, we assume $K\ge M$ in the following.


For the unconstrained optimal parameter as in~\eqref{def:popuestimate}, we can compute its singular value  decomposition (SVD) $\bmB_0 = \bmU_0\bmD_0\bmV_0\trans$, where $\bmU_0\in\bbR^{K\times M}$ and $\bmV_0\in\bbR^{M\times M}$ are two orthonormal matrices with the singular vectors in their columns, and $\bmD_0 = \diag(\bar{\sigma}_{01},\dots, \bar{\sigma}_{0M})$ is a diagonal matrices with non-increasing singular values. 
It is well-known that the best rank-$R$ approximation to $\bmB_0$ (in terms of Frobenius norm) is obtained by truncating the SVD and only keeping the leading $R$ singular values with the associated vectors. 
We usually interpret the quantity $\sum_{r=R+1}^{M} \bar{\sigma}_{0r}^2$ as the rank-$R$ approximation error. The following lemma provides a more precise bound on the spline approximation error $\sE(\mathbb{S}_K)$ and rank-$R$ constraint error $\sE(\mathbb{M})$.

\begin{lemma} \label{lemma:reduced:approxerror}
	(i) The squared spline approximation error $\sE(\mathbb{S}_K)$ is bounded by		
	\begin{equation} \label{lemma:reduced:approxerror1}
		\big\{\sE(\mathbb{S}_K)\big\}^2 =  \sum_{m=1}^{M}
		\Big\{\|\bar{\beta}_{0m}- \beta_{0m}\|_{X}^2 
		+\eta_1\|\bar{\beta}_{0m} \|_{\Gamma}^2\Big\}
		\lesssim M(K^{-2\tau} + \eta_1  K^{2(d-\nu)_{+}}),
	\end{equation}
	with $\tau = \nu \wedge (\mathfrak{o}+1) + \{q \wedge (\mathfrak{o}+1)\}/2$. 
	
	\noindent(ii) When $R<M$, the rank-$R$ constraint error $\sE(\mathbb{M})$ satisfies
	\begin{equation} \label{lemma:reduced:approxerror2}
		\big\{\sE(\mathbb{M})\big\}^2= \sum_{m=1}^{M}\Big\{
		\|\bar{\beta}_m- \bar{\beta}_{0m}\|_{X}^2 +  	\eta_1\|\bar{\beta}_m - \bar{\beta}_{0m} \|_{\Gamma}^2 \Big\} \lesssim \sum_{r=R+1}^{M}\bar{\sigma}_{0r}^2+ 	\eta_1\sum_{m=1}^{M}
		\| \bar{\beta}_{0m} \|_{\Gamma}^2.
	\end{equation}
	Otherwise, when $R=M$, the constraint error is zero, i.e., $\sE(\mathbb{M}) = 0$.
\end{lemma}

For the constrained optimal parameter $\bmB\in\mathbb{M}$ in~\eqref{def:popuestimate:constraint},
suppose it has compact SVD $\bmB = \bmU\bmD\bmV\trans$ where $\bmD\in\bbR^{R\times R}$ is a diagonal matrix of strictly positive singular values. 
According to Proposition~2.1 of \cite{vandereycken2013low}, the tangent space of the manifold $\mathbb{M}$ at $\bmB$ is 
\begin{align}
	 T_{\bmB} \mathbb{M} = \Big\{\mDelta\in&\bbR^{K\times M}:\ \mDelta = \bmU\mM\bmV\trans + \mU_p\bmV\trans + \bmU\mV_p\trans,\nonumber \\ &\mM\in\bbR^{R\times R}, \mU_p\in\bbR^{K\times R}, 
	\mV_p\in\bbR^{M\times R},\mU_p\trans \bmU = \vzero,\mV_p\trans \bmV = \vzero\Big\}. \label{eqn:reduced:tangentspace}
\end{align}
Applying Theorem~\ref{theorem:main} requires us to quantify the local complexity of the tangent space. 
For the reduced model~\eqref{eqn:method:reduced} with penalty~\eqref{eqn:method:penaltyP}, the norm $\sQ_{\eta_1}$ defined in~\eqref{eqn:defineQnorm} has the explicit expression as $\sQ_{\eta_1}(\mDelta) = \| (\mI + \eta_1 \mGamma)^{1/2}\mDelta\|_F$.	
The next lemma presents the complexity upper bound of the local set~\eqref{eqn:generalNeighbor} in the tangent space.

\begin{lemma} \label{lemma:reduced:gamma}
	Consider the local neighbor of the tangent space 
	\begin{align}  \label{eqn:reduced:neighbor}
		\mathbb{N}(\bmB,\delta) = \big\{ \mDelta\in T_{\bmB}\mathbb{M}:\ \| (\mI + \eta_1 \mGamma)^{1/2}\mDelta\|_F\le \delta \big\}
	\end{align}
	for some $\delta$. We have the following order of  complexity 
	\begin{align*}
		\gamma_{2}(\mathbb{N}(\bmB,\delta), \|\cdot\|_F) \lesssim R^{1/2}\{K^{1/2}\wedge\eta_1^{-1/(4d+4q)} +(M-R)^{1/2}\} \delta .
	\end{align*}
\end{lemma}

The above lemma implies the critical radius satisfies
\begin{equation} \label{eqn:reduced:criticalradius}
	\widehat{\delta}_N \lesssim \frac{R^{1/2}\{K^{1/2}\wedge\eta_1^{-1/(4d+4q)} +(M-R)^{1/2}\} }{N^{1/2}} .
\end{equation}
It remains to check the second fundamental form of the fixed rank manifold satisfies Condition~\ref{ass:manifold}. 

\begin{lemma} \label{lemma:reduced:second}
	The second fundamental form $\rm II(\cdot,\cdot)$ for the rank-$R$ manifold at $\bmB$ is 
	\begin{align} \label{eqn:reduced:secondFund}
		\rm {II}(\mDelta_1,\mDelta_2) =  \mathrm{P}_{\bmB}^{\perp}
		(\mDelta_1 \bmB^+\mDelta_2+ \mDelta_2 \bmB^+\mDelta_1),
	\end{align}
	where $\bmB^+ = \bmV \bmD^{-1} \bmU\trans$ is the generalized inverse, and
	$ \mathrm{P}_{\bmB}^{\perp}(\cdot)$ is the projection onto the normal space $N_{\bmB}\mathbb{M} = ( T_{\bmB}\mathbb{M})^{\perp}$. Condition~\ref{ass:manifold} is satisfied if the $R$-th singular value of $\bmB$ is bounded away from zero and $\eta_1K^{2(d-\nu)_+}$ is bounded from above.
\end{lemma}

The  result~\eqref{eqn:reduced:secondFund} can be derived from the adjoint relation \citep[Equation (8.4) of][]{Lee2018} between the second fundamental form and the Weingarten map. 
The Weingarten map of the fixed-rank manifold has been developed in~\cite{absil2013extrinsic}. The above discussion leads to the following theorem. 

\begin{theorem} \label{thm:reducedmodel}	
	For the reduced multi-task regression~\eqref{eqn:method:reduced} with fixed $M$, 
	suppose Conditions~\ref{ass:beta}--\ref{ass:lossfunConvex} hold and the $R$-th singular value of $\bmB$ is bounded away from zero. Then, with $\tau = \nu \wedge (\mathfrak{o}+1) + \{ q \wedge (\mathfrak{o}+1)\}/2$, we have the following upper bound of the convergence rate
	\begin{align}
		&	\frac{1}{M}\sum_{m=1}^{M} \Big\{\|\hat{\beta}_m-\beta_0\|_{X} + \eta^{1/2}_1\|\hat{\beta}_m\|_{\Gamma}\Big\}	 \nonumber\\
	    & \qquad \qquad =  O_p\Big(
		\frac{R^{1/2}\{K^{1/2}\wedge\eta_1^{-{1}/{(4d+4q)}} +(M-R)^{1/2}\}}{M^{1/2}N^{1/2}}	+ K^{-\tau} \nonumber\\
		& \qquad \qquad \qquad \qquad \qquad + \eta_1^{1/2}  K^{(d-\nu)_{+}} + \Big\{\frac{1}{M}\sum_{r=R+1}^{M}\bar{\sigma}_{0r}^2\Big\}^{1/2}\Big),
		\label{eqn:thm:reducedbound}
	\end{align}
provided the right hand side of~\eqref{eqn:thm:reducedbound} converges to zero as $N\to\infty$.
\end{theorem}

%

The bound on the right hand side of~\eqref{eqn:thm:reducedbound} is a direct consequence of plugging~\eqref{lemma:reduced:approxerror1}, \eqref{lemma:reduced:approxerror2}, and~\eqref{eqn:reduced:criticalradius} into~\eqref{eqn:mainbound}.  
From Theorem~\ref{thm:reducedmodel}, we can derive the rates of convergence of the penalized estimator according to different configurations of the parameters $\eta_1$ and $K$. 
The next two corollaries assumes the last term (the rank-$R$ constraint error) in~\eqref{eqn:thm:reducedbound} is negligible.
Corollary~\ref{corollary:reduced:rate1} addresses the case where the penalty derivative order $d$ is no greater than the true smoothness order $\nu$ of the slope functions (i.e., $d\le \nu$), while
Corollary~\ref{corollary:reduced:rate2} presents the result for $d>\nu$. 

\begin{corollary} \label{corollary:reduced:rate1}
	Under the same assumptions of Theorem~\ref{thm:reducedmodel}, consider the case when penalty derivative order is smaller or equal to the smoothness order of the slope functions (i.e., $d \le \nu$). Define $\iota = q+d$, and suppose the rank-$R$ approximation error ($\sum_{r=R+1}^{M}\bar{\sigma}_{0r}^2/M$) is negligible.  Then,\\
	(i)  we have the rate of convergence
	$$
	\frac{1}{M}\sum_{m=1}^{M} \Big\{\|\hat{\beta}_m-\beta_0\|_{X} + \eta^{1/2}_1\|\hat{\beta}_m\|_{\Gamma}\Big\} = O_p\Big(	(MN/R)^{-\tau/(2\tau + 1) }\Big),
	$$
	when $\eta_1 \lesssim (MN/R)^{-2(\iota\vee \tau)/(2\tau + 1)} $ and $K \asymp (MN/R)^{1/(2\tau + 1)}$;\\
	(ii)  we have the rate of convergence
	$$
	\frac{1}{M}\sum_{m=1}^{M} \Big\{\|\hat{\beta}_m-\beta_0\|_{X} + \eta^{1/2}_1\|\hat{\beta}_m\|_{\Gamma}\Big\} = O_p\Big(		(MN/R)^{-\iota/(2\iota+1)}\Big),
	$$
	when $\eta_1 \asymp (MN/R)^{-2\iota/(2\iota+1)}$ and $K^{(\iota\wedge\tau)} \gtrsim (MN/R)^{\iota / (2\iota + 1)}$. 
\end{corollary}

Corollary~\ref{corollary:reduced:rate1} has two subcases. Conclusion~(i) corresponds to the asymptotic behavior of the regression spline estimator, where the roughness penalty controlled by $\eta_1$ is relatively weak and the number $K$ of knots is tuned to be optimal. 
On the other hand, Conclusion~(ii) corresponds to the asymptotic behavior of the smoothing spline estimator, where the number $K$ of knots diverges fast to infinity and the penalty parameter $\eta_1$ is tuned to be optimal.

It is worth to mention that, when $M=R=1$,  Corollary~\ref{corollary:reduced:rate1} reduces to the convergence result for the single-task functional linear regression. In this case, \cite{yuan2010reproducing} have shown the optimal rate of converges is  
$N^{-(q+\nu)/\{2(q+\nu)+1\}}$ under the setting of RKHS. This rate can be achieved in Conclusion~(ii) of Corollary~\ref{corollary:reduced:rate1} by setting $d=\nu$. On the other hand, as $\tau= \nu \wedge (\mathfrak{o}+1) + \{ q \wedge (\mathfrak{o}+1)\}/2< q+\nu$, the rate $N^{-\tau/(2\tau+1)}$ obtained in Conclusion~(i) is slower. 
This is because, for the regression spline estimator, the order of spline approximation error $K^{-\tau}$ is relatively larger.

Corollary~\ref{corollary:reduced:rate2} below focuses on the case of $d>\nu$. 
It also has two subcases behave either like using regression splines or smoothing splines, respectively.

\begin{corollary} \label{corollary:reduced:rate2}
	Under the same conditions of Theorem~\ref{thm:reducedmodel}, consider the case when penalty order is larger than the smoothness order (i.e., $d > \nu$). 
	Define $\iota = q+d$ and $\kappa = \tau+d-\nu$, then: \\
	(i)  we have the rate of convergence
	$$
	\frac{1}{M}\sum_{m=1}^{M} \Big\{\|\hat{\beta}_m-\beta_0\|_{X} + \eta^{1/2}_1\|\hat{\beta}_m\|_{\Gamma}\Big\} = O_p\Big( (MN/R)^{-\tau/(2\tau + 1) }\Big),
	$$
	when $\eta_1 \lesssim (MN/R)^{-2\iota/(2\tau + 1)}$ and $K \asymp (MN/R)^{1/(2\tau + 1)}$;\\
	(ii)  we have the rate of convergence
	$$
	\frac{1}{M}\sum_{m=1}^{M} \Big\{\|\hat{\beta}_m-\beta_0\|_{X} + \eta^{1/2}_1\|\hat{\beta}_m\|_{\Gamma}\Big\} = O_p\Big(		(MN/R)^{-\iota\tau/(\kappa +2\iota\tau)}\Big),
	$$
	when $\eta_1 \asymp (MN/R)^{-2\iota \kappa/(\kappa+2\iota\tau)}$ and $K\asymp  (MN/R)^{\iota/(\kappa+2\iota\tau)} $.
\end{corollary}

The above rates of convergence are summarized in Table~\ref{tbl:reduceRates}.
The rows are divided into two groups depending on $d \le \nu$ (Corollary~\ref{corollary:graph:rate1}) and $d > \nu$ (Corollary~\ref{corollary:graph:rate2}). 
Each group has its own cases corresponding to asymptotic behaviors like the regression spline estimator and like the smoothing spline estimator, depending on the values of $\eta_1$ and $K$.

\section{Application II: The Graph Regularized Multi-task Model}
\label{sec:graphrate}
In this section, Theorem~\ref{theorem:main} is applied to derive the upper bound of the convergence rate for the graph regularized model~\eqref{eqn:method:graph} as both $N$ and $M$ diverge to infinity. 
The penalty in~\eqref{eqn:method:penaltyGraph} involves a graph Laplacian matrix $\mOmega$. 
The convergence of the graph Laplacian has been studied in a large amount of works, e.g., \cite{hein2005geometrical,hein2007graph,belkin2007convergence,von2008consistency}. 
Generally speaking, as $M$ increases to infinity, it is known that the graph Laplacian converges to the Laplace-Beltrami operator over a manifold. 

\subsection{Convergence of Laplacian Matrix}
\label{sec:graphrate:review}
For the model~\eqref{eqn:method:graph},	suppose the	auxiliary variables $\vs_1,\cdots, \vs_M \in \bbR^{s}$ are concentrated on  a manifold $\sS$ ($\subset \bbR^{s}$), which is a compact Riemannian submanifold of $\bbR^{s}$ with intrinsic dimension $\mu$ and without boundary.
Besides, $\sS$ is assumed to satisfy certain regularity conditions as imposed in \cite{trillos2020error}, and it is endowed with a metric $\langle\cdot,\cdot\rangle_{\vs}$ and a corresponding Riemannian volume form $\intd V_{\vs}$. 
Suppose the	variables $\vs_1,\cdots, \vs_M$ are random sampled in accordance with a density $p$ defined over $\sS$. The density is Lipschitz continuous with Lipschitz constant $L_p$, and it is bounded from below and above ($C_p > p(\vs) > 1/C_p$ for some constant $C_p>0$).

Consider a smooth function $f(\vs): \sS\to \mathbb{R}$,  and let $\vf = \big(f(\vs_1),\cdots, f(\vs_M)\big)\trans$ be a $M$-dimensional vector containing the function evaluations at $\vs_1,\cdots, \vs_M$. The graph Laplacian matrix $\mOmega$ corresponds to the discrete Dirichlet form 
\begin{align*}
	b(f):=(1/M)\sum_{v,v'=1}^{M} w_{vv'} \big\{f(\vs_{v}) - f(\vs_{v'})\big\}^2 = \vf\trans\mOmega\vf/M,
\end{align*}	
where $w_{vv'}$ is defined as in \eqref{eqn:edgeweights} with a bandwidth parameter $h$.
Denote $\lambda_{k}^{\uparrow}(\mOmega)$ as the $k$-th eigenvalue of $\mOmega$ in increasing order. 
We can check $\lambda_{k}^{\uparrow}(\mOmega)$ is also the $k$-th eigenvalue of $b(f)$ with respect to the normalized Euclidean norm $\|\vf\|_{M} = \big\{\sum_{m=1}^M f^2(\vs_m)/M\big\}^{1/2}$. 
The continuous counterpart of the Dirichlet form has the expression 
\begin{align*} 
	D(f):=\int_{\sS} \big\langle \nabla_{\vs} f(\vs), \nabla_{\vs} f(\vs) \big\rangle_{\vs}\cdot p^2(\vs) \intd V_{\vs} =			\int_{\sS} f(\vs) \big\{\Delta_{\vs}f(\vs) \big\}
	p(\vs) \intd V_{\vs},
\end{align*}
where the Laplace-Beltrami operator $\Delta_{\vs}f:= -({1}/{p})\mathrm{div}(p^2\nabla_{\vs} f)$, with manifold divergence $\mathrm{div}(\cdot)$ and gradient $\nabla_{\vs}$.
We let $\lambda_{k}^{\uparrow}(\Delta_{\vs})$ denote the $k$-th eigenvalue of the quadratic form $D(f)$ with respect to the weighted $L_2$ norm $\|f\|_{\mathbb{L}_2(\sS,p)} = \big\{\int_{\sS} f^2(\vs) p(\vs)\intd V_{\vs}\big\}^{1/2}$.

\begin{figure}
	\centering
	\includegraphics[width=0.32\textwidth]{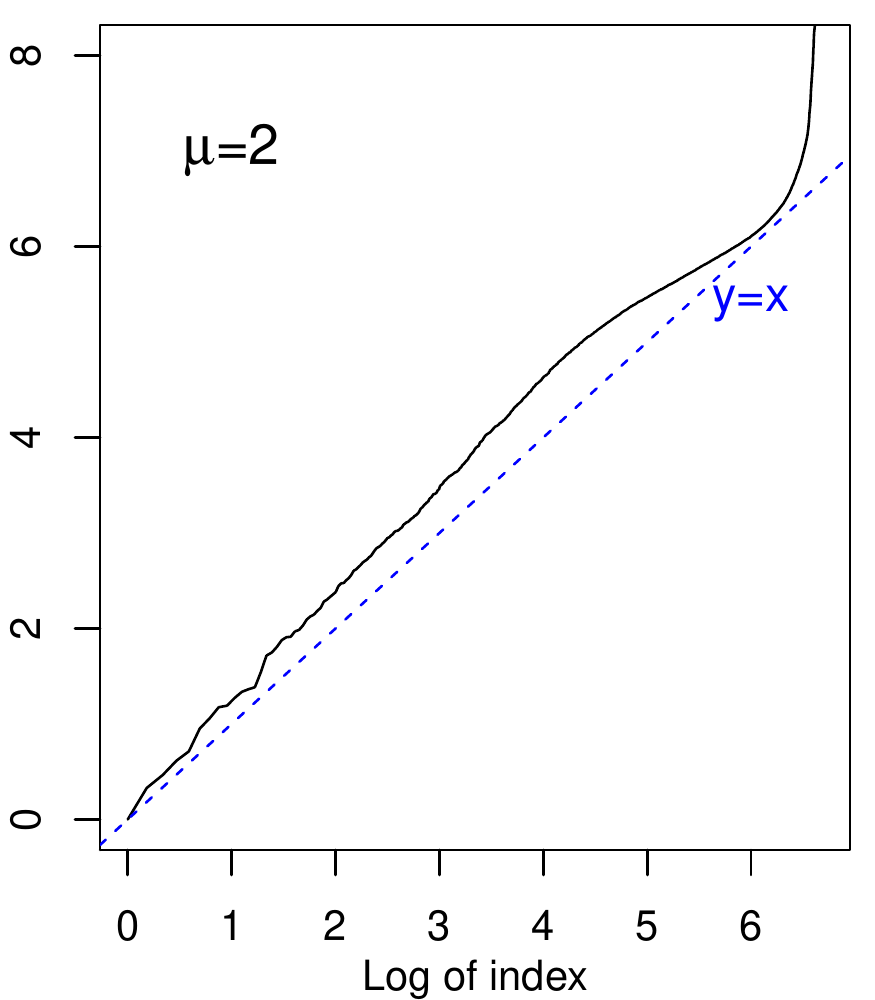}
	\includegraphics[width=0.32\textwidth]{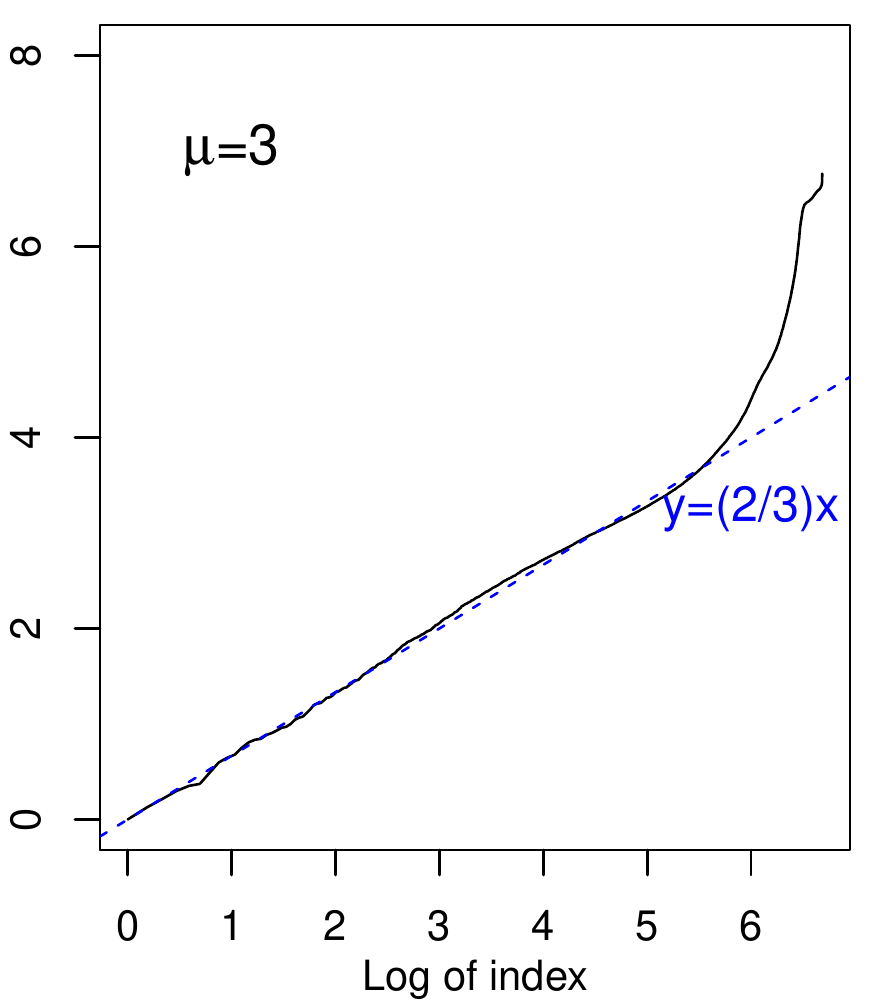}
	\includegraphics[width=0.32\textwidth]{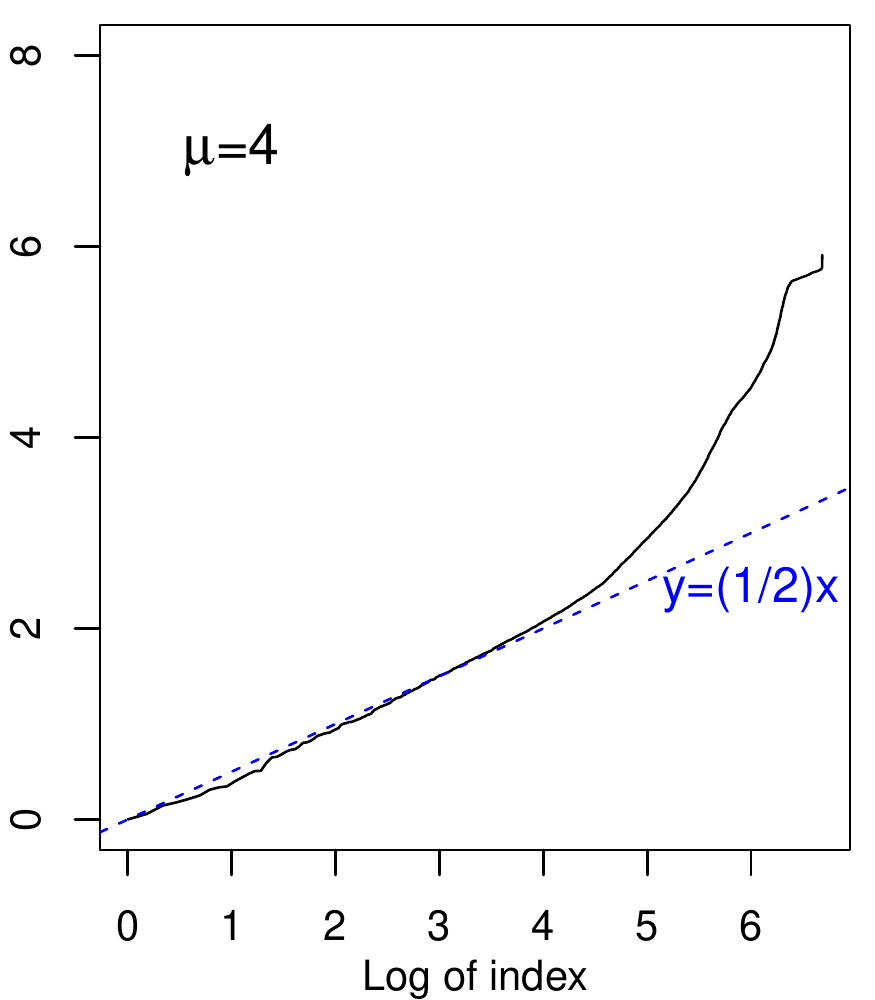}
	\caption{Plots of graph Laplacian eigenvalues. The graph Laplacian $\mOmega$ is computed from a collection of points of size $4000$, which are uniformly sampled from a sphere witch intrinsic dimension $\mu$. The three panels from left to right correspond to $\mu=2$, $3$, and $4$, respectively. 
	The first four smallest eigenvalues are excluded. The black solid curves represent the logarithm of the eigenvalues $\log\{\lambda_m^{\uparrow}(\mOmega)\}-\log\{\lambda_5^{\uparrow}(\mOmega)\}$ versus the logarithm of the index $\log(m/5)$, for $m=5,6,\cdots$.
	The blue dashed lines ($y=(2/\mu) x$) indicate the theoretical growth rate for the eigenvalues of the Laplace-Beltrami operator. \label{fig:lapeigv}}
\end{figure}

For a given first-order smooth $f$ and some properly chosen bandwidth $h$ in \eqref{eqn:edgeweights}, the relative magnitude between $b(f)$ and $D(f)$ can be bounded as $M\to\infty$. 
The bound  leads to the eigenvalue convergence in \cite{trillos2020error}. Corollary~1 of~\cite{trillos2020error} implies, for each given $m$, it holds that
\begin{equation} \label{eqn:trillosEigenConv}
	|\lambda_m^{\uparrow}(\mOmega) - \lambda_m^{\uparrow}(\Delta_{\vs})|\, /\, \lambda_m^{\uparrow}(\Delta_{\vs}) = o_p(1),
\end{equation}
when the kernel bandwidth $h$ satisfies that $h\to 0$ and $hM^{1/\mu}/\log(M)^{\zeta_\mu} \to\infty$, with $\zeta_\mu=3/4$ if $\mu=2$ and $\zeta_\mu=1/\mu$ if $\mu\ge 3$. According to the Weyl's law (see Equation (2.8) of~\cite{grigor2006heat}), the Laplacian-Beltrami operator $\Delta_{\vs}$ has discrete non-negative spectrum satisfying $\lambda_{m}^{\uparrow}\big( \Delta_{\vs} \big) \asymp m^{2/\mu}$. 
The eigenvalue convergence~\eqref{eqn:trillosEigenConv} implies $\lambda_m^{\uparrow}(\mOmega)\asymp m^{2/\mu}$ for a given $m$ and a large enough $M$. Although the convergence in~\eqref{eqn:trillosEigenConv} is not uniform in $m$, numerical results show $\lambda_m^{\uparrow}(\mOmega)$ can be lower bounded by $C m^{2/\mu}$ for some constant $C$ in lots of cases. See Figure~\ref{fig:lapeigv} for example, where 4000 points $\vs_1, \cdots, \vs_{4000}$ are sampled over a standard Euclidean sphere with intrinsic dimension $\mu$ ($=2,3,4$). The Laplacian matrix $\mOmega$ is computed by the procedure in Section~\ref{sec:method:graph} with kernel $G(u) \propto \exp(-u) I(u\in(0,1))$. 
The first four smallest eigenvalues of $\mOmega$ are excluded from the plot. The horizontal axis shows the logarithm $\log(m/5)$ of the index $m$. The black solid curves are the logarithm of the eigenvalues $\log\{\lambda_m^{\uparrow}(\mOmega)\} - \log\{\lambda_5^{\uparrow}(\mOmega)\}$ for $m=5,6,\cdots$. 
The blue dashed lines ($y = (2/\mu) x$) indicate the theoretical growth rate for the eigenvalues of the Laplace-Beltrami operator. From the figure, we can see $C m^{2/\mu}$ is a reasonable lower bound for $\lambda_m^{\uparrow}(\mOmega)$ in these empirical examples. 
In the following, we impose such lower bound assumption on the growing order of $\lambda_k^{\uparrow}(\mOmega)$.

\begin{condition}  \label{ass:laplacianEigenvalue}
	The eigenvalues of the graph Laplacian matrix grow at least with the order of $2/\mu$, i.e., $\lambda_m^{\uparrow}(\mOmega) \gtrsim  m^{2/\mu}$.
\end{condition}

\begin{remark} \label{remark:fixedgraph}
In this work, the graph $\mathcal{G}$ is constructed from the auxiliary random covariates $\vs_1,\cdots, \vs_M$ and is embedded in the manifold $\mathcal{S}$. 
We may also consider the setting where a large graph $\mathcal{G}$ is pre-given and is not embedded in any manifold as in \cite{kirichenko2017estimating}.
The work of \cite{kirichenko2017estimating} also adopts the  eigenvalue lower bound $\lambda_m^{\uparrow}(\mOmega) \gtrsim  m^{2/\mu}$.  They show that the growing order of eigenvalues is satisfied by various graph types, such as grid, discrete tori, lollipop graph, Watts-Strogatz ``small world'' graph, etc. 
See the detailed conditions and the discussions in \cite{kirichenko2017estimating}. 
This means our convergence rates developed in this section can also be verified in a similar fixed graph setting.
\end{remark}

\subsection{Convergence Rates}

Condition~\ref{ass:beta} has assumed certain smoothness for each true slope function $\beta_{0m}$. 
For the graph regularized model, a stronger condition characterizing the smoothness of the slope functions between distinct tasks is required. Suppose the true slope function for each task is a slice of a smooth surface $\beta_0(t,\vs)$ defined over $\sT\times \sS$. 
In other words, $\beta_{0m}(t) \equiv \beta_0(t,\vs_m)$ for the $m$-th task with auxiliary variable $\vs_m$. 
We denote $\beta^{(k)}_0(t,\vs) = (\partial^k / \partial t^k) {\beta_0(t,\vs)}$ as the $k$-th partial derivative with respect to $t$, and $\nabla_{\vs}\beta_0(t,\vs)$ as the gradient with respect to $\vs$ over the manifold $\sS$. 
The next condition states that the function $\beta_0(t,\vs)$ is continuously differentiable in both $t$ and $\vs$.  

\begin{condition} \label{ass:betagraph}
	For $k=0,1,\cdots,\nu$, both the derivatives of the true slope surface $\beta_0^{(k)}(t,\vs)$ and the manifold gradients $\nabla_{\vs} \beta^{(k)}_0(t,\vs)$ are continuous with respect to $t$ and $\vs$.
\end{condition}

Based on the additional Conditions~\ref{ass:laplacianEigenvalue} and \ref{ass:betagraph}, we now apply Theorem~\ref{theorem:main} to the model~\eqref{eqn:method:graph}. Recall the penalty in~\eqref{eqn:method:penaltyGraph} for the graph regularized model is
\begin{equation*} 
	\sP_{\veta}(\mB) = \eta_1 \tr\big(\mB\trans\mGamma\mB\big) +
	\eta_2 \tr\big(\mB\mOmega\mB\trans\widehat{\mSigma}\big) + \eta_1\eta_2
	\tr\big(\mB\trans\mGamma\mB\mOmega\big).
\end{equation*}
The pooled covariance matrix $\widehat{\mSigma}$ in the above is expected to converge to $\mI$ as $N\to \infty$, since we have used a simultaneous diagonalization technique to the spline basis (see Proposition~\ref{prop:splineStruct}). 
The next lemma rigorously establishes the limit of $\sP_{\veta}(\mB) $ based on Proposition~\ref{proposition:empiricalnorm}. 

\begin{lemma} \label{lemma:graph:Qbound} 
	Suppose $\{K\wedge \eta_{1}^{-1/(2q+2d)}\}/(MN)\to 0$ as $N\to \infty$. The scaled penalty $\sP_{\veta}(\mB)/M$ converges in probability to
	\begin{align}
		&	\big\{\eta_1 \tr\big(\mB\trans\mGamma\mB\big)  + \eta_2\tr\big(\mB\mOmega\mB\trans\big) + \eta_1 \eta_2
		\tr\big(\mB\trans\mGamma\mB\mOmega\big) \big\}/M \nonumber \\
		& \qquad =  \frac{\eta_1}{M} \sum_{m=1}^M
		\int_{\sT} \big\{\beta_m^{(d)}(t)  \big\}^2 \intd t  \label{eqn:aux:graphpenalty}\\
		&\qquad	\qquad \qquad +	\frac{\eta_2}{M} \sum_{v,v'=1}^M w_{vv'}
		\Big[\Expect_{x_{nm}}\langle x_{nm},\beta_{v} -  \beta_{v'} \rangle^2  
		+ \eta_1 \int \big\{\beta_{v}^{(d)}(t) -  \beta_{v'}^{(d)}(t) \big\}^2 \intd t
		\Big].  \nonumber
	\end{align} 
	Moreover, with probability at least
	$1-\exp\big(-  K\wedge \eta^{-\frac{1}{2q+2d}}_1 \big)$, it holds that
	\begin{align}  \label{eqn:graph:normbound}
		(1/2) \| (\mI + \eta_1 \mGamma)^{1/2} \mB (\mI + \eta_2 \mOmega)^{1/2}\|_F \le
		\sQ_{\veta}(\mB ) \le 2
		\| (\mI + \eta_1 \mGamma)^{1/2} \mB (\mI + \eta_2 \mOmega)^{1/2}\|_F.
	\end{align}		
\end{lemma}

In the above, the first term in~\eqref{eqn:aux:graphpenalty} measures the roughness of slope functions. 
Meanwhile, the second and third terms in~\eqref{eqn:aux:graphpenalty} measure the similarity of slope functions based on their differences of the predictive error and the $d$-th derivatives, respectively. 
The second conclusion~\eqref{eqn:graph:normbound} implies the norm in terms of $\sQ_{\veta}(\mB)$ is equivalent to the that of $\| (\mI + \eta_1 \mGamma)^{1/2} \mB (\mI + \eta_2 \mOmega)^{1/2}\|_F$ with high probability.

We next quantify the approximation error term $\sE(\mathbb{S}_K)$ defined by~\eqref{eqn:perror:spline}.
In particular, the bound of the penalty term $\sP_{\veta} (\bar{\vbeta}_0)$ in $\sE(\mathbb{S}_K)$ is derived, which is achieved by constructing spline approximation of the true slope function $\beta_{0m}(t) =\beta_0(t,\vs_m)$ at each $\vs_m$.
When the true slope surface $\beta_0(t,\vs)$ is sufficiently smooth and satisfies Condition~\ref{ass:betagraph}, it is reasonable to expect that the magnitude of the penalty $\sP_{\veta} (\bar{\vbeta}_0)$ can be controlled.

\begin{lemma} \label{lemma:graph:bias}
	The spline approximation error $\sE(\mathbb{S}_K)$ for the graph regularized model satisfies
	$$
	\big\{\sE(\mathbb{S}_K) \big\}^2 =O_p\Big( M\big\{K^{-2\tau} +
	\eta_1  K^{2(d-\nu)_{+}}  + \eta_2 + \eta_1\eta_2K^{2(d-\nu)_{+}}\big\}\Big).
	$$
	when the kernel bandwidth $h$ satisfies $h\to 0$ and $h M^{1/\mu}/\log(M)^{\zeta_\mu} \to\infty$, with $\zeta_\mu=3/4$ if $\mu=2$ and $\zeta_\mu=1/\mu$ if $\mu\ge 3$.
\end{lemma}
In the above,  the two terms $K^{-2\tau} +
\eta_1  K^{2(d-\nu)_{+}}$ are due to spline approximation error and the roughness penalty, as in Proposition~\ref{proposition:general:splineerror}. The additional two terms $\eta_2 + \eta_1\eta_2K^{2(d-\nu)_{+}}$ are attributed to the graph regularization.

Because $\mathbb{M}=\bbR^{K\times M}$ for the graph regularized model~\eqref{eqn:method:graph}, the solutions $\bmB_0$ to~\eqref{def:popuestimate} and $\bmB$ to~\eqref{def:popuestimate:constraint} are identical.
We therefore have null manifold approximation error, i.e., $\sE(\mathbb{M})=0$ for~\eqref{eqn:perror:manifold}. In this case, the tangent space at $\bmB$ is identical to the full Euclidean space $\bbR^{K\times M}$. 
As for the local neighbor set~\eqref{eqn:generalNeighbor}, we are considering
$$\mathbb{N}(\bmB,\delta) = \big\{ \mDelta\in\bbR^{K\times M}:\ 	\sQ_{\veta}(\mDelta )\le \delta \big\}.$$
Based on Lemma~\ref{lemma:graph:Qbound}, the local neighbor $\mathbb{N}(\bmB,\delta)$ has the same complexity level as the set
$$\mathbb{N}_1(\bmB,\delta) = \big\{ \mDelta\in\bbR^{K\times M}:\ 	\| (\mI + \eta_1 \mGamma)^{1/2} \mDelta (\mI + \eta_2 \mOmega)^{1/2}\|_F\le \delta \big\}.$$
with high probability. In the above, $\mGamma=\mathrm{diag}(\gamma_1,\cdots, \gamma_K)$ is a diagonal matrix specified in Proposition~\ref{prop:splineStruct}. 
Given Condition~\ref{ass:laplacianEigenvalue}, we can also diagonalize $\mOmega$ to a matrix containing its eigenvalues $\lambda_m^{\uparrow}(\mOmega) \gtrsim  m^{2/\mu}$. These lead to the complexity upper bound in the next lemma.


\begin{lemma} \label{lemma:graph:gamma}
	Consider the local neighbor $\mathbb{N}(\bmB,\delta) = \big\{ \mDelta\in\bbR^{K\times M}:\ \sQ_{\veta}(\mB )\le \delta \big\}$
	for some $\delta >0$. Under Condition~\ref{ass:laplacianEigenvalue}, we have the following complexity upper bound 
	\begin{align*}
		\gamma_{2}(\mathbb{N}(\bmB,\delta), \|\cdot\|_F)	=O_p\Big(
		(M^{1/2}\wedge \eta_2^{-\mu/4})\times\{K^{1/2}\wedge\eta_1^{-{1}/{(4d+4q)}} \} \, \delta\Big).
	\end{align*}
\end{lemma}

The above lemma implies the critical radius $\widehat{\delta}_N$ satisfies
\begin{equation} \label{eqn:graph:criticalradius}
	\widehat{\delta}_N =O_p\Big( \frac{ (M^{1/2}\wedge \eta_2^{-\mu/4})\times \{K^{1/2}\wedge\eta_1^{-1/(4d+4q)} \} }{N^{1/2}} \Big).
\end{equation}
Combining the above discussions, we use the unified result in Theorem~\ref{theorem:main} to get the upper bound of the convergence rate for the penalized estimator of the graph regularized model~\eqref{eqn:method:graph}.

\begin{theorem} \label{thm:mainGraph}
	Suppose Conditions~\ref{ass:beta}--\ref{ass:betagraph} holds. In addition, the kernel bandwidth $h$ is chosen such that $h\to 0$ and $h M^{1/\mu}/\log(M)^{\zeta_\mu} \to\infty$, with $\zeta_\mu=3/4$ if $\mu=2$ and $\zeta_\mu=1/\mu$ if $\mu\ge 3$. 
    The rate of convergence for the penalized estimator of the the graph regularized model~\eqref{eqn:method:graph} has an upper bound as
	\begin{align*}
		& (1/M) \bigg\{\sum_{m=1}^{M} \|\hat{\beta}_m-\beta_{0m}\|_{X} +  \sP_{\veta}^{1/2}(\hat{\vbeta})\bigg\}  \\
        & \qquad \qquad = O_p\Big(	
		\frac{(M^{1/2}\wedge \eta_2^{-\mu/4})\times\{K^{1/2}\wedge\eta_1^{-1/(4d+4q)} \}}{M^{1/2}N^{1/2}}	+ K^{-\tau} \\
		& \qquad \qquad \qquad \qquad \qquad + \eta_1^{1/2}  K^{(d-\nu)_{+}} +  \eta_2^{1/2} +  \eta_1^{1/2}\eta_2^{1/2}K^{(d-\nu)_{+}} \Big),
	\end{align*}
	with $\tau = \nu \wedge (\mathfrak{o}+1) + \{ q \wedge (\mathfrak{o}+1)\}/2$, provided the right hand side of the above converges to zero as $M,N,K\to \infty$, and $\eta_1,\eta_2\to 0$.
\end{theorem}

Distinct rates of convergence will appear as we vary the configuration of the related parameters for the graph regularized model. Corollary~\ref{corollary:graph:rate1} and Corollary~\ref{corollary:graph:rate2} below deal with the cases when the graph regularization is weak ($\eta_2 \lesssim M^{-2/\mu}$) and 	strong ($\eta_2 \gtrsim M^{-2/\mu}$), respectively. Only the cases of $d\le \nu$ are presented. The cases of $d>\nu$ can be analyzed similarly and are omitted.

\begin{corollary} \label{corollary:graph:rate1}
	(Weak graph regularization) Under the same conditions of Theorem~\ref{thm:mainGraph}, consider the case when the penalty derivative order is smaller or equal to the smoothness order (i.e., $d \le \nu$). 
	Denote $\iota = q+d$.
	Suppose the graph regularization is weak as $\eta_2 \lesssim M^{-2/\mu}$, then it holds that:	\\	
	(i)	when $\eta_1 \lesssim N^{-2(\iota \vee \tau)/(2\tau + 1)}$, $\eta_2 \lesssim M^{-2/\mu} \wedge N^{-2\tau/(2\tau + 1)}$, and $K \asymp N^{1/(2\tau + 1)}$, the rate of convergence satisfies
	$$
	(1/M) \bigg\{\sum_{m=1}^{M} \|\hat{\beta}_m-\beta_{0m}\|_{X} +  \sP_{\veta}^{1/2}(\hat{\vbeta})\bigg\} =O_p\Big(	 N^{-\tau/(2\tau + 1)}\Big);
	$$
	(ii) when $\eta_1 \asymp N^{-2\iota /(2\iota +1)}$, $\eta_2 \lesssim M^{-2/\mu} \wedge N^{-2\iota/(2\iota+1)} $, $K^{( \iota \vee \tau)} \gtrsim N^{\iota/(2\iota+1)}$, the rate of convergence satisfies
	$$
	(1/M) \bigg\{ \sum_{m=1}^{M} \|\hat{\beta}_m-\beta_{0m}\|_{X} +  \sP_{\veta}^{1/2}(\hat{\vbeta}) \bigg\} =O_p\Big(	  N^{-\iota/(2\iota+1)}\Big). 
	$$
\end{corollary}

The above results mean that, when the graph regularization is weak ($\eta_2 \lesssim M^{-2/\mu}$), the estimation behaves as if each task is estimated independently. 
In particular, Conclusion~(i) of Corollary~\ref{corollary:graph:rate1} corresponds to the asymptotic behavior of the regression spline estimator. 
Conclusion~(ii) of Corollary~\ref{corollary:graph:rate1} reflects  the asymptotic behavior of the smoothing spline estimator. 

The next result shows the rate of convergence can be further improved when the graph regularization is strong ($\eta_2 \gtrsim M^{-2/\mu}$) and the number $M$ of task is relatively large.

\begin{corollary} \label{corollary:graph:rate2}
	(Strong graph regularization) Under the same conditions of Theorem~\ref{thm:mainGraph}, consider the case when the penalty derivative order is smaller or equal to the smoothness order (i.e., $d \le \nu$). Suppose the graph regularization is strong as $\eta_2 \gtrsim M^{-2/\mu}$. Consider the convergence  bound of the form
	\begin{equation} \label{corollary:graph:rate2:bound}
		(1/M) \bigg\{\sum_{m=1}^{M} \|\hat{\beta}_m-\beta_{0m}\|_{X} +  \sP_{\veta}^{1/2}(\hat{\vbeta})\bigg\} = O_p\Big(	 (MN)^{-r_{mn}}\Big).
	\end{equation}
    Then, the rate of convergence $r_{mn}$ can be identified in the following cases:	\\	
	(i) we have the convergence rate $r_{mn}= {\tau}/\{{\tau(2+\mu) + 1}\}$,
	when the tuning parameters are configured as $\eta_1 \lesssim (MN)^{-2(\iota\vee\tau)r_{mn}/\tau}$, $\eta_2 \asymp (MN)^{-2r_{mn}}$, and $K \asymp (MN)^{r_{mn} /\tau}$.\\
	(ii) we have the convergence rate $r_{mn}= {\iota}/\{{\iota(2+\mu) + 1} \}$,
	when the tuning parameters are configured as $\eta_1 \asymp (MN)^{-2r_{mn}}$, $\eta_2 \asymp (MN)^{-2r_{mn}}$, and $K \gtrsim (MN)^{r_{mn}/(\iota \wedge \tau)}$.
\end{corollary}

In Corollary~\ref{corollary:graph:rate2}, it is implicitly required that the optimal tuning	$\eta_2 \asymp (MN)^{-2r_{mn}}$ satisfies the strong graph regularization bound $\eta_2 \gtrsim M^{-2/\mu}$. 
This is equivalent to saying the number $M$ of tasks should be large enough such that $M\gtrsim N^{\mu r_{mn}/(1-\mu r_{mn})}$, and Corollary~\ref{corollary:graph:rate2} implies two different scenarios accordingly.
\begin{enumerate}
    \item When $M\asymp N^{\mu r_{mn}/(1-\mu r_{mn})}$, the convergence  bound \eqref{corollary:graph:rate2:bound} has the equivalent expression
	\begin{equation} \label{eqn:remark:graph01}
		(1/M) \bigg\{\sum_{m=1}^{M} \|\hat{\beta}_m-\beta_{0m}\|_{X} +  \sP_{\veta}^{1/2}(\hat{\vbeta})\bigg\}
		=O_p\Big(	 N^{-r_{mn}/(1-\mu r_{mn})} \Big).
	\end{equation}
	In Conclusion~(i) of Corollary~\ref{corollary:graph:rate2} with 
	$r_{mn}= \tau/\{\tau(2+\mu) + 1\}$, we can find that
	\begin{equation} \label{eqn:remark:graph02}
		\frac{r_{mn}}{1-\mu r_{mn}} = \frac{\tau}{2\tau+1};
	\end{equation}
	while in Conclusion~(ii) of Corollary~\ref{corollary:graph:rate2} with 
	$r_{mn}= \iota/\{\iota(2+\mu) + 1\}$, we can find that
	\begin{equation} \label{eqn:remark:graph03}
		\frac{r_{mn}}{1-\mu r_{mn}} = \frac{\iota}{2\iota+1}.
	\end{equation}
	The above means the results in Corollary~\ref{corollary:graph:rate2} reduce to those in
	Corollary~\ref{corollary:graph:rate1} when $ M\asymp N^{\mu r_{mn}/(1-\mu r_{mn})}$, i.e., the rates of convergence are the same as estimating the slope functions independently.
    \item When $M\gg N^{\mu r_{mn}/(1-\mu r_{mn})}$, the advantage of graph regularization kicks in. In this case, Corollary~\ref{corollary:graph:rate2} implies a much faster convergence rate than that of~\eqref{eqn:remark:graph01}--\eqref{eqn:remark:graph03}. The result reveals that, compared with estimating each slope function individually, the graph regularization can considerably improve the estimation when the number $M$ of slope functions grows fast enough.
\end{enumerate}

We conclude by considering the opposite case, where the number of tasks is not large enough $ M\ll N^{\mu r_{mn}/(1-\mu r_{mn})}$, but we still require strong graph regularization $\eta_2 \gtrsim M^{-2/\mu}$. In this case, the optimal choice of $\eta_2$ is the lower bound $\eta_2 \asymp M^{-2/\mu}$. This will cause the graph regularization bias to dominate the upper bound of the convergence rate. We will get a slower rate of convergence
\begin{align} \label{corollary:graph:rate3}
	(1/M) \bigg\{\sum_{m=1}^{M} \|\hat{\beta}_m-\beta_{0m}\|_{X} +  \sP_{\veta}^{1/2}(\hat{\vbeta})\bigg\}
	=O_p\Big( \eta_2^{1/2} \Big) =O_p\Big( M^{-{1}/{\mu}}\Big).
\end{align}

The above convergence rates in Corollary~\ref{corollary:graph:rate1}, Corollary~\ref{corollary:graph:rate2}, and~\eqref{corollary:graph:rate3} are summarized in Table~\ref{tbl:graphRates}.
The rows are divided into two groups: weak graph regularization ($\eta_2\lesssim M^{-2/\mu}$ in Corollary~\ref{corollary:graph:rate1}) and strong graph regularization ($\eta_2\gtrsim M^{-2/\mu}$ in Corollary~\ref{corollary:graph:rate2}). 
Each group has its own subcases corresponding to behaving either like the regression spline estimator or the smoothing spline estimator asymptotically, depending on the values of $\eta_1$ and $K$.
For example, in the strong graph regularization group (Settings (iii)--(vi)), the regression spline asymptotic scenario has two settings ((iii) and (iv)), according to the distinct orders of $M$.
The smoothing spline asymptotic scenario also has two similar settings ((v) and (vi)).

\section{Discussion}\label{sec:discussion}
In this work, we have proposed a general model with double regularization for multi-task functional linear regression models. Two folds of regularization include the matrix submanifold constraint and a penalization as the composite sum of  quadratic forms.  Through a comprehensive study of the properties of penalized splines in the scope of functional linear regression models, we show the composite quadratic penalty can induce a specific norm to quantify the manifold curvature and bound the complexity of the local set for the estimator using the technique of generic chaining.
All these tools lead to the unified upper bound of the convergence rate for the proposed general model. 
We further apply the unified upper bound to two specific multi-task functional linear  regression models with reduced rank and graph regularization, and figure out the convergence rates and the phase transition behaviors of the penalized spline estimators.  Although we have illustrated our framework using these two specific cases, our analysis tools can have implications for future research on similar problems.

Our theoretical study assumes the independent realizations of the functional covariates are fully observed. However, in practice, there are some situations that the functional data are sparsely observed. Extending the current results to the sparse-observed functional covariates is a potential future research topic. Moreover, the first regularization of the manifold constraint set is assumed to be known in this work. How to incorporate the manifold learning theory when the constraint set is unknown to our current results is also of interest, and needs further investigation.

\clearpage

\appendix

\textbf{\center \Large Appendix}

\section{Discussion on Condition~\ref{ass:Kx}}
\label{sec:covexamples}
In this work, we have employed Condition~\ref{ass:Kx}, which is a mild condition and includes many interesting examples. 
Condition~\ref{ass:Kx} is indeed a simplification of the Sacks-Ylvisaker condition, since it lists out the essential properties of the covariance function (i.e. the eigenvalue decay rate, the covariance function smoothness, and the possible existence of the null space).
In particular, the value of $p$ is not specified in this condition. The value of $p$ can be $0$ (i.e., $\mathbb{P}_{p}$ is an empty set) or a large number depending on the particular functional data of interest.  The covariance functions of many random processes satisfy Condition~\ref{ass:Kx}:

\begin{enumerate}
	\item[(1)] In the first example, $x(t)$ is a random process following the standard Brownian motion over the interval $\sT=[0,1]$. It then holds that $x(0)=0$ and its covariance function is $\mathcal{C}(s,t) = \min (s,t)$. This covariance function satisfies Condition~\ref{ass:Kx} with $q=1$ and $p=1$ such that $\mathbb{P}_p=\{x:\ x(t)\equiv c \text{ for some } c\in\mathbb{R}\}$ is the space of constant functions.  Because the standard Brownian motion starts with $x(0)=0$, it can be seen that this random process does not have variability in the space of constant functions. 
	
	\item[(2)] In the second example, we can  specify the covariance  as $\mathcal{C}(s,t) = a+ b\min (s,t)$ for some positive constants $a, b> 0$. In this case, the random process $x(t)$ has variability over the \textit{full} Sobolev space $\mathbb{L}_2^{1}(\sT)$, including the space of constant functions. In other words, we directly have
	$\mathbb{L}_2^1(\sT) = H(\mathcal{C})$, and $\mathbb{P}_p = \emptyset$ an empty set with $p=0$. 
	
	\item[(3)] Discussion similar to (1) and (2) above also applies to other random processes with $q=1$. For example, the Ornstein-Uhlenbeck process with covariance $\mathcal{C}(s, t)=c_1 \exp(- c_2|s - t|)$, where $c_1,c_2 > 0$ \citep[we have $H(\mathcal{C}) = \mathbb{L}_2^1 (\sT)$ according to][]{muller1996optimal};  %
	the sum $y_1(t)+y_2(1-t)$ of two independent Brownian motions whose covariance is $\mathcal{C}(s, t) = 1 - |s - t|$ \citep[we have $H(\mathcal{C}) = \mathbb{L}_2^1 (\sT)$ according to][]{muller1998spatial,ritter2000average}.
	For the Brownian Bridge with $\mathcal{C}(s, t) = \min(s, t) - st$, we have $H(\mathcal{C}) \oplus \mathbb{P}_2 = \mathbb{L}_2^{1}(\sT)$ with $\mathbb{P}_2=\{x:\ x(t)= c_1+c_2 t \text{ for some } c_1, c_2\in\mathbb{R}\}$  according to \cite{ritter1995}.
	
	\item[(4)] In the fourth example, we may consider the classical covariance kernel function with $q \geq 1$,
	$$\mathcal{C}(s, t)=\int_0^1 \frac{(s-u)_{+}^{q-1}(t-u)_{+}^{q-1}}{\{(q -1) !\}^2} \intd u.$$ 
	When the random process $x(\cdot)$ has the above $\mathcal{C}(s, t)$ as its covariance function, the random process does not have variability in $\mathbb{P}_{q}$, which consists of polynomials of order $q$. In fact, this covariance function corresponds to the $(q-1)$-fold integrated Brownian motion.
	
	\item[(5)] Continue the above example, when the covariance function is the Sobolev reproducing kernel
	$$\mathcal{C}(s, t)=
	\sum_{\ell=0}^{q-1}  c_{\ell}\frac{s^{\ell}}{\ell !} \frac{t^{\ell}}{\ell !} + c_q \int_0^1 \frac{(s-u)_{+}^{q-1}(t-u)_{+}^{q-1}}{\{(q -1) !\}^2} \intd u$$
	for some positive constants $c_0, c_1,\ldots, c_q (>0)$,  the corresponding random process $x(\cdot)$ has variability in the \textit{full}
	Sobolev space $\mathbb{L}_2^{q}(\sT)$ with an empty $\mathbb{P}_p = \emptyset$ (i.e., $p=0$). 
\end{enumerate}

\section{Technical Proofs of Section~\ref{sec:prelim}}

\subsection{Proof of Proposition~\ref{proposition:general:splineerror}}

	\begin{proof}
	Let $\tilde{\beta}_{0m}$ be the projection  of $\beta_{0m}$ onto the spline space $\mathbb{S}_K$. The projection is defined in the sense of $L_2$ norm such that
	$$
	\tilde{\beta}_{0m} = \argmin_{\beta\in\mathbb{S}_K} \| \beta - \beta_{0m}\|_{L_2}^2.
	$$
	The first order optimality condition implies that  the residual $\tilde{\beta}_{0m} - \beta_{0m}$ is orthogonal to the spline basis  $\vphi(t)$ of the spline space $\mathbb{S}_K$,
	\begin{equation} \label{proof:general:splineerror:ortho}
		\int \vphi(t) \cdot \{\tilde{\beta}_{0m}(t) - \beta_{0m}(t)\} \intd t  = \vzero.
	\end{equation}
	Let  $\bar{\mathcal{C}}$ be the optimal projection of $\mathcal{C}$ onto the tensor product spline space $\mathbb{S}_K^2$. Then,
	\begin{align*}
		\int \int \{\tilde{\beta}_{0m}(s) - \beta_{0m}(s)\} 
		\bar{\mathcal{C}}(s, t) \{\tilde{\beta}_{0m}(t) - \beta_{0m}(t) \} \intd s \intd t =0,
	\end{align*}
	due to the orthogonality~\eqref{proof:general:splineerror:ortho} and that $\bar{\mathcal{C}}(s, t)  = \vphi(s)\trans \mK\vphi(t)$ for some matrix $\mK$.
	It follows that
	\begin{align*}
		\|\tilde{\beta}_{0m} - \beta_{0m} \|_{X}^2 
		&=  \Expect \langle x_{nm}, \tilde{\beta}_{0m} -\beta_{0m} \rangle^2 \\
		&= \int \int \{\tilde{\beta}_{0m}(s) - \beta_{0m}(s)\} \mathcal{C}(s, t) \{\tilde{\beta}_{0m}(t) - \beta_{0m}(t)\} \intd s \intd t \\
		&= \int \int \{\tilde{\beta}_{0m}(s) - \beta_{0m}(s)\} 
		\{\mathcal{C}(s, t) - \bar{\mathcal{C}}(s, t) \}\{\tilde{\beta}_{0m}(t) - \beta_{0m}(t) \} \intd s \intd t \\
    	& \le \|\mathcal{C} - \bar{\mathcal{C}}\|_{L_2} \cdot \|\beta_{0m} - \tilde{\beta}_{0m}\|_{L_2}^2 \\
    	& \asymp K^{-2\tau},
	\end{align*}
	where $\tau = \nu \wedge (\mathfrak{o}+1) +[ q \wedge (\mathfrak{o}+1)]/2$.  The last inequality holds due to the following reasons. From Condition~\ref{ass:beta} and  Theorem~6.25 of \cite{schumaker2007spline}, we have 	$\|\beta_{0m} - \tilde{\beta}_{0m}\|_{L_2} \asymp K^{-\nu\wedge(\mathfrak{o}+1)}$.
	From (ii) of Condition~\ref{ass:Kx}  and Theorem~12.7 and Theorem~13.18 of \cite{schumaker2007spline}, we have $\|\mathcal{C} - \bar{\mathcal{C}}\|_{L_2} \asymp  K^{- q \wedge (\mathfrak{o}+1)}$.
	
		Meanwhile, also due to Theorem~6.25 of \cite{schumaker2007spline}, it holds that $\|\tilde{\beta}_{0m}\|_{\Gamma} \asymp K^{(d - \nu)_{+}}$. Combine the above results together to get the conclusion of Proposition~\ref{proposition:general:splineerror}.
\end{proof}

\subsection{Proof of Proposition~\ref{prop:splineStruct}}

\begin{proof}
	We can apply a two-step procedure to simultaneously diagonalize  $\|\cdot\|_{X}$ and  $\|\cdot\|_\Gamma$.\\
	
	\noindent	\textbf{Step One.} The first step of transformation relies on 
	Lemma A3 of \cite{claeskens2009asymptotic}. It is presented  as the next result for an ordinary B-spline basis $\widetilde{\vphi}$. Define two matrices
	$$
	\widetilde{\mN} = \int \widetilde{\vphi}(t) \widetilde{\vphi}\trans(t)  \intd t,\quad 
	\widetilde{\mGamma} = \int \widetilde{\vphi}^{(d)}(t) \big\{\widetilde{\vphi}^{(d)}(t) \big\}\trans \intd t.
	$$
	
	\begin{lemma}[Lemma A3 of \cite{claeskens2009asymptotic}] \label{lemma:splinestepone}
		Consider the eigen decomposition of $\widetilde{\mN}^{-1/2}\widetilde{\mGamma}\widetilde{\mN}^{-1/2} = \mV_1 \mW_1\mV_1\trans$, where $\mW_1$ is a diagonal matrix containing the eigenvalues in increasing order, and the matrix $\mV_1$ has the eigenvectors in its columns. The diagonal elements of the  matrix $\mW_1 = \diag(\omega_1, \cdots, \omega_{K})$ satisfy that
		\begin{equation} \label{eqn:gamma:diagonalscale}
			\omega_1=\cdots = \omega_d = 0, \quad 
			\omega_k \asymp (k - d)^{2d},\; \text{ for }\; k=d+1,\cdots, K.
		\end{equation}
	\end{lemma}

	Based on this lemma, we can construct an intermediate basis $\vvarphi(\cdot) = \mV_1\trans\widetilde{\mN}^{-1/2}\widetilde{\vphi}(\cdot)$. With this intermediate basis, a function 
	$\beta(u) = \widetilde{\vb}\trans\widetilde{\vphi}(\cdot)$ expressed by the B-spline basis $\widetilde{\vphi}$
	has an equivalent representation
	\begin{equation} \label{eqn:beta:equivalent1}
		\beta(u) = \ve\trans\vvarphi(\cdot) = \widetilde{\vb}\trans\widetilde{\vphi}(\cdot), 
	\end{equation}
	with an adjusted coefficient $\ve = \mV_1\trans\mN^{1/2}\widetilde{\vb}$. 
	In terms of the spline representation via $\vvarphi(\cdot)$, it holds that
	\begin{align}
		\int \{\beta^{(d)}(u)\}^2 \intd u
		&=\ve\trans \Big[\int \vvarphi^{(d)}(t) \big\{\vvarphi^{(d)}(t) \big\}\trans \intd t\Big] \ve\nonumber \\
		&=\ve\trans\big(\mV_1\trans\widetilde{\mN}^{-1/2}\widetilde{\mGamma}
		\widetilde{\mN}^{-1/2}	\mV_1 \big) \ve
		=\ve\trans \mW_1 \ve. \label{eqn:vvarphi:gamma}
	\end{align}
	At the same time, the intermediate  basis $\vvarphi(t)$ has become an orthonormal basis because
	\begin{equation} \label{eqn:vvarphi:orthonormal}
		\int \vvarphi(t)\vvarphi\trans(t)\intd t = \mV_1\trans\widetilde{\mN}^{-1/2}
		\Big\{\int \widetilde{\vphi}(t) 
		\widetilde{\vphi}\trans(t) \intd t\Big\}
		\widetilde{\mN}^{-1/2}\mV_1=  \mV_1\trans\mV_1=\mI.
	\end{equation}
	From the above~\eqref{eqn:vvarphi:gamma} and~\eqref{eqn:vvarphi:orthonormal}, we have actually
	that  the intermediate  basis $\vvarphi(t)$ simultaneously diagonalizes the pair of  norms $\|\cdot\|_{L_2}$ and  $\|\cdot\|_\Gamma$.
	
	\vspace{20pt}		
	
	\noindent \textbf{Step Two.}  A further step of transformation is then applied to the intermediate  basis $\vvarphi(t)$, such that it simultaneously diagonalizes $\|\cdot\|_{X}$ and  $\|\cdot\|_\Gamma$. For this purpose, define a matrix
	\begin{equation} \label{eqn:simuDiag:F}
		\mF := \Var\Big\{\int x_{nm}(t) \vvarphi(t) \intd t\Big\} = \int \int  \mathcal{C}(s,t) 
		\vvarphi(s)  \vvarphi\trans(t) \intd t \intd s,
	\end{equation}
	where $\mathcal{C}$ is the common covariance function for all functional predictors $x_{nm}$.

    Denote $\lambda_j^{\downarrow}\big(\mF\big)$ as the $j$-th eigenvalue of the matrix $\mF$ in decreasing order. It is upper bounded by
    \begin{equation} \label{eqn:simuDiag:lambda}
    	\lambda_j^{\downarrow}\big(\mF\big)\le \lambda_{0j},
    \end{equation}
    where $\lambda_{0j}$ is the $j$-th eigenvalue of $\mathcal{C}$ in decreasing order.
    To see this, suppose $\vu_j$ is the $j$-th eigenvector  of $\mF$ associated with the eigenvalue $\lambda_j^{\downarrow}(\mF)$. We can define a related function $\theta_j(\cdot) = \vu_j\trans\vvarphi(\cdot)$  for each $j$. 
    Due to the orthonormalilty of the intermediate basis $\vvarphi$ in~\eqref{eqn:vvarphi:orthonormal}, these functions $\theta_j(\cdot)$'s are also orthonormal to each other, i.e.,
    $\int \theta_j(t)\theta_{j'}(t) dt = \delta_{jj'}$. The first $j$ of them span a subspace $\mathcal{M}_j = \text{span}\{\theta_1,\theta_2,\cdots, \theta_j\}$ of $\mathbb{L}_2(\sT)$, where $\mathbb{L}_2(\sT)$ is the set of all squared integrable function over domain $\sT$. For  the first $j$  eigenvector  of $\mF$, they also span a $j$ dimensional subspace $\sU_j = \text{span}\{\vu_{1},\vu_{2},\cdots, \vu_{j}\}$ of the Euclidean space. By max-min principle of eigenvalues, it holds that
    \begin{align*}
	    \lambda_{0j} \stackrel{(i)}{=} & \max_{\Theta:\ \text{dim}(\Theta) =j\;} \min_{\theta\in\Theta} \int\int \mathcal{C}(s,t)\theta(s)\theta(t) \intd s \intd t \\
	    \ge &  \min_{\theta\in\mathbb{M}_j} \int\int \mathcal{C}(s,t)\theta(s)\theta(t) \intd s \intd t \\
    	= & \min_{\vu\in \sU_j} \vu\trans\mF\vu = \lambda_j^{\downarrow}(\mF).
    \end{align*}
    In first line~(i) of the above, $\Theta$ varies among any $j$ dimensional subspace of  $\mathbb{L}_2(\sT)$, and $\theta\in\Theta$ is an arbitrary function with unit norm.
	
    Suppose the eigen-decomposition of the matrix $\mF$ is $\mV_f\mW_f\mV_f\trans$, where $\mW_f$ is a diagonal matrix with decreasing eigenvalues. Recall when the covariance function $\mathcal{C}$ satisfies  Condition~\ref{ass:Kx},   the support of the probability measure of the random $x_{nm}$ may not be the full Sobolev space  $\mathbb{L}_2^q(\sT)$, but up to an additional finite dimensional subspace $\mathbb{P}_p$. This implies the (semi-)norm $\|\beta\|_{X}$ could possibly be zero for a non-zero $\beta\in\mathbb{S}_K$ in the spline space. Equivalently, the quadratic form $\|\beta\|_{X}^2$ can have zero eigenvalues with respect to  $\|\beta\|_{L_2}^2$ for $\beta \in\mathbb{S}_K$. Denote $\bar{p}$ as  the replicate number of the zero eigenvalues, and we can know $\bar{p}$ also equals the replicate number of zero eigenvalues of the matrix $\mathbf{F}$. Define
	\begin{align*}
		\widetilde{\mF} = \mV_f \big\{\mW_f + \vzero_{K-\bar{p}} \oplus
		\mathrm{diag}(\underbrace{K^{-2q},\cdots, K^{-2q}}_{\text{repeated }\bar{p} \text{ times }}) \big\} \mV_f\trans.
	\end{align*}
    Note $\widetilde{\mF}$ has the same eigenvectors and eigenvalues with $\mF$, except that the zero eigenvalues of $\mF$ is replaced by $K^{-2q}$.  Thereby, $\widetilde{\mF}$ is invertible.
	
	Take the eigen-decomposition of the matrix product 
	$\widetilde{\mF}^{-1/2}\mW_1\widetilde{\mF}^{-1/2}$ with $\mW_1$ defined in Lemma~\ref{lemma:splinestepone}. Suppose its decomposed matrix is 
	$$\widetilde{\mF}^{-1/2}\mW_1\widetilde{\mF}^{-1/2}=\mV_2\mGamma\mV_2\trans,$$
	where $\mGamma$ is a diagonal matrix with eigenvalues ordered increasingly. Then, apply one more step of  transformation to the intermediate basis $\vvarphi$ to get
	\begin{equation} \label{eqn:simuDiag:transform2}
		\vphi(\cdot) = \mV_2\trans\widetilde{\mF}^{-1/2} \vvarphi(\cdot).
	\end{equation}
	Given this transformation, we have another representation for the function
	$\beta(u) = \widetilde{\vb}\trans\widetilde{\vphi}(\cdot)$ in the spline space $\mathbb{S}_{K}$. Together with~\eqref{eqn:beta:equivalent1}, we can find that
	$$
	\beta(u) = \ve\trans\vvarphi(\cdot) = \vb\trans\vphi(\cdot),
	$$
	with the corresponding coefficients $\vb = \mV_2\trans\widetilde{\mF}^{1/2}\ve$.
	
	Now, we get $\vphi(\cdot)$ which has the desired  simultaneous diagonalization property.  To see this, we can verify that
	\begin{align*}
		\Var\left\{\int x_{nm}(t) \vphi(t) \intd t\right\} &=
		\mV_2\trans\widetilde{\mF}^{-1/2} \Var\left\{\int x_{nm}(t) \vvarphi(t) \intd t\right\}
		\widetilde{\mF}^{-1/2}\mV_2 \\
		&= \mV_2\trans (\mI_{K-\bar{p} }\oplus \vzero_{\bar{p}})\mV_2 = \mI_{K-\bar{p} }\oplus \vzero_{\bar{p}}.
	\end{align*}
	In the above, the second equality has employed~\eqref{eqn:simuDiag:F}.
	Then,  it holds for $\beta(u) =  \vb\trans\vphi(\cdot)$ that
	$$
	\|\beta\|_{X}^2 =  \Expect_{x^*}
	\langle x^*, \beta\rangle^2 =\vb\trans
	\left[\Var\left\{\int x_{nm}(t) \vphi(t) \intd t\right\}\right]\vb = \vb\trans(\mI_{K-\bar{p} }\oplus \vzero_{\bar{p}})\vb.
	$$
	Meanwhile, we have
	\begin{align}
		\int \vphi^{(d)} (t) \big\{ \vphi^{(d)} (t)\big\} \trans\intd t  &=
		\mV_2\trans\widetilde{\mF}^{-1/2} 
		\Big[\int \vvarphi^{(d)} (t) \big\{\vvarphi^{(d)} (t)\big\}\trans\intd t\Big]  \widetilde{\mF}^{-1/2}\mV_2 \nonumber\\
		&=  \mV_2\trans\widetilde{\mF}^{-1/2} \mW_1 \widetilde{\mF}^{-1/2}\mV_2 =
		\mGamma.
	\end{align}
	Therefore
	$$
	\|\beta\|_\Gamma^2 =  \int
	\big\{\beta^{(d)}(t)\big\}^2\intd t =
	\vb\trans \Big[\int \vphi^{(d)} (t) \big\{\vphi^{(d)} (t)\big\}\trans\intd t\Big] \vb=
	\vb\trans\mGamma\vb.
	$$
	
	Summarizing Step One and Step Two above,  we get a transformation of the original basis 
	$\vphi(\cdot) = \mQ\widetilde{\vphi}(\cdot)$,
	with $\mQ = \mV_2\trans\widetilde{\mF}^{-1/2}  \mV_1\trans\widetilde{\mN}^{-1/2}$. The transformed $\vphi(\cdot)$ simultaneously diagonalized $\|\cdot\|_X$ and $\|\cdot\|_\Gamma$. This verifies the first part of 
	Proposition~\ref{prop:splineStruct}.
	
	It remains to quantify the magnitude of the diagonal elements of $\mGamma$.
	For any $i,j,k$ satisfying $j+k\le i+1$ and $i=1,\cdots, K$, it holds that
	\begin{align*}
		\lambda_{K-i+1}^{\downarrow}(\mGamma) &=
		\lambda_{K-i+1}^{\downarrow}\big(\widetilde{\mF}^{-1/2}\mW_1\widetilde{\mF}^{-1/2}\big)\\
		&= \lambda_{K-i+1}^{\downarrow}\big(\mW_1\widetilde{\mF}^{-1}\big)  \\
	    & \stackrel{(i)}{\ge} \lambda_{K-j+1}^{\downarrow}\big(\mW_1\big)\times 
		\lambda_{K-k+1}^{\downarrow}\big(\widetilde{\mF}^{-1}\big) \\
		&= \lambda_{K-j+1}^{\downarrow}\big(\mW_1\big)\times 
		\big\{\lambda_{k}^{\downarrow}\big(\widetilde{\mF}\big)\big\}^{-1} \\
		&\stackrel{(ii)}{\gtrsim}  \lambda^{\downarrow}_{K-j+1}\big(\mW_1\big)\times  k^{2q} \,.
	\end{align*}
	The inequality~(i) uses Result 6.75(b) on Page 119 of \cite{seber2008matrix}.  The inequality (ii) uses~\eqref{eqn:simuDiag:lambda}
	and (i) of Condition~\ref{ass:Kx}. 	
	
	When $i>2d$, by exploiting Lemma~\ref{lemma:splinestepone}, we can set $j=k=\floor{i/2}$ to get
	$$
	    \gamma_{i}=	\lambda^{\downarrow}_{K-i+1}(\mGamma) \gtrsim
	    (\floor{i/2}-d)^{2d} \cdot  
	    \big(\floor{i/2}\big)^{2q} \gtrsim i^{2d+2q};
	$$
	when $i\le 2d$, it holds that $\gamma_i \ge 0$. 
\end{proof}

\subsection{Proof of Proposition~\ref{proposition:empiricalnorm}}

The result is established based on \cite{dirksen2015tail}. They consider the concentration bound for
$$
    Z(f) =  \Big|\frac{1}{N}\sum_{n=1}^{N} \big\{f^2(X_i) - \Expect f^2(X_i)  \big\} \Big|,
$$
for $f\in \mathcal{F}$ in a function class $\mathcal{F}$ and $X_1,\cdots, X_N$ are some i.i.d. random variables. 
\begin{lemma}[Corollary~5.7 of \cite{dirksen2015tail}] \label{lemma:squaredconcen}
	Suppose for $\sigma_1, G_1$ such that
	$$
	\sup_{f\in\mathcal{F}}  \frac{1}{N}\sum_{n=1}^{N}  
	\Expect \big\{f^2(X_i) - \Expect f^2(X_i)  \big\}^q \le \frac{q!}{2} \sigma_1^2 G_1^{q-2},
	$$
	for $q=2,3,\cdots$. Then, for any $u\ge 1$, 
	$$
	\Pr\bigg[ \sup_{f\in\mathcal{F}} Z(f) \ge C_1\Big\{
	\frac{1}{N} \gamma_{2}^2(\mathcal{F},d_{\psi_2} )+ 	
	\frac{\mathrm{diam}_{\psi_2}(\mathcal{F})}{\sqrt{N} } \gamma_{2}(\mathcal{F},d_{\psi_2}) \Big\} + 
	c_1\Big\{\sqrt{u}\frac{\sigma_1}{\sqrt{N}} + u\frac{G_1}{N}\Big\}\bigg] \le e^{-u},
	$$
	for some constants $C_1, c_1>0$.
\end{lemma}

\noindent \textbf{Proof of Proposition~\ref{proposition:empiricalnorm}}
    For $\epsilon\in (0,1)$, the expression
	\begin{equation}  \label{proof:empiricalTarget}
		(1+\epsilon )\big\{\|\beta\|_{X}^2   +\eta_1 \| \beta \|_{\Gamma}^2 \big\} \ge \|\beta\|_{Nm}^2 +\eta_1 \| \beta \|_{\Gamma}^2 \ge 
		(1-\epsilon )\big\{\|\beta\|_{X}^2   +\eta_1 \| \beta \|_{\Gamma}^2 \big\},
	\end{equation}
	is equivalent to
	\begin{equation} \label{proof:empiricalTarget:equiv}
		\big|\|\beta\|_{Nm}^2 - \|\beta\|_{X}^2 \big|
		\le   \epsilon \cdot \big\{ \|\beta\|_{X}^2   + \eta_1\| \beta \|_{\Gamma}^2 \big\}.
	\end{equation} 
	
	Under the complementary event  of~\eqref{proof:empiricalTarget:equiv}, there exists a $\beta \in\mathbb{S}_K$ such that
	\begin{equation} \label{proof:empirical:Am:complement}
		\big|	\|\beta\|_{Nm}^2 - \|\beta\|_{X}^2 \big|
		\ge  \epsilon \cdot \big\{\|\beta\|_{Nm}^2   +\eta_1 \| \beta\|_{\Gamma}^2\big\}.
	\end{equation} 
	Denote the set
	$$
	E=\{\beta\in\mathbb{S}_K:\ \|\beta\|_{X}^2   +\eta_1 \| \beta \|_{\Gamma}^2\le 1\}.
	$$
	It is not difficult to see that~\eqref{proof:empirical:Am:complement} is equivalent to that 
	there exists a $\beta \in E$ such that
	\begin{equation} \label{proof:empirical:Am:complement2}
		\big| \|\beta\|_{Nm}^2 - \|\beta\|_{X}^2 \big|	\ge  \epsilon .
	\end{equation} 
	We consider to apply Lemma~\ref{lemma:squaredconcen} to control the probability
	\begin{equation} \label{proof:empirical:interProb}
		\Pr\Big(\sup_{\beta\in E }  \Big|  \|\beta\|_{N m}^2 - \|\beta\|_{X}^2 \Big| \ge \epsilon \Big) .
	\end{equation}
	Set $f_{\beta}(x) =  \langle x, \beta \rangle$ as a class of functions indexed by $\beta\in E$. Correspondingly, denote
	$$
	Z(f_{\beta}) = 
	\Big|  \frac{1}{N} \sum_{n=1}^{N} \big( \langle x_{nm}, \beta\rangle^2-
	\Expect \langle x_{nm}, \beta\rangle^2 \big) \Big|.
	$$
	In order to apply Lemma~\ref{lemma:squaredconcen}, we also need to derive upper bounds for the constants $\sigma_1$  and $G_1$ under Condition~\ref{ass:subGaussian}.  
	For a fixed $\beta$, denote $Z_{nm} =\big| \langle x_{nm}, \beta\rangle^2 - \Expect \langle x_{nm}, \beta\rangle^2\big|$. By the basic properties of the Orlicz norm, we have for $\beta\in E$ that
	\begin{align*}
		\|Z_{nm}\|_{\psi_1} & \le  \|\langle x_{nm}, \beta \rangle^2 \|_{\psi_1} +\| \Expect \langle x_{nm}, \beta \rangle^2\|_{\psi_1}  \\
		& \le 2  \|\langle x_{nm}, \beta \rangle^2 \|_{\psi_1} = 2\|\langle x_{nm}, \beta \rangle \|^2_{\psi_2} \le 2 C_g \| \beta\|_{X}^2\le 2C_g.
	\end{align*}
	Besides, from the definition of the Orlicz norm, we can find that
	\begin{align*}
		2\ge \Expect \exp\Big( Z_{mn}/\|Z_{mn}\|_{\psi_1}\Big) \ge 1 + \frac{\Expect Z_{mn}^q}{q!\times \|Z_{mn}\|_{\psi_1}^q}, 
	\end{align*}
	which means that
	\begin{align*}
		\Expect Z_{mn}^q \le q!\times \|Z_{mn}\|_{\psi_1}^q \le \frac{q!}{2} (8C_g^2) \times (2C_g)^{q-2}.
	\end{align*}
	From the above, we can set $\sigma_1= 4C_g$ and $G_1  =4 C_g$ for  Lemma~\ref{lemma:squaredconcen}. In addition, 
	$$\mathrm{diam}_{\psi_2}(\mathcal{F})  = \sup_{\beta,\beta' \in E }
	\| \langle x_{nm}, \beta - \beta'\rangle \|_{\psi_2} \le 
	C_g \sup_{\beta,\beta' \in E }		\|  \beta - \beta' \|_{X} 
	\le 2 C_g.$$
	Also recall from~\eqref{eqn:complexity0:keybound}, we have
	for some constant $C_\gamma$ that
	\begin{equation} 	
		\gamma_{2}(E,d)  \le C_\gamma(K^{\frac{1}{2}}\wedge \eta_1^{-\frac{1}{4(q+d)}}).
	\end{equation}
	Now, applying Lemma~\ref{lemma:squaredconcen}, we have
	\begin{align*}
		&  \Pr\bigg[ \sup_{\beta \in E } Z(f_{\beta} )
		\ge C_1\Big\{ \frac{C_\gamma^2 ( K\wedge \eta_1^{-\frac{1}{2(q+d)}})}{N} 
		+ 	 \frac{2C_gC_\gamma ( K^{\frac{1}{2}}\wedge \eta_1^{-\frac{1}{4(q+d)}})}{\sqrt{N} } \Big\} + c_1\Big\{\sqrt{u}\frac{4C_g}{\sqrt{N}} + u\frac{4C_g}{N}\Big\}\bigg] \\
		\le & \exp\big(- u \big).
	\end{align*}
	Setting $u = K\wedge \eta^{-\frac{1}{2(q+d)}}_1$, it holds
	$$
	C_1\Big\{ \frac{C_\gamma^2 ( K\wedge \eta_1^{-\frac{1}{2(q+d)}})}{N} 
	+ 	 \frac{2C_gC_\gamma ( K^{\frac{1}{2}}\wedge \eta_1^{-\frac{1}{4(q+d)}})}{\sqrt{N} } \Big\} + c_1\Big\{\sqrt{u}\frac{4C_g}{\sqrt{N}} + u\frac{4C_g}{N}\Big\} \le C_e
	\frac{ K^{\frac{1}{2}}\wedge \eta_1^{-\frac{1}{4(q+d)}}}{\sqrt{N} } ,
	$$
	for $N$ large enough and for some constant $C_e$, since 
	$K\wedge \eta_1^{-\frac{1}{2(q+d)}}/N\to 0$. This implies a probability bound for~\eqref{proof:empirical:interProb} as
	\begin{align*}
		& \Pr\Big(	\sup_{\beta\in E }  \Big|  \|\beta\|_{N m}^2 - \|\beta\|_{X}^2 \Big| \ge \epsilon \Big) \\
		=&\Pr\bigg[	\sup_{\beta \in E }  \Big|  \frac{1}{N} \sum_{n=1}^{N} \big\{\langle x_{nm}, \beta \rangle^2 -
		\Expect \langle x_{nm}, \beta \rangle^2 \big\}\Big| \ge \epsilon \bigg] \\
		\le &  \Pr\bigg[	\sup_{\beta \in E } Z(f_{\beta} )
		\ge C_1\Big\{ \frac{C_\gamma^2 ( K\wedge \eta_1^{-\frac{1}{2(q+d)}})}{N} 
		+ 	 \frac{C_gC_\gamma ( K^{\frac{1}{2}}\wedge \eta_1^{-\frac{1}{4(q+d)}})}{\sqrt{N} } \Big\} + c_1\Big\{\sqrt{u}\frac{4C_g}{\sqrt{N}} + u\frac{4C_g}{N}\Big\}\bigg] \\
		\le & \exp\big\{- K\wedge \eta^{-\frac{1}{2(q+d)}}_1 \big\}.
	\end{align*}
	In the above, $\epsilon = C_e \big\{K\wedge \eta_1^{-1/ (2q+2d) } / N\big\}^{1/2}$.	$\hfill \blacksquare$

\section{Technical Proofs for Section~\ref{sec:manifold}}

\subsection{Proof of Lemma~\ref{lemma:metriccompatible}}

\textbf{Proof of~\eqref{lemma:metriccompatible1} and~\eqref{eqn:manifold:secondorderappr}.}
Consider a sufficiently small neighbor such that $\sQ_{\veta}(\mDelta) \le R_{M} \le 1/(2C_{\rm II})$. 	Let $\gamma$ be the geodesic over $\mathbb{M}$ starting at $\gamma(0)  = \bmB$ with the normalized  velocity $\dot{\gamma}(0) = \mDelta/\sQ_{\veta} \big( \mDelta\big)$. It satisfies that $\gamma(t)  = \exp_{\bmB}(\mDelta)$ at $t= \sQ_{\veta} \big( \mDelta\big)$. 	As $\mathbb{M}$ is a submanifold of $\bbR^{K\times M}$, we can view $\gamma(s)$ as a curve in $\bbR^{K\times M}$. Via Taylor expansion, it holds that
	\begin{align*}
		\exp_{\bmB}(\mDelta) - \bmB - \mDelta &=
		\gamma(t) - \gamma(0) - \dot\gamma(0) t \\
		& = \int_0^t (t-s) \ddot{\gamma}(s) \intd s \\
		& = \int_0^t (t-s) \mathrm{II}( \dot{\gamma}(s), \dot{\gamma}(s)) \intd s.
	\end{align*}
	The last inequality is due to Proposition~8.10 of~\cite{Lee2018}.
	Therefore, by Condition~\ref{ass:manifold},
	for $t$ small enough, it holds that
	\begin{align}
		\sQ_{\veta} (	\exp_{\bmB}(\mDelta) - \bmB - \mDelta) 
		& \le  \int_0^t (t-s) \sQ_{\veta} \big( \mathrm{II}( \dot{\gamma}(s), \dot{\gamma}(s))\big) \intd s \nonumber\\
		& \le  C_{\mathrm{II}} \int_0^t (t-s) \sQ_{\veta}^2 \big( \dot{\gamma}(s)\big) \intd s \label{proof:metric:velocityintegral}.
	\end{align}
	To bound~\eqref{proof:metric:velocityintegral}, we need to derive the upper bound for the size of $\dot{\gamma}(s)$ as measured by 
	$\sQ_{\veta} ( \cdot)$. This can be achieved by taking the first order derivative of the squared norm, i.e.,
	\begin{align*}
		\frac{\intd}{\intd s}	 \sQ^2_{\veta} \big( \dot{\gamma}(s)\big) &= 2 \tr\big\{\ddot{\gamma}(s)\trans\dot{\gamma}(s)\big\} +
		\sum_{j=1}^{P} \eta_j \big[\tr\big\{\ddot{\gamma}(s)\trans\mPi_{j1}\dot{\gamma}(s)\mPi_{j2}\big\}
		+ \tr\big\{\dot{\gamma}(s)\trans\mPi_{j1}\ddot{\gamma}(s)\mPi_{j2}\big\} \big] \\
		&\le 2 \| \ddot{\gamma}(s)\|_{F} \cdot\|\dot{\gamma}(s)\|_F +
		2\sum_{j=1}^{P} \eta_j
		\|\mPi_{j1}^{1/2}\ddot{\gamma}(s)\mPi_{j2}^{1/2} \|_F
		\cdot\|\mPi_{j1}^{1/2}\dot{\gamma}(s)\mPi_{j2}^{1/2}\|_F\\
		&\le 2 \Big\{\| \ddot{\gamma}(s)\|_{F}^2 +
		\sum_{j=1}^{P} \eta_j
		\|\mPi_{j1}^{1/2}\ddot{\gamma}(s)\mPi_{j2}^{1/2} \|_F^2\Big\}^{1/2}\\
		&\qquad \quad\times
		\Big\{\|\dot{\gamma}(s)\|_F^2 +
		\sum_{j=1}^{P} \eta_j
		\|\mPi_{j1}^{1/2}\dot{\gamma}(s)\mPi_{j2}^{1/2}\|_F^2\Big\}^{1/2} \\
		&\le 2\sQ_{\veta}(\ddot{\gamma}(s)) \sQ_{\veta}(\dot{\gamma}(s))  \\
		&=2\sQ_{\veta}(\mathrm{II}( \dot{\gamma}(s), \dot{\gamma}(s))) \sQ_{\veta}(\dot{\gamma}(s)) \\
		&\le 2C_{\mathrm{II}} \sQ_{\veta}^3( \dot{\gamma}(s)) .
	\end{align*}
The above inequality implies
$$
    \frac{\intd}{\intd s}	\big\{\sQ_{\veta}^2( \dot{\gamma}(s))\big\}^{-1/2}
    \ge -C_{\rm II},
$$
which further implies that
	$$
	\big\{\sQ_{\veta}^2( \dot{\gamma}(s))\big\}^{-1/2}- 	\big\{\sQ_{\veta}^2( \dot{\gamma}(0))\big\}^{-1/2}
\ge -C_{\rm II}s.
$$
Recall from the beginning of this proof, we have normalized the initial velocity such that $\dot{\gamma}(0) = \mDelta/\sQ_{\veta} \big( \mDelta\big)$. This means
$$
\sQ_{\veta}( \dot{\gamma}(s)) \le \big(1 - C_{\rm II}\cdot s \big)^{-1/2} \le 2,
$$
for $s \le t \le R_M \le 1/(2C_{\rm II})$.  Plug the above bound into~\eqref{proof:metric:velocityintegral}, it follows that
\begin{align} 
	\sQ_{\veta} (	\exp_{\bmB}(\mDelta) - \bmB - \mDelta) 
	\le 2
	\sQ_{\veta}^2 \big( \mDelta\big).
\end{align}
This exactly is the conclusion~\eqref{eqn:manifold:secondorderappr}.

Furthermore, by sub-additivity of $\sQ_{\veta} (\cdot)$,  it follows  that
\begin{equation*}
	\{1 - 2
	\sQ_{\veta} \big( \mDelta\big) \} Q_{\veta} \big( \mDelta\big)\le
	\sQ_{\veta} \big( \exp_{\bmB}(\mDelta) - \bmB  \big) \le \{1 + 2	\sQ_{\veta} \big( \mDelta\big) \} Q_{\veta} \big( \mDelta\big).
\end{equation*}
For $\sQ_{\veta}(\mDelta)\le R_M$, with  small enough $R_M$, we have
\begin{equation}
	(1/2) Q_{\veta} \big( \mDelta\big)\le 	\sQ_{\veta} \big(	\exp_{\bmB}(\mDelta) - \bmB  \big) \le 2 Q_{\veta} \big( \mDelta \big).
\end{equation}
This is the  conclusion~\eqref{lemma:metriccompatible1} of the lemma. $\hfill \blacksquare$ \\

\noindent \textbf{Proof of~\eqref{lemma:metriccompatible2}.} 	We can set a geodesic $\gamma$ to connect between
$\exp_{\bmB}(\mDelta_1)$ and $\exp_{\bmB}(\mDelta_2)$. Then, with a similar argument leading to the establishment of~\eqref{lemma:metriccompatible1}, we can derive a bound between their geodesic distance  $d_{\mathbb{M}}$ and Euclidean distance $\|\cdot\|_F$,
	\begin{align} 
		(1/2)d_{\mathbb{M}} \big(\exp_{\bmB}(\mDelta_1), \exp_{\bmB}(\mDelta_2)\big) & \le
		\big\|\exp_{\bmB}(\mDelta_1) - \exp_{\bmB}(\mDelta_2)\big\|_F \nonumber\\
		&\le 2  d_{\mathbb{M}} \big(\exp_{\bmB}(\mDelta_1), \exp_{\bmB}(\mDelta_2)\big).\label{proof:metric:deltaDiff1}
	\end{align}
The above result~\eqref{proof:metric:deltaDiff1} can also be verified by Proposition~6 of \cite{Smolyanov2007}.

	When $\mDelta_1$ and $\mDelta_2$ are linearly dependent, the conclusion~\eqref{lemma:metriccompatible2} directly follows from~\eqref{proof:metric:deltaDiff1}, because $d_{\mathbb{M}} \big(\exp_{\bmB}(\mDelta_1), \exp_{\bmB}(\mDelta_2)\big) = \| \mDelta_1 -\mDelta_2\|_F$. We now proceed to establish the second conclusion~\eqref{lemma:metriccompatible2} for linearly independent $\mDelta_1, \mDelta_2$ in a neighbor of $0$ in $\sD_{\bmB}(\subset T_{\bmB}\mathbb{M})$.  
	Notice, according to~(11) of~\cite{meyer1989toponogov}, we have
	\begin{align}
		&d_{\mathbb{M}} \big(\exp_{\bmB}(\mDelta_1), \exp_{\bmB}(\mDelta_2)\big) \nonumber\\
		& \qquad \qquad = \|\mDelta_1-\mDelta_2\|_F
		\, \big\{1-(1/12)K(\mDelta_1,\mDelta_2)
		(1+\langle\mDelta_1,\mDelta_2\rangle)\big\} 
		+o(\|\mDelta_1-\mDelta_2\|_F), \label{proof:metric:deltaDiff2}
	\end{align}
	where $K(\mDelta_1,\mDelta_2)$ is the sectional curvature of $\mathbb{M}$ for the plane spanned by $\mDelta_1$ and $\mDelta_2$.
	We can apply Gram-Schmidt procedure over 
	$\mDelta_1$ and $\mDelta_2$ to get two orthonormal tangent vectors
	$\widetilde{\mDelta}_1$ and $\widetilde{\mDelta}_2$.
	The sectional curvature can be upper bounded via
	\begin{align*}
		K(\mDelta_1,\mDelta_2)&= \|\mDelta_1\|_F	 \|\mDelta_2\|_F \, K\Big(\frac{\mDelta_1}{\|\mDelta_1\|_F},
		\frac{\mDelta_2}{\|\mDelta_2\|_F} \Big) \\
		&= \|\mDelta_1\|_F	 \|\mDelta_2\|_F \,
		K(\widetilde{\mDelta_1}, \widetilde{\mDelta_2} )\\
		&= \|\mDelta_1\|_F	 \|\mDelta_2\|_F \, \langle R(\widetilde{\mDelta}_1,\widetilde{\mDelta}_2)\widetilde{\mDelta}_2,\widetilde{\mDelta}_1\rangle \\
		&= \|\mDelta_1\|_F	 \|\mDelta_2\|_F \,
		\big\{\langle \mathrm{II}(\widetilde{\mDelta}_1,\widetilde{\mDelta}_1), \mathrm{II}(\widetilde{\mDelta}_2,\widetilde{\mDelta}_2)\rangle- \langle \mathrm{II}(\widetilde{\mDelta}_1,\widetilde{\mDelta}_2), \mathrm{II}(\widetilde{\mDelta}_2,\widetilde{\mDelta}_1)\rangle \big\}.
	\end{align*}
	In the above, the third equality uses Proposition 8.29 of \cite{Lee2018} and $R$ is the curvature tensor of $\mathbb{M}$. The last equality above uses the Gaussian Equation \citep[see Theorem~8.5 of][]{Lee2018}.
	Together with the identity (due to the bilinearity of the second fundamental form),
	\begin{align*}
		4\mathrm{II}(\widetilde{\mDelta}_1,\widetilde{\mDelta}_2) = \mathrm{II}(\widetilde{\mDelta}_1+\widetilde{\mDelta}_2,\widetilde{\mDelta}_1+\widetilde{\mDelta}_2)-
		\mathrm{II}(\widetilde{\mDelta}_1-\widetilde{\mDelta}_2,\widetilde{\mDelta}_1-\widetilde{\mDelta}_2),
	\end{align*}
	we can get an upper bound 
	\begin{align*}
		K(\mDelta_1,\mDelta_2) &= \|\mDelta_1\|_F	 \|\mDelta_2\|_F \cdot |\langle R(\widetilde{\mDelta}_1,\widetilde{\mDelta}_2)
		\widetilde{\mDelta}_2,\widetilde{\mDelta}_1\rangle| \\
		& \le  \|\mDelta_1\|_F \|\mDelta_2\|_F \Big\{\big|\langle \mathrm{II}(\widetilde{\mDelta}_1,\widetilde{\mDelta}_2), \mathrm{II}(\widetilde{\mDelta}_2,\widetilde{\mDelta}_1)\rangle\big|+ 
		\big|\langle \mathrm{II}(\widetilde{\mDelta}_1,\widetilde{\mDelta}_1), \mathrm{II}(\widetilde{\mDelta}_2,\widetilde{\mDelta}_2)\rangle\big|
		\Big\}\\
		& \le \|\mDelta_1\|_F	 \|\mDelta_2\|_F \Big\{ 
		\|\mathrm{II}(\widetilde{\mDelta}_1+\widetilde{\mDelta}_2,\widetilde{\mDelta}_1+\widetilde{\mDelta}_2)\|_F^2 
		+ \| \mathrm{II}(\widetilde{\mDelta}_1-\widetilde{\mDelta}_2,\widetilde{\mDelta}_1-\widetilde{\mDelta}_2)\|_F^2 \\
		& \qquad\qquad\qquad\qquad\qquad+
		\| \mathrm{II}(\widetilde{\mDelta}_1,\widetilde{\mDelta}_1)\|_F
		\| \mathrm{II}(\widetilde{\mDelta}_2,\widetilde{\mDelta}_2)\|_F \Big\} \\
		& \le C_{\rm II}^2 \|\mDelta_1\|_F	 \|\mDelta_2\|_F \Big\{ \sQ_{\veta}^4 \big(
		\widetilde{\mDelta}_1+\widetilde{\mDelta}_2\big) +
	    \sQ_{\veta}^4 \big(\widetilde{\mDelta}_1-\widetilde{\mDelta}_2\big) +
		\sQ_{\veta}^2 \big( \widetilde{\mDelta}_1\big)  \times\sQ_{\veta}^2 \big(\widetilde{\mDelta}_2\big)
		\Big\}.
	\end{align*}
	In the last inequality, we have applied Condition~\ref{ass:manifold}, i.e., 
	$\| \mathrm{II}(\widetilde{\mDelta}_1,\widetilde{\mDelta}_1)\|_F \le 
	\sQ_{\eta_1}\big(\mathrm{II}(\widetilde{\mDelta}_1,\widetilde{\mDelta}_1)\big)\le
	C_{\rm II}\sQ_{\veta}^2 \big( \widetilde{\mDelta}_1\big)$. The above upper bound with~\eqref{proof:metric:deltaDiff2} implies that, when $\sQ_{\veta}\big(\mDelta_1\big)$ and $\sQ_{\veta}\big(\mDelta_2\big)$ are small enough,  it holds that
	$$
	(1/2)\|\mDelta_1-\mDelta_2\|_F \le
	d_{\mathbb{M}} \big(\exp_{\bmB}(\mDelta_1), \exp_{\bmB}(\mDelta_2)\big) \nonumber\\
	\le 2\|\mDelta_1-\mDelta_2\|_F.
	$$
	Then, combined with~\eqref{proof:metric:deltaDiff1}, we arrive at the second conclusion~\eqref{lemma:metriccompatible2} of the lemma. $\hfill \blacksquare$

\subsection{Proof of Proposition~\ref{lemma:epupper}}

\begin{proof}
	Due to the basic properties of the Orlicz norm  $\|\cdot\|_{\psi_2}$ and Proposition~2.6.1 of~\cite{vershynin2018high}, there exists a constant $C_{\psi}$ $(>1)$ such that
	\begin{align}
		\| \sV(\mB) -	\sV(\mB')\|_{\psi_2}^2 
		&\le   \frac{C_{\psi}}{N^2} \sum_{m=1}^M\sum_{n=1}^N \| \ell_m(y_{nm}, \vx_{nm}\trans\vb_m) -
		\ell_m(y_{nm}, \vx_{nm}\trans\vb_m')\|_{\psi_2}^2 \nonumber\\
		&\qquad \qquad + \frac{C_{\psi}}{N^2} \sum_{m=1}^M\sum_{n=1}^N 
		\| \Expect\, \ell_m(y_{nm}, \vx_{nm}\trans\vb_m) - \Expect\, \ell_m(y_{nm}, \vx_{nm}\trans\vb_m')\|_{\psi_2}^2  \nonumber\\
		& \stackrel{(i)}{\le}  \frac{C_{\psi} C_L^2}{N^2} \sum_{m=1}^M\sum_{n=1}^N \Big\{ \|\langle x_{nm}, \beta_{m} - \beta_{m}'\rangle\|_{\psi_2}^2 + 4\Expect \langle x_{nm}, \beta_{m} - \beta_{m}'\rangle^2
		\Big\} \nonumber \\
		&\stackrel{(ii)}{\le} \frac{C_{\psi} C_L^2(C^2_g+4)}{N} \sum_{m=1}^M \|\beta_m -\beta_{m}'\|_{X}^2\nonumber\\
		& =  \frac{C_{\psi} C_L^2(C^2_g+4)}{N} \|\mB - \mB'\|_F^2.\nonumber
	\end{align}
	The  inequality~(i) uses the Lipschitz condition~\eqref{ass:lossfun:lip}, and the inequality~(ii) uses the sub-Gaussian Condition~\ref{ass:subGaussian}.
	The above bound means $\sV(\mB)$ is a sub-Gaussian process with the metric $d(\mB,\mB') = \|\mB - \mB'\|_F$.

For $\delta < R_M$, recall from~\eqref{eqn:gammaalpha:mani2tan} that
\begin{align*}
	\gamma_2\big( \exp_{\bmB}(\mathbb{N}(\bmB,\delta)),\; \|\cdot\|_F\big) 
	\le 8 \gamma_2( \mathbb{N}(\bmB,\delta),\; \|\cdot\|_F).
\end{align*}
It can also be easily checked that
\begin{align*}
    \mathrm{diam}\big\{ \exp_{\bmB}(\mathbb{N}(\bmB,\delta)) \big\} 
    &= \sup_{\mDelta,\mDelta'\in \mathbb{N}(\bmB,\delta)} 
    \big\| \exp_{\bmB}(\mDelta) -\exp_{\bmB}(\mDelta') \big\|_F \\
    &\le 4\sup_{\mDelta,\mDelta'\in \mathbb{N}(\bmB,\delta)}
    \big\| \mDelta -\mDelta' \big\|_F \\
    &\le 8\sup_{\mDelta\in \mathbb{N}(\bmB,\delta)} \big\| \mDelta  \big\|_F \\
    &\le 8\sup_{\mDelta\in \mathbb{N}(\bmB,\delta)} 
    \sQ_{\veta}( \mDelta )\\
    &\le 8\delta.
\end{align*}
According to the result of generic chaining \citep[see Theorem~8.5.5 of][]{vershynin2018high}, with probability at least $1 - 2 \exp(-u^2)$, it holds for $\mB = \exp_{\bmB}(\mDelta)$ with $\mDelta \in \mathbb{N}(\bmB, \delta)$ that
\begin{align}
	& \Big|\sV(\mB) - \sV(\bmB) \Big| \nonumber \\
	& \qquad \le   \frac{C_0 C_{\psi}^{1/2}C_L(C_g+2)}{N^{1/2}} \big[ \gamma_2\big( \exp_{\bmB}(\mathbb{N}(\bmB,\delta)),\; \|\cdot\|_F\big)   + u\cdot \mathrm{diam}\big\{ \exp_{\bmB}(\mathbb{N}(\bmB,\delta)) \big\} \big] \nonumber \\
	& \qquad \le \frac{C_{V}C_L(C_g+2)}{N^{1/2}} \big\{ \gamma_2( \mathbb{N}(\bmB,\delta),\; \|\cdot\|_F) + u \delta\big\},
\end{align}
where $C_0, C_V$ are positive constants. This completes the proof.
\end{proof}

\section{Technical Proofs for Section~\ref{sec:maintheory}}
\label{sec:appendix:mainproof}

\subsection{Check Condition~\ref{ass:lossfunConvex} for Quantile Regression}
\label{sec:appendix:mainproof:quantile}

Recall the setting of functional linear quantile regression from Section~\ref{sec:singleregression}.
Denote $P_{\epsilon |  X}(e)=P(Y_{m} - \langle X_m ,\beta_{0m} \rangle  \leq e|  X)$
as the cumulative distribution function for the true residual $Y_{m} - \langle X_m ,\beta_{0m} \rangle $, where $\beta_{0m}$ is the true slope function.  Similar to \cite{huang2021asymptotic},
we impose the following regularity condition.

\begin{condition} \label{ass:forquantile}
	There exist constants $L_2 \ge L_1 >0$ and $B$, such that for any interval $A=[a,b]\subseteq [-B,B]$, it holds that
	\begin{equation}
		L_1 (b-a) \le  P_{\epsilon |  X}(A).
	\end{equation}
Beside for any interval $A=[a,b]\subseteq \bbR$ on the real line, it holds that
	\begin{equation}
  P_{\epsilon |  X}(A) \le L_2 (b-a).
\end{equation}
\end{condition}

A sufficient condition for Condition~\ref{ass:forquantile} is: there exist constants $L_1,L_2$ and a density function $f(e)$ of  $P_{\epsilon |  X}(e)$ such that 
$$f(e) > L_1 >0  \text{ for }  e\in [-B,B], \quad  
\text{ and  }\quad f(e) < L_2 \text{ for all }  e\in\bbR.$$  
That is, the density is upper bounded on the whole real line $\bbR$ and bounded away from zero on the interval $[-B,B]$.

\vspace{10pt}

Recall from Section~\ref{sec:singleregression}, the loss function for $w$-quantile regression is $\ell_m(y, u) = (y-u)\times \{w - I(y<u) \}$, where $I(\cdot)$ is the indicator function. The Knight  identity  \citep{knight1998limiting} is
\begin{align}
\ell_m(y, u) - \ell_m(y, u')  = & 
\underbrace{(u-u')  \cdot \{I(y-u'\le 0) - w\} }_{D_1} \nonumber \\
&\qquad \qquad+
\underbrace{\int_0^{u-u'} \big\{ I(y-u'\le s) - I(y-u'\le 0)\big\} \intd s}_{D_2}. \label{eqn:Knight}
\end{align}
Setting $y = y_{1m}$, $u'=\langle x_{1m}, \beta'\rangle$  and $u=\langle x_{1m}, \beta\rangle$ in~\eqref{eqn:Knight} and taking expectation, we can see that
\begin{align} \label{eqn:quantile:lossD}
\bar{\sL}_{m}(\beta) - \bar{\sL}_{m}(\beta') = \Expect(D_1) + \Expect(D_2).
\end{align}
In the above, $\Expect(D_1)$ and  $\Expect(D_2)$ are the expected  value of the two terms on the right hand side of~\eqref{eqn:Knight}. Now, we have 
\begin{align*}
\Expect(D_2 ) &= \Expect \bigg[ \int_0^{ \langle x_{1m}, \beta -\beta' \rangle} 
	\Big\{ I\big(y_{1m}-\langle x_{1m},  \beta_{0m} \rangle\le 
	\langle x_{1m},  \beta'-\beta_{0m} \rangle + s \big)  \\
	&\qquad\qquad\qquad\qquad  -
	I\big(y_{1m}-\langle x_{1m},  \beta_{0m} \rangle\le 
	\langle x_{1m},  \beta'-\beta_{0m} \rangle \big)\Big\} \intd s \bigg] \\
	&= \Expect \bigg[\int_0^{ \langle x_{1m}, \beta -\beta' \rangle} 
	\Big\{ P_{\epsilon|X}\big(  \langle x_{1m},  \beta'-\beta_{0m} \rangle + s \big) -
	P_{\epsilon|X}\big(  \langle x_{1m},  \beta'-\beta_{0m} \rangle  \big) \Big\}  \intd s \bigg]
\end{align*}
The second equality takes the conditional expectation of $y_{1m}$ given $x_{1m}$.
By Condition~\ref{ass:forquantile}, it is easy to see that
\begin{align} \label{eqn:quantile:normUpper}
\Expect(D_2)  &\le L_2 \cdot \Expect \Big\{\int_0^{| \langle x_{1m}, \beta -\beta' \rangle |} 
	s  \intd s \Big\} =
	L_2 \cdot \Expect \Big\{
	\langle x_{1m}, \beta -\beta' \rangle^2/2\Big\} = \frac{L_2}{2} \| \beta - \beta'\|_X^2.
\end{align}

\vspace{5pt}
On the other hand, let $E_1$ be the event that
$$
\big| \langle x_{1m}, \beta -\beta' \rangle \big| < \frac{B}{2} \quad \text{ and }\quad 
\big| \langle x_{1m}, \beta_{0m} -\beta' \rangle \big| < \frac{B}{2} .
$$
We then have the lower bound
\begin{align} \label{eqn:quantile:norml}
    \Expect( D_2 )  &\ge \Expect \Big\{\int_0^{ \langle x_{1m}, \beta -\beta' \rangle} 
	L_1 s  \intd s \cdot I(E_1) \Big\} =
    L_1\cdot \Expect \Big\{	 \langle x_{1m}, \beta -\beta' \rangle^2/2 \cdot I(E_1) \Big\}.
\end{align}
Notice that
\begin{align} \label{eqn:quantile:xnormdecompose}
	\| \beta - \beta'\|_X^2 =
	\Expect \Big\{	\langle x_{1m}, \beta -\beta' \rangle^2 \cdot I(E_1) \Big\} +
	\Expect \Big\{	 \langle x_{1m}, \beta -\beta' \rangle^2 \cdot I(E_1^c) \Big\}
\end{align}
Recall Condition~\ref{ass:lossfunConvex} is stated for  $\beta,\beta'$ in a local  neighbor of the true $\beta_{0m}$, which means $\|\beta - \beta_{0m}\|_X$  and $\|\beta' - \beta_{0m}\|_X$ are small. When the neighbor is sufficiently small, the sub-Gaussian assumption in Condition~\ref{ass:subGaussian} means
\begin{align} \label{eqn:quantile:normtail}
	\Expect \Big\{ \langle x_{1m}, \beta -\beta' \rangle^2 \cdot I(E_1^c) \Big\} \le \| \beta - \beta'\|_X^2 /3.
\end{align}
Equations \eqref{eqn:quantile:xnormdecompose} and~\eqref{eqn:quantile:normtail} imply that
\begin{align}
	\Expect \Big\{\langle x_{1m}, \beta -\beta' \rangle^2 \cdot I(E_1) \Big\}
    \ge (2/3)	\| \beta - \beta'\|_X^2.
\end{align}
Together with~\eqref{eqn:quantile:norml}, we get
\begin{align}  \label{eqn:quantile:normLower}
	\Expect(D_2 ) & \ge  (2L_1/3) \| \beta - \beta'\|_X^2.
\end{align}


From the above, we can see $\Expect(D_2)$ in \eqref{eqn:quantile:lossD} is a second order term.
The Fr\'echet derivative  corresponds to the term $\Expect(D_1)$ in~\eqref{eqn:quantile:lossD}, which is 
$$
    D \bar{\sL}_{m}(\beta')[\beta - \beta'] = \Expect(D_1) =\Expect \Big[ 
    \langle x_{1m}, \beta-\beta'\rangle   \cdot
    \big\{\,I\big(y_{1m} \le  \langle x_{1m}, \beta'\rangle \big) - w\big\} \Big].
$$
Thus, collecting~\eqref{eqn:quantile:normUpper} and 
\eqref{eqn:quantile:normLower}, we have verified 
$$
(2L_1/3) 	\| \beta - \beta'\|_X^2\le 
E_m(\beta,\beta')=\bar{\sL}_{m}(\beta) - \bar{\sL}_{m}(\beta') - 	D \bar{\sL}_{m}(\beta')[\beta - \beta'] \le \frac{L_2}{2} 	\| \beta - \beta'\|_X^2.
$$ 
This means Condition~\ref{ass:lossfunConvex} indeed holds for the quantile regression loss under Condition~\ref{ass:forquantile}.

\subsection{Proof of Theorem~\ref{theorem:main}}

The proof of Theorem~\ref{theorem:main} is based on the following Lemmas~\ref{lemma:nablabound} and~\ref{lemma:main:estimationerror}, whose proofs are deferred to the next subsection. Because $\bmB$ is the solution to the manifold constraint optimization problem~\eqref{def:popuestimate:constraint}, the first-order optimality condition implies that
$\big\langle \mDelta, \nabla \bar{\sL} (\bmB)+\nabla\sP_{\veta} (\bmB) \big\rangle=\vzero$ for any $\mDelta\in T_{\bmB}\mathbb{M}$. The Euclidean gradient $\nabla \bar{\sL} (\bmB)+\nabla\sP_{\veta} (\bmB)$ at $\bmB$ is generally non-zero.
However, the magnitude of $\nabla \bar{\sL} (\bmB)+\nabla\sP_{\veta} (\bmB)$
provides an avenue  to measure the closeness between $\bmB$ and $\bmB_0$.
From the definition of the manifold constraint error~\eqref{eqn:perror:manifold}, we can see 
$\sE(\mathbb{M} ) = \sQ_{\veta}(\bmB - \bmB_0)$.
The following lemma bounds $\sQ_{\veta}(\bmB - \bmB_0)$ by the magnitude of $\nabla \bar{\sL} (\bmB)+\nabla\sP_{\veta} (\bmB)$. The latter is measured by the dual norm $\sQ_{\veta}^{*}(\cdot)$  of $\sQ_{\veta}(\cdot)$. The dual norm for a matrix $\mA$ is given by
$\sQ_{\veta}^{*}(\mA) = \sup_{\mK:\, \sQ_{\veta}(\mK)\le 1} \langle \mK, \mA\rangle$.

\begin{lemma} \label{lemma:nablabound}
	For the two optimal parameters $\bmB_0$ in~\eqref{def:popuestimate} and 
	$\bmB$ in~\eqref{def:popuestimate:constraint}, it holds that
	\begin{align} \label{eqn:nablabound2}
		\frac{1}{2\bar{C}_c}	\sQ_{\veta}^*\big(\nabla \bar{\sL} (\bmB)+ \nabla \sP_{\eta}(\bmB) \big) 
		\le  \sQ_{\veta}(\bmB - \bmB_0) 
		\le 		\frac{1}{2\ubar{c}_c}	\sQ_{\veta}^*\big(\nabla \bar{\sL} (\bmB)+ \nabla \sP_{\eta}(\bmB) \big) 
	\end{align}
	for $\ubar{c}_c = \min\{c_c,1\}$ and $\bar{C}_c = \max\{C_c, 1\}$.
\end{lemma}

\begin{lemma} \label{lemma:main:estimationerror}
	Set $\ubar{c}_c = \min\{c_c,1\}$ and 
	\begin{equation} \label{eqn:defineDeltaU}
			\delta_u:= (9/\ubar{c}_c)\big\{\widehat{\delta}_N +C_{V}C_L(C_g+2)u/N^{1/2} \big\}.
	\end{equation}	
	With probability at least  $1- 2\exp(-u^2)$, there exists a local optimal solution $\widehat{\mB}$ such that 
	$Q_{\veta} (\widehat{\mB} - \bmB)  \le \delta_u$.
\end{lemma}

\noindent \textbf{Proof of Theorem~\ref{theorem:main}}
	According to Lemma~\ref{lemma:main:estimationerror}, there exists a local estimate such that
	\begin{align*}
		&\Big\{ \sum_{m=1}^{M}
		\|\hat{\beta}_m- \bar{\beta}_m\|_{X}^2 + \sP_{\veta}(\hat{\vbeta} - \bar{\vbeta})  \Big\}^{1/2} 
		\le  (9/\ubar{c}_c) \big\{\widehat{\delta}_N +C_{V}C_L(C_g+2)u/N^{1/2} \big\}.
	\end{align*}
	Via the Cauchy-Schwartz inequality and the triangular inequality, it follows
	\begin{align*}
		& M^{-1/2} \sum_{m=1}^{M} \|\hat{\beta}_m-\beta_{0m}\|_{X} +  \sP_{\veta}^{1/2}(\hat{\vbeta}) \\
		\le  & \Big( \sum_{m=1}^{M} \|\hat{\beta}_m-\beta_{0m}\|_{X}^2\Big)^{1/2} +  \sP_{\veta}^{1/2}(\hat{\vbeta}) \\
		\le & \sqrt{2}\Big\{\sum_{m=1}^{M} \|\hat{\beta}_m-\beta_{0m}\|_{X}^2 +  \sP_{\veta}(\hat{\vbeta}) \Big\}^{1/2} \\
		\le & 
		\sqrt{2} \Big\{ \sum_{m=1}^{M}
		\|\hat{\beta}_m - \bar{\beta}_m\|_{X}^2 + \sP_{\veta}(\hat{\vbeta} - \bar{\vbeta})  \Big\}^{1/2}
		+ \sqrt{2}\Big\{ \sum_{m=1}^{M}
		\|\bar{\beta}_m- \bar{\beta}_{0m}\|_{X}^2 +  	\sP_{\veta}(\bar{\vbeta} - \bar{\vbeta}_0)   \Big\}^{1/2} 	\\
		&\qquad\qquad  +	\sqrt{2} \Big\{ \sum_{m=1}^{M}	\|\bar{\beta}_{0m}-\beta_{0m}\|_{X}^2   +
		\sP_{\veta}( \bar{\vbeta}_0) \Big\}^{1/2} \\
		\le & C_U/\ubar{c}_c \big[\ubar{c}_c \big\{\sE(\mathbb{M}) + \sE(\mathbb{S}_K)\big\}
		+  \widehat{\delta}_N + C_L(C_g+2)u/N^{1/2}\big],
	\end{align*}
where $C_u$ is an absolute constant.
$\hfill \blacksquare$

\subsection{Proof of Lemma~\ref{lemma:nablabound}}
\begin{proof} Note the penalty being quadratic and the convexity property of the loss function~\eqref{ass:lossfun:convex} means that
	\begin{align*}
		\sP_{\eta}(\mB) - 	\sP_{\eta}(\mB') &= \langle \mB - \mB',  \nabla \sP_{\eta}(\mB')  \rangle 	+ \sP_{\veta} (\mB - \mB'), \\
		\bar{\sL} (\mB)  - \bar{\sL} (\bmB) &\le \sum_{m=1}^{M}\Big\{
		\langle \vb_m - \vb'_m, \  \nabla \bar{\sL}_{m} (\vb'_m) \rangle
		+	C_c \| \vb_m - \vb'_m\|_{2}^2 \Big\},
	\end{align*}
	for all $\mB, \mB'\in\bbR^{K\times M}$.
	Combining the above together, we get
	\begin{align}
		&\bar{\sL} (\mB) + \sP_{\eta}(\mB) - \big\{ \bar{\sL} (\mB') +\sP_{\eta}(\mB') \big\} \nonumber \\
		& \qquad \le 
		\langle \mB - \mB', \  \nabla \bar{\sL} (\mB')+ \nabla \sP_{\eta}(\mB')  \rangle
		+	C_c \| \mB - \mB'\|_{F}^2 
		+ \sP_{\veta} (\mB - \bmB)  \nonumber  \\
		& \qquad \le 
		\langle \mB - \mB', \  \nabla \bar{\sL} (\mB')+ \nabla \sP_{\eta}(\mB')  \rangle
		+	\bar{C}_c \sQ_{\veta}^2( \mB - \mB'), \label{proof:nablabound:general}
	\end{align}
	for $\bar{C}_c = \max\{C_c, 1\}$.
	
	According to definition~\eqref{def:popuestimate}, $\bmB_0$ is the optimal solution among $\bbR^{K\times M}$. From~\eqref{proof:nablabound:general}, we know
	for any $\mB, \mB'\in\bbR^{K\times M}$ and for this specific $\bmB_0$ that
	\begin{align}
		\hspace{-10pt}\bar{\sL} (\bmB_0) + \sP_{\eta}(\bmB_0) 
		&\le		\bar{\sL} (\mB) + \sP_{\eta}(\mB)  \nonumber \\
		& \le \bar{\sL} (\mB') +\sP_{\eta}(\mB') +
		\langle \mB - \mB', \  \nabla \bar{\sL} (\mB')+ \nabla \sP_{\eta}(\mB')  \rangle
		+	\bar{C}_c \sQ_{\veta}^2( \mB - \mB').  \label{proof:nablabound:bmB0}
	\end{align}
	We specify the values of $\mB$ and $\mB'$ in~\eqref{proof:nablabound:bmB0} as follows:
	\begin{itemize} 
		\item Suppose $\bar{\mK}\in\bbR^{K\times M}$ is the matrix achieving the supremum in the definition of the dual norm, i.e.,
		\begin{equation*}
			\sQ_{\veta}^*\big(\nabla \bar{\sL} (\mB')+ \nabla \sP_{\eta}(\mB') \big) =
			\sup_{\sQ_{\veta}(\mK)\le 1}  	\langle \mK, \  \nabla \bar{\sL} (\mB')+ \nabla \sP_{\eta}(\mB')  \rangle =\langle \bar{\mK}, \  \nabla \bar{\sL} (\mB')+ \nabla \sP_{\eta}(\mB')  \rangle.
		\end{equation*}
		We can take $\mB = \bmB - 	\{1/(2\bar{C}_c)\} \sQ_{\veta}^*\big(\nabla \bar{\sL} (\mB')+ \nabla \sP_{\eta}(\mB') \big) \cdot \bar{\mK}$ in~\eqref{proof:nablabound:bmB0}. 
		
		\item  We can set $\mB' = \bmB$ as the optimal solution~\eqref{def:popuestimate:constraint} with the constraint $\mathbb{M}$.
	\end{itemize}
	Then, \eqref{def:popuestimate:constraint} leads to
	\begin{align*}
		\bar{\sL} (\bmB_0) + \sP_{\eta}(\bmB_0) 
		\le  \bar{\sL} (\bmB) +\sP_{\eta}(\bmB) - 1/(4\bar{C}_c)\big\{\sQ_{\veta}^*\big(\nabla \bar{\sL} (\bmB)+ \nabla \sP_{\eta}(\bmB) \big) \big\}^2,
	\end{align*}
	which implies that
	\begin{align*}
		1/(4\bar{C}_c) \big\{\sQ_{\veta}^*\big(\nabla \bar{\sL} (\bmB)+ \nabla \sP_{\eta}(\bmB) \big) \big\}^2
		&\le  \bar{\sL} (\bmB) +\sP_{\eta}(\bmB) - 	\bar{\sL} (\bmB_0) - \sP_{\eta}(\bmB_0) \\
		&\le \bar{C}_c\sQ_{\veta}^2(\bmB - \bmB_0).
	\end{align*}
	In the above, the last inequality follows from~\eqref{proof:nablabound:general} by setting  $\mB' = \bmB_0$, $\mB = \bmB$, and meanwhile noticing that
	$\nabla \bar{\sL} (\bmB_0)+ \nabla \sP_{\eta}(\bmB_0)=\vzero$ due to the optimality of $\bmB_0$. The above result completes the first inequality of~\eqref{eqn:nablabound2}, and the second inequality of~\eqref{eqn:nablabound2} can be proved similarly.
\end{proof}

\subsection{Proof of Lemma~\ref{lemma:main:estimationerror}}
\begin{proof}
	Our goal is to show that for all $\mDelta\in\sD_{\bmB}\subset T_{\bmB}\mathbb{M}$ with $\sQ_{\veta}(\mDelta) = \delta_u$, it holds that
	$$
	\sL(\exp_{\bmB}(\mDelta)) + \sP_{\veta}(\exp_{\bmB}(\mDelta))
	- \{\sL(\bmB) - \sP_{\veta}(\bmB) \} >0.
	$$
	with the required probability $1- 2\exp(-u^2)$.	Consider the difference of the objective function for $\mB = \exp_{\bmB}(\mDelta)$ with $\mDelta \in \mathbb{N}(\bmB, \delta_u)$,
\begin{align*}
	\sL( \mB) + \sP_{\veta} (\mB) - 
	\sL(\bmB) -\sP_{\veta} (\bmB) 
	=\bar{\sL} (\mB) - \bar{\sL} (\bmB) + \sP_{\eta}(\mB) 
	- \sP_{\eta}(\bmB) +  \sV(\mB) - \sV(\bmB) .
\end{align*}
Note the penalty being quadratic means that
\begin{align*}
	\sP_{\eta}(\mB) - 	\sP_{\eta}(\bmB) = \langle \mB - \bmB,  \nabla \sP_{\eta}(\bmB)  \rangle 
	+ \sP_{\veta} (\mB - \bmB),
\end{align*}
together with the lower bound~\eqref{ass:lossfun:convex} of $\bar{\sL}$, we have
\begin{align*}
	&\bar{\sL} (\mB) + \sP_{\eta}(\mB) - \bar{\sL} (\bmB) -\sP_{\eta}(\bmB) \\
	\ge & \sum_{m=1}^{M}\Big\{
	\langle \vb_m - \bvb_m, \  \nabla \bar{\sL}_{m} (\bvb_m) \rangle
	+	c_c \| \vb_m - \bvb_m\|_{2}^2 \Big\}	+ 	\langle \mB - \bmB,  \nabla \sP_{\eta}(\bmB)  \rangle 
	+ \sP_{\veta} (\mB - \bmB) \\
= &\;  \big\langle \mB - \bmB, \nabla \bar{\sL} (\bmB)+\nabla\sP_{\veta} (\bmB) \big\rangle
	+ \ubar{c}_c Q_{\veta}^2 (\mB - \bmB) ,
\end{align*}
where $\ubar{c}_c = \min\{c_c,1\}$.  Due to the optimality of $\bmB$ to the manifold constraint optimization problem~\eqref{def:popuestimate:constraint}, it holds that
$\big\langle \mDelta, \nabla \bar{\sL} (\bmB)+\nabla\sP_{\veta} (\bmB) \big\rangle=\vzero$ for $\mDelta\in T_{\bmB}\mathbb{M}$. It follows that
	\begin{align*}
			&\bar{\sL} (\mB) + \sP_{\eta}(\mB) - \bar{\sL} (\bmB) -\sP_{\eta}(\bmB) \\
	\ge &\;  \big\langle \mB - \bmB -\mDelta, \nabla \bar{\sL} (\bmB)+\nabla\sP_{\veta} (\bmB) \big\rangle
	+ \ubar{c}_c \sQ_{\veta}^2 (\mB - \bmB)\\
	\ge &  -
	\sQ_{\veta} (\mB - \bmB-\mDelta) \cdot
	\sQ_{\veta}^{*} \Big(\nabla \bar{\sL} (\bmB)+\nabla\sP_{\veta} (\bmB) \Big)	
	+ \ubar{c}_c \sQ_{\veta}^2 (\mB - \bmB)\\
		\stackrel{(i)}{\ge} &  -
	2\bar{C}_c	\sQ_{\veta} (\mB - \bmB-\mDelta) \cdot
    \sQ_{\veta}(\bmB - \bmB_0)
	+ \ubar{c}_c \sQ_{\veta}^2 (\mB - \bmB)\\
	\stackrel{(ii)}{\ge} &  -
	4\bar{C}_c   \sQ_{\veta}^2 \big( \mDelta\big) \cdot
	\sQ_{\veta}(\bmB - \bmB_0)
	+ (\ubar{c}_c/4) \sQ^2_{\veta}(\mDelta) .
\end{align*}
In (i) of the above, we have used Lemma~\ref{lemma:nablabound}. In (ii), we used the inequalities~\eqref{lemma:metriccompatible1} and~\eqref{eqn:manifold:secondorderappr}. 
In Theorem~\ref{theorem:main}, we have assumed that
the manifold approximation error $\sE(\mathbb{M})$ is sufficiently small with
$$
\sE(\mathbb{M})= \sQ_{\veta}(\bmB - \bmB_0)\le
\frac{\ubar{c}_c}{32\bar{C}_c } .
$$
Therefore, it follows
	\begin{align*}
	\bar{\sL} (\mB) + \sP_{\eta}(\mB) - \bar{\sL} (\bmB) -\sP_{\eta}(\bmB) 
	\ge   (\ubar{c}_c/8) \sQ^2_{\veta}(\mDelta).
\end{align*}
When $N$ is large enough, we also have $\delta_u < \min\{R_M, \mathrm{inj}(\bmB) \}$. Apply Proposition~\ref{lemma:epupper} to get
\begin{align*}
	& \sL( \mB) + \sP_{\veta} (\mB) - 
	\sL(\bmB) -\sP_{\veta} (\bmB) \\
	& \qquad = \bar{\sL} (\mB) - \bar{\sL} (\bmB) + \sP_{\eta}(\mB) 
	- \sP_{\eta}(\bmB) +  \sV(\mB) - \sV(\bmB) \\
	& \qquad \ge  (\ubar{c}_c/8) Q_{\veta}^2 (\mDelta) - \frac{C_{V}C_L(C_g+2)}{N^{1/2}} \big\{ \gamma_2( \mathbb{N}(\bmB,\delta_u),\; \|\cdot\|_F)   + u \delta_u\big\} \\	
    & \qquad = \delta_u \Big[ (\ubar{c}_c/8) \delta_u 
    - \frac{C_{V}C_L(C_g+2)}{N^{1/2}} \big\{ \gamma_2( \mathbb{N}(\bmB,\delta_u),\; \|\cdot\|_F)/\delta_u   + u  \big\}  \Big].
\end{align*}
Recalling from~\eqref{eqn:defineDeltaU} that $\delta_u \ge \widehat{\delta}_N$, and
$\delta_u \le \text{inj}(\bmB)$ for sufficiently large $N$. Therefore, we have
$$
    \frac{C_{V} C_L(C_g+2)}{N^{1/2}} 
    \frac{\gamma_2( \mathbb{N}(\bmB,\delta_u),\; \|\cdot\|_F)}{\delta_u  } 
    =
    \frac{C_{V}C_L(C_g+2)}{N^{1/2}} 
    \frac{\gamma_2( \mathbb{N}(\bmB, \widehat{\delta}_N),\; \|\cdot\|_F)}{ \widehat{\delta}_N } 
    \le \widehat{\delta}_N,
$$
as $\gamma_2( \mathbb{N}(\bmB,\delta),\; \|\cdot\|_F) / \delta $ is a constant in $\delta\in(0,\mathrm{inj}(\bmB))$. It follows that
\begin{align*}
	\sL( \mB) + \sP_{\veta} (\mB) - 
	\sL(\bmB) -\sP_{\veta} (\bmB) 
	& \ge  \delta_u \Big\{ (\ubar{c}_c/8) \delta_u 
	-  \widehat{\delta}_N - C_{V}C_L(C_g+2)u/N^{1/2} \Big\}\\
	&> 0.
\end{align*}
The conclusion of Lemma~\ref{lemma:main:estimationerror} follows from the above result.
\end{proof}

\section{Technical Proofs for Section~\ref{sec:reducedrate}}

\subsection{Proof of Lemma~\ref{lemma:reduced:approxerror}}

\begin{proof} In the proof, we only need to verify \eqref{lemma:reduced:approxerror1} and \eqref{lemma:reduced:approxerror2}.\\
	\textbf{Verification of~\eqref{lemma:reduced:approxerror1}}.
	Let $\tilde{\beta}_m$ be an optimal spline approximation to the true $\beta_{0m}$ such that $\tilde{\beta}_m$ satisfies Proposition~\ref{proposition:general:splineerror}, i.e.
	\begin{equation} \label{proof:reduced:approxerror:feasible0}
		\| \tilde{\beta}_m - \beta_{0m}\|_{X}  +\eta_1^{1/2} \| \tilde{\beta}_m\|_{\Gamma}
		\asymp K^{-\tau} + \eta^{1/2}_1 K^{(d-\nu)_{+}} .
	\end{equation}
	Define the vector $\tilde{\vbeta} = (\tilde{\beta}_1,\cdots, \tilde{\beta}_M)\trans$. Because $\bar{\vbeta}_0$ is the optimal solution to~\eqref{def:popuestimate} in the spline space $\mathbb{S}_{K}$, it holds that
	$$
	\bar{\sL} (\bar{\vbeta}_0)
	+ \sP_{\veta}(\bar{\vbeta}_0) \le
	\bar{\sL} (\tilde{\vbeta})
	+ \sP_{\veta}(\tilde{\vbeta}).
	$$
	We then have
	\begin{align}
		\ubar{c}_c\big\{\|\bar{\vbeta}_0 - \vbeta_0\|_X^2+\sP_{\veta}(\bar{\vbeta}_0)\big\}
		&\stackrel{(i)}{\le}  \bar{\sL} (\bar{\vbeta}_0) - \bar{\sL}(\vbeta_0)
		+ \sP_{\veta}(\bar{\vbeta}_0) \nonumber \\
		&\stackrel{(ii)}{\le}  \bar{\sL} (\tilde{\vbeta}) -  \bar{\sL}(\vbeta_0)
		+ \sP_{\veta}(\tilde{\vbeta}) \nonumber  \\
		&\stackrel{(iii)}{\le}	\bar{C}_c \big\{ \|\tilde{\vbeta}_0 - \vbeta_0\|_X^2+\sP_{\veta}(\tilde{\vbeta}_0)\big\}. \label{proof:reduced:approxerror:base0}
	\end{align}
	In the above, (i) and (iii) follow from the convexity assumption in Condition~\ref{ass:lossfunConvex}. As the true $\vbeta_0$ is the minimizer of $\bar{\sL}(\vbeta)$ in the  space $\mathbb{L}_2^{\nu}(\sT)$, we have $D \bar{\sL}_{m}(\beta_{0m})[\beta - \beta_{0m}]=0$
	for any $\beta\in \mathbb{L}_2(\sT)$. The  conclusion~\eqref{lemma:reduced:approxerror1} of Lemma~\ref{lemma:reduced:approxerror} is based on~\eqref{proof:reduced:approxerror:feasible0} and~\eqref{proof:reduced:approxerror:base0}.

	\vspace{10pt}
	
	\textbf{Verification of~\eqref{lemma:reduced:approxerror2}}. In the rest of this proof,	denote $[\mA]_{1:R}$ as the first $R$ columns of a matrix $\mA\in\bbR^{K\times M}$ and denote $[\mA]_{(R+1):M}$ as its remaining columns.
	Let $\widetilde{\mA} = [\bmU_0\bmD_0]_{1:R}$ and $\widetilde{\mV} = [\bmV_0]_{1:R}$ have the first $R$ singular vectors in their columns.  We can see 
	the rank-$R$ matrix	$\widetilde{\mB} = \bmA\bmV\trans$ is a feasible solution to the optimization problem in~\eqref{def:popuestimate:constraint}. The optimality of $\bmB$ to~\eqref{def:popuestimate:constraint} means 
	\begin{equation}
		\bar{\sL} (\bmB) + \sP_{\veta}(\bmB) \le  \bar{\sL} (\widetilde{\mB})
		+ \sP_{\veta}(\widetilde{\mB}). \label{proof:reduced:approxerror:baseinequal}
	\end{equation}
	Hence
	\begin{align}
		\ubar{c}_c \big\{\|\bmB  - \bmB_0\|_F^2 + \sP_{\eta_1}(\bmB - \bmB_0)\big\}
		&\stackrel{(i)}{\le} \bar{\sL} (\bmB) + \sP_{\veta}(\bmB) -\bar{\sL} (\bmB_0) - \sP_{\veta}(\bmB_0) \nonumber\\
		&\stackrel{(ii)}{\le} \bar{\sL} (\widetilde{\mB}) + \sP_{\veta}(\widetilde{\mB})-\bar{\sL} (\bmB_0) - \sP_{\veta}(\bmB_0)\nonumber\\
		&\stackrel{(iii)}{\le}	\bar{C}_c \big\{\|\widetilde{\mB}  - \bmB_0\|_F^2 + \sP_{\eta_1}(\widetilde{\mB} - \bmB_0)\big\} . \label{proof:reduced:approxerror:base2}
	\end{align}
	In the above (iii) is based the inequality~\eqref{proof:nablabound:general} by setting $\mB = \widetilde{\mB}$,  $\mB' = \bmB_0$ and noting $\nabla (\bar{\sL} (\bmB_0) + \sP_{\veta}(\bmB_0))  = \vzero$. The inequality~(i) can be easily checked by a similar argument. The inequality~(ii) follows from~\eqref{proof:reduced:approxerror:baseinequal}. 
	
	It immediately follows that $\|\widetilde{\mB}  - \bmB_0\|_F^2 =\sum_{r=R+1}^{M}\bar{\sigma}_{0r}^2$ from the construction of $\widetilde{\mB}$. Besides,
	\begin{align*}
		\sP_{\eta_1}(\widetilde{\mB} - \bmB_0) &= \eta_1\tr\big\{(\widetilde{\mB} - \bmB_0)\trans \mGamma (\widetilde{\mB} - \bmB_0) \big\}\\
		& = \eta_1\tr\big\{ [\bmU_0]_{(R+1):M}\trans \mGamma [\bmU_0]_{(R+1):M}[\bmD_0]_{(R+1):M}^2\big\} \\
		&\le \eta_1\sum_{r=R+1}^{M} \bar{\sigma}_{0r}^2 \bar{\vu}_{0r}\trans\mGamma \bar{\vu}_{0r} \\
		&\le \eta_1\sum_{r=1}^{M} \bar{\sigma}_{0r}^2 \bar{\vu}_{0r}\trans\mGamma \bar{\vu}_{0r}
		=\sP_{\eta_1}(\bmB_0).
	\end{align*}
	Combining the above discussion with~\eqref{proof:reduced:approxerror:base2} leads to the conclusion~\eqref{lemma:reduced:approxerror2}.		
\end{proof}

\subsection{Proof of Lemma~\ref{lemma:reduced:gamma}}

\begin{proof}
	Set $\bmA = \bmU\bmD$,  $\mDelta_{a1} = \bmU\mM + \mU_p$ and 
	$\mDelta_{v1} = \mV_p$. Then $\bmB = \bmA\bmV\trans$ and the tangent  vector in~\eqref{eqn:reduced:tangentspace} is re-parameterized as $\mDelta= \mDelta_{a1}\bmV\trans + \bmU\mDelta_{v1}\trans $. Because $\sD_{\bmB}$ is a subset of $ T_{\bmB}\mathbb{M}$, the neighborhood $\mathbb{N}(\bmB,\delta)$ in~\eqref{eqn:reduced:neighbor} is contained inside a larger set $\sN_1(\bmB,\delta)$ where
	\begin{align}
		\sN_1(\bmB,\delta) = \big\{ \mDelta\in T_{\bmB}\mathbb{M}:\ \| (\mI + \eta_1 \mGamma)^{1/2}\mDelta\|_F\le \delta \big\}.
	\end{align}
	It follows that $\gamma_{2}(\mathbb{N}(\bmB,\delta), \|\cdot\|_F)\le 	\gamma_{2}(\sN_1(\bmB,\delta), \|\cdot\|_F)$ and we proceed to derive bound for the latter.
	
	For any two matrices $\mDelta,\mDelta'\in \sN_1(\bmB,\delta)$ in the local neighborhood, the Frobenius norm of their difference has an upper bound
	\begin{align}
		d(\mDelta,\mDelta')=&	 \|\mDelta - \mDelta'\|_F  
		\le   \|\mDelta_{a1} - \mDelta_{a1}'\|_F + \|\mDelta_{v1} - \mDelta_{v1}'\|_F, \label{eqn:complexity1:dbbupper}
	\end{align}
	with the structure $\mDelta  = \mDelta_{a1}\bmV\trans + \bmU\mDelta_{v1}\trans$ and $\mDelta'
	= \mDelta_{a1}'\bmV\trans + \bmU(\mDelta_{v1}')\trans$. 
    For $\mDelta\in\sN_1(\bmB,\delta)$, the inequality $\|(\mI+\eta_1\mGamma)^{1/2} \mDelta\|_F  \le \delta$ implies
	$$
	\|(\mI+\eta_1\mGamma)^{1/2} \mDelta\|_F^2=
	\|(\mI+\eta_1\mGamma)^{1/2} \mDelta_{a1} \bmV\trans\|_F^2 + 
	\|(\mI+\eta_1\mGamma)^{1/2} \bmU\mDelta_{v1}\trans\|_F^2 \le \delta^2.
	$$
	The above equality holds because $\bmV\trans \mDelta_{v1} = \vzero$.
	This means, when $\mDelta\in\sN_1(\bmB,\delta)$, the parameters ($\mDelta_{a1},\mDelta_{v1}$) of the tangent vector  $\mDelta  = \mDelta_{a1}\bmV\trans + \bmU\mDelta_{v1}\trans$  satisfy that
	\begin{equation} \label{eqn:complexity1:tangentUpper}
		\|(\mI+\eta_1\mGamma)^{1/2} \mDelta_{a1} \|_F \le \delta \text{ and }
		\|\mDelta_{v1}\|_F \le \delta.
	\end{equation}
	Let us denote
	$$
	T_a = \{\mDelta_{a1}:\ \|(\mI+\eta_1\mGamma)^{1/2} \mDelta_{a1} \|_F \le \delta \},\quad\text{and}\quad T_v = \{\mDelta_{v1}:\ \|\mDelta_{v1}\|_F \le \delta, \mDelta_{v1}^T\bmV = \vzero\}.
	$$
	From the analysis related to~\eqref{eqn:complexity1:tangentUpper}, we find a super-set $T$ for $	\mathbb{N}(\bmB,\delta) $ as
	\begin{equation} \label{eqn:complexity1:seqsubsets}
		\mathbb{N}(\bmB,\delta)  \subset		\sN_1(\bmB,\delta)  \subset T:= \{\mDelta:\, \mDelta  = \mDelta_{a1}\bmV\trans + \bmU\mDelta_{v1}\trans,\ \mDelta_{a1}\in T_a,\ \mDelta_{v1} \in T_b\}.
	\end{equation}

	Our goal turns to find a complexity bound for the set $T$. 	
Suppose	the metric $d(\mA, \mB)  = \|\mA - \mB\|_F$ is the Frobenius norm.
	Let $\{T_{a,n}\}_{n=1}^\infty$ be an admissible sequence of subsets for $T_{a}$, such that 
	$$ \sup_{\mDelta_{a1} \in T_a} \sum_{n=1}^\infty
	2^{n/2} d(\mDelta_{a1}, T_{a,n-1})\le  4 \gamma_{2}(T_a, d).$$	
	Similarly, let  $\{T_{v,n}\}_{n=1}^\infty$ be an admissible sequence of subsets for $T_{v}$, such that $$ \sup_{\mDelta_{v1} \in T_v} \sum_{n=1}^\infty
	 2^{n/2} d(\mDelta_{v1}, T_{v,n-1})\le  4 \gamma_{2}(T_v, d).$$ 
	 Then,  for $n=1,2\cdots$, the sets
	$$
	T_{n} = \{\mDelta:\, \mDelta  = \mDelta_{a1}\bmV\trans + \bmU\mDelta_{v1}\trans,\ \mDelta_{a1}\in T_{a,n-1},\ \mDelta_{v1} \in T_{b,n-1}\}
	$$
	constitute an admissible sequence of subsets for $T$, because $|T_n| \le |T_{a,n-1}| \cdot |T_{b,n-1}|\le 2^{2^n}$. Therefore, we have
	\begin{align}
		\gamma_{2}(\mathbb{N}(\bmB,\delta), \|\cdot\|_F) b 
		&	\stackrel{(i)}{\lesssim} 	\gamma_{2}(T, \|\cdot\|_F) \nonumber  \\
		& \le \sup_{\mDelta \in T} \sum_{n=1}^\infty 2^{n/2}d(\mDelta, T_n) + d(\mDelta,T_0)  \nonumber  \\
		& \stackrel{(ii)}{\le}  \sup_{\mDelta_{a1} \in T_a} \sum_{n=1}^\infty
		2^{n/2} d(\mDelta_{a1}, T_{a,n-1}) \nonumber  \\
		& \qquad \qquad +  \sup_{\mDelta_{v1} \in T_v}  \sum_{n=1}^\infty 2^{n/2} d(\mDelta_{v1}, T_{v,n-1}) + \text{diam}(T)  \nonumber \\
		& \le 4\big\{ \gamma_{2}(T_a, d) + \gamma_{2}(T_v, d)+\delta\big\}, \label{eqn:complexity1:gammaT}
	\end{align}
	where $\text{diam}(T)$ is the diameter of the set $T$ as measured by the metric $d(\mDelta, \mDelta')  = \|\mDelta - \mDelta'\|_F$.
	In the above, the inequality~(i) uses the relation~\eqref{eqn:complexity1:seqsubsets}, and the inequality~(ii) is based on~\eqref{eqn:complexity1:dbbupper}.

	Note each column of $\mDelta_{a1}\in T_a$ belongs to an ellipsoid of the type~\eqref{eqn:complexity0:ellipsoid} but with different radius.  As a result, similar to~\eqref{eqn:complexity0:keybound}, we can show that
	\begin{equation} \label{eqn:complexity1:gammaTa}
		\gamma_{2}(T_a, d) \lesssim R^{1/2} 
		\big\{ K^{\frac{1}{2}}\wedge \eta_1^{-\frac{1}{4(q+d)}}\big\} \delta.
	\end{equation}
	In addition, for $\mDelta_{v1}\in\bbR^{M\times R}$ in the set $T_v$, each column of $\mDelta_{v1}$ is orthonormal to $\bmV$. As $\dim(\text{span}(\bmV)) = R$, each column of $\mDelta_{v1}$  belongs to a $M-R$ subspace. Therefore, $\mDelta_{v1}$  belongs to a Euclidean space of dimension $R(M-R)$. For $\epsilon\le \delta$ and some constant $C_v$, the entropy number  of $T_v$ is known as
	$\log N(\epsilon, T_v, d)\leq \{ R(M-R) \}\log(C_v /\epsilon)$. We find the upper bound
	\begin{equation}
		\gamma_{2}(T_v, d) \lesssim \int_0^{\delta} \sqrt{\log N(\epsilon, T_v, d)} \intd \epsilon \lesssim
		R^{1/2}\big(M-R\big)^{1/2}\delta. \label{eqn:complexity1:gammaTv}
	\end{equation}
	The conclusion follows by combining~\eqref{eqn:complexity1:gammaT}, \eqref{eqn:complexity1:gammaTa}, and \eqref{eqn:complexity1:gammaTv}.
\end{proof}

\subsection{Proof of Lemma~\ref{lemma:reduced:second}}
\begin{proof} The second fundamental form and the Weingarten map have the adjoint relation \citep[see (8.4) of][]{Lee2018}. Suppose $\mN\in N_{\bmB}\mathbb{M}$ is a normal vector and $\mDelta_1,\mDelta_2\in T_{\bmB}\mathbb{M}$ are two tangent vectors, then it holds that 
	\begin{equation} \label{proof:second:adjoint}
		\langle W_{\mN}(\mDelta_1), \mDelta_2\rangle = \langle \mN, \rm II(\mDelta_1, \mDelta_2) \rangle,
	\end{equation}
	where $W_{\mN}$ is the Weingarten map in the direction of $\mN$. From Section~4.5 of~\cite{absil2013extrinsic}, we can find
	the Weingarten map of the fixed rank manifold is
	\begin{equation} \label{proof:second:weingarten}
		W_{\mN}(\mDelta_1) = \mN\mDelta_1\trans(\bmB^{+})\trans
		+(\bmB^{+})\trans\mDelta_1\trans\mN.
	\end{equation} 
	As a result, \eqref{eqn:reduced:secondFund} can be concluded by combining~\eqref{proof:second:adjoint} and~\eqref{proof:second:weingarten}.
	
	Next, we check that Condition~\ref{ass:manifold} is satisfied if the $R$-th singular value of $\bmB$ is bounded away from zero. Because the second fundamental form is bi-linear with respect to its two arguments, to verify Condition~\ref{ass:manifold}, it suffices to show
	$$
	\sQ_{\veta} ( \mathrm{II}(\mDelta,\mDelta) ) \le C_{\mathrm{II}},
	$$
	for some constant $C_{\mathrm{II}}$ and all $\mDelta$ satisfying  $\sQ_{\eta_1} (\mDelta)\le 1$.
	
	Consider a tangent vector $\mDelta$ specified by~\eqref{eqn:reduced:tangentspace} and satisfying  $\sQ_{\eta_1} (\mDelta)\le 1$, we can check that
	\begin{align*}
		1\ge \sQ_{\eta_1}^2 (\mDelta)&=  \| (\mI + \eta_1 \mGamma)^{1/2}\mDelta\|_F^2 \\
    & \ge \| \mDelta\|_F^2 \\
		&= \| \mM\|_F^2 + \|\mU_p\|_F^2 + \|\mV_p\|_F^2.
	\end{align*}
	This means the Frobenius norm of the matrices $\mM$, $\mU_p$, and $\mV_p$ in~\eqref{eqn:reduced:tangentspace} are all bounded. In particular, we have
	\begin{equation} \label{proof:second:M}
		\|\mM\|_F\le 1.
	\end{equation}
	On the other hand, we have
	\begin{align*}
		1\ge \sQ_{\eta_1}^2 (\mDelta) &=  \| (\mI + \eta_1 \mGamma)^{1/2}\mDelta\|_F^2 \\
		&= \|(\mI+\eta_1\mGamma)^{1/2}(\bmU\mM\bmV\trans + \mU_p\bmV\trans)\|_F^2 +
		\|(\mI+\eta_1\mGamma)^{1/2}(\bmU\mV_p\trans)\|_F^2,
	\end{align*}
	which implies $\|(\mI+\eta_1\mGamma)^{1/2}(\bmU\mM\bmV\trans + \mU_p\bmV\trans)\|_F$ is also bounded, i.e.,
	\begin{equation} \label{proof:second:UMVUV}
		\|(\mI+\eta_1\mGamma)^{1/2}(\bmU\mM\bmV\trans + \mU_p\bmV\trans)\|_F \le 1.
	\end{equation}
	
	For an arbitrary $\mDelta$ in~\eqref{eqn:reduced:tangentspace}, the second fundamental form can be simplified as
	\begin{align*}
		\rm II(\mDelta,\mDelta) =  P_{\bmB}^{\perp}
		(\mDelta \bmB^+\mDelta+ \mDelta \bmB^+\mDelta) = 2\mU_p \bmD^{-1}\bmV_p\trans.
	\end{align*}
	Then, 
	\begin{align*}
		\sQ_{\veta} ( \mathrm{II}(\mDelta,\mDelta) ) &\le
		2\|(\mI+\eta_1\mGamma)^{1/2}\mU_p \bmD^{-1}\bmV_p\trans\|_F\\
		&\le 2 \|(\mI+\eta_1\mGamma)^{1/2}\mU_p\|_F \times\| \bmD^{-1}\|\times \|\bmV_p\trans\| \\
		&\le 2 \bar{\sigma}_R^{-1}\|(\mI+\eta_1\mGamma)^{1/2}\mU_p\|_F,
	\end{align*}
	where $\bar{\sigma}_R$ is the $R$-th singular value of $\bmB$.
	To bound the right hand side of the above, 
	we have by~\eqref{proof:second:M}, \eqref{proof:second:UMVUV}, and the triangular inequality that
	\begin{align*}
		\|(\mI+\eta_1\mGamma)^{1/2}\mU_p\|_F
		&\le \|(\mI+\eta_1\mGamma)^{1/2}(\bmU\mM\bmV\trans + \mU_p\bmV\trans)\|_F + \|(\mI+\eta_1\mGamma)^{1/2}\bmU\mM\bmV\trans\|_F \\
		&\le 1 + \|(\mI+\eta_1\mGamma)^{1/2}\bmU\|_F \|\mM\|_F \|\bmV\| \\
		&\le 1 + \|(\mI+\eta_1\mGamma)^{1/2}\bmU\|_F .
	\end{align*}
	The proof is thus completed if we can show the right hand side of the last line in the above is bounded. 
    This can be verified by observing $\|(\mI+\eta_1\mGamma)^{1/2}\bmB_0\|$ is bounded due to Lemma~\ref{lemma:reduced:approxerror} and under the assumption that $\eta_1K^{2(d-\nu)_+}$ is bounded.
Note that $\bmU$ contains the left singular vectors of $\bmB$, and we have
	\begin{align*}
		\bar{\sigma}_R \|(\mI+\eta_1\mGamma)^{1/2}\bmU\|_F  
		&\le \|(\mI+\eta_1\mGamma)^{1/2}\bmU\bmD\bmV\trans\|_F   \\
		&= \sQ_{\eta_1} (\bmB) \le \sQ_{\eta_1} (\bmB - \bmB_0) + \sQ_{\eta_1} (\bmB_0) < \infty,
	\end{align*}
    which finishes the proof.
\end{proof}

\subsection{Proof of Corollary~\ref{corollary:reduced:rate1}}	

\textbf{Subcase~(i).}  When $K \lesssim \eta_1^{-\frac{1}{2(\iota\vee\tau)}} \le \eta_1^{-\frac{1}{2(\iota\wedge\tau)}}$, the upper bound of the convergence rate in this case can be written as
$$
\frac{1}{M}\sum_{m=1}^{M} \Big\{\|\hat{\beta}_m-\beta_0\|_{X} + \eta^{1/2}_1\|\hat{\beta}_m\|_{\Gamma}\Big\} =O_p\Big( \Big(\frac{R}{MN}\Big)^{1/2}K^{1/2} + K^{-\tau}
\Big).
$$ 
The right hand side is optimized with $K \asymp (MN/R)^{1/(2\tau + 1)}$, which leads to
$$
\frac{1}{M}\sum_{m=1}^{M} \Big\{\|\hat{\beta}_m-\beta_0\|_{X} + \eta^{1/2}_1\|\hat{\beta}_m\|_{\Gamma}\Big\} =O_p\Big( (MN/R)^{-\tau/(2\tau+1)}\Big),$$
when $\eta_1 \lesssim K^{-2(\iota\vee \tau)}\asymp (MN/R)^{-2(\iota\vee \tau)/(2\tau + 1)} $.\\

\noindent \textbf{Subcase~(ii).} 	When $ \eta_1^{-\frac{1}{2(\iota\vee\tau)}} \le \eta_1^{-\frac{1}{2(\iota\wedge\tau)}} \lesssim K$, the upper bound of the convergence rate in this case can be written as 
$$
\frac{1}{M}\sum_{m=1}^{M} \Big\{\|\hat{\beta}_m-\beta_0\|_{X} + \eta^{1/2}_1\|\hat{\beta}_m\|_{\Gamma}\Big\} =O_p\Big( \Big(\frac{R}{MN}\Big)^{1/2}\eta_1^{-{1}/{(4\iota)}} + \eta_1^{1/2}\Big).
$$ 
The right hand side is optimized with $\eta_1 \asymp (MN/R)^{-2\iota/(2\iota + 1)}$, which leads to
$$
\frac{1}{M}\sum_{m=1}^{M} \Big\{\|\hat{\beta}_m-\beta_0\|_{X} + \eta^{1/2}_1\|\hat{\beta}_m\|_{\Gamma}\Big\} =O_p\Big( (MN/R)^{-\iota/(2\iota+1)}\Big),$$
provided $K \gtrsim \eta_1^{-\frac{1}{2(\iota\wedge\tau)}} \asymp (MN/R)^{\frac{\iota}{(2\iota + 1)(\iota\wedge\tau)}} $.\\

\noindent \textbf{Subcase~(iii).} When $ \eta_1^{-\frac{1}{2(\iota\vee\tau)}} \lesssim K \lesssim \eta_1^{-\frac{1}{2(\iota\wedge\tau)}}$. The discussion below is conducted separately for $\tau\ge\iota $ and $ \tau < \iota $. We find they correspond to Subcases~(ii) and~(i), respectively.

\begin{enumerate}
	\item  When $\tau \le \iota $, the upper bound of the convergence rate is then of order 
	$$
	(1/M) \Big\{\sum_{m=1}^{M} \|\hat{\beta}_m-\beta_{0m}\|_{X} +  \sP_{\veta}^{1/2}(\hat{\vbeta})\Big\} =O_p\Big( \Big(\frac{R}{MN}\Big)^{1/2}\eta_1^{-{1}/{(4\iota)}} + K^{-\tau}\Big).
	$$
	After optimizing the right hand side with respect to $K$ in the range $\eta_1^{-\frac{1}{2(\iota \vee \tau)}} \lesssim K \lesssim  \eta_1^{-\frac{1}{2(\iota \wedge \tau)}}$, we find the optimal $K$ is its upper bound in this range, i.e.,
	$K\asymp  \eta_1^{-\frac{1}{2(\iota \wedge \tau)}}\asymp  \eta_1^{-\frac{1}{2\tau}}$. Plugging in this value to the upper bound of the convergence rate, we obtain
	$$
	\frac{1}{M}\sum_{m=1}^{M} \Big\{\|\hat{\beta}_m-\beta_0\|_{X} + \eta^{1/2}_1\|\hat{\beta}_m\|_{\Gamma}\Big\} =O_p\Big( \Big(\frac{R}{MN}\Big)^{1/2}\eta_1^{-{1}/{(4\iota)}} + \eta_1^{1/2}\Big).
	$$ 
    It is seen that the above has the same form as in Subcase~(ii), and we accordingly get the rate of convergence as
	$$
	\frac{1}{M}\sum_{m=1}^{M} \Big\{\|\hat{\beta}_m-\beta_0\|_{X} + \eta^{1/2}_1\|\hat{\beta}_m\|_{\Gamma}\Big\} = O_p\Big( (MN/R)^{-\iota/(2\iota+1)}\Big),$$
	when $\eta_1 \asymp (MN/R)^{-2\iota/(2\iota + 1)}$ and
	$K \asymp (MN/R)^{\frac{\iota}{\tau(2\iota + 1)}} $.

	\item When $\tau > \iota $, the upper bound of the convergence rate is then of order 
	$$
	(1/M) \Big\{\sum_{m=1}^{M} \|\hat{\beta}_m-\beta_{0m}\|_{X} +  \sP_{\veta}^{1/2}(\hat{\vbeta})\Big\} =O_p\Big( \Big(\frac{R}{MN}\Big)^{1/2}K^{1/2} + \eta_1^{1/2}\Big).
	$$
	After optimizing the right hand side with respect to $\eta_1$ in the range $\eta_1^{-\frac{1}{2(\iota \vee \tau)}} \lesssim K \lesssim  \eta_1^{-\frac{1}{2(\iota \wedge \tau)}}$, we find that the optimal $\eta_1$ satisfies
	$\eta_1 \asymp  K^{-2\tau}$. Plugging in this value to the upper bound of the convergence rate, we obtain
	$$
	\frac{1}{M}\sum_{m=1}^{M} \Big\{\|\hat{\beta}_m-\beta_0\|_{X} + \eta^{1/2}_1\|\hat{\beta}_m\|_{\Gamma}\Big\}	 = O_p\Big( \Big(\frac{R}{MN}\Big)^{1/2}K^{1/2} + K^{-\tau}\Big).
	$$ 
    It is seen that the above has the same form as in Subcase~(i), and we accordingly get the rate of convergence as
	$$
	\frac{1}{M}\sum_{m=1}^{M} \Big[\|\hat{\beta}_m-\beta_0\|_{X} + \eta^{1/2}_1\|\hat{\beta}_m\|_{\Gamma}\Big]	 =O_p\Big((MN/R)^{-\tau/(2\tau+1)}\Big),$$
	when $K \asymp (MN/R)^{1/(2\tau + 1)}$ and $\eta_1 \asymp (MN/R)^{-2\tau /(2\tau + 1)} $.
\end{enumerate}

\subsection{Proof of Corollary~\ref{corollary:reduced:rate2}}

\begin{proof}
    When $\mathfrak{o}+1\ge d > \nu$, notice  $\tau = \nu + \{q \wedge (\mathfrak{o}+1)\}/2$. The rate of convergence is upper bounded by
    $$
        O_p\Big(\frac{R^{1/2}\{K^{1/2}\wedge \eta_1^{-1/4\iota} + (M-R)^{1/2}\} }{M^{1/2} N^{1/2}} + K^{-\tau} + \eta_1^{1/2} K^{(d-\nu)}\Big).
    $$
\noindent \textbf{Subcase (i).} When $K \lesssim \eta_1^{-1/(2\iota)} < \eta_1^{-1/(2\tau)}$, it holds 
	that $\eta_{1}^{1/2}\lesssim K^{-\iota}$. This means
	$$
	K^{(d-\nu)}\eta_{1}^{1/2}\lesssim K^{-\iota+d-\nu} = K^{-q-\nu} \lesssim K^{-\tau}.
	$$
	The upper bound of the convergence rate in this case is then of the order as  
	$$
	\frac{1}{M}\sum_{m=1}^{M} \Big\{\|\hat{\beta}_m-\beta_0\|_{X} + \eta^{1/2}_1\|\hat{\beta}_m\|_{\Gamma}\Big\}	=O_p\Big(		\frac{R^{1/2} K^{1/2}}{M^{1/2} N^{1/2}} + K^{-\tau}\Big).
	$$
	Therefore, by setting $\eta_1 \lesssim (MN/R)^{-2\iota/(2\tau + 1)}$ and $K \asymp (MN/R)^{1/(2\tau + 1)}$, we obtain the rate of convergence $(MN/R)^{-\tau/(2\tau + 1)}$. \\

\noindent 	\textbf{Subcase (ii).} When $\eta_1^{-1/(2\iota)} \lesssim K$, the upper bound of the convergence rate becomes 
	$$
	\frac{1}{M}\sum_{m=1}^{M} \Big\{\|\hat{\beta}_m-\beta_0\|_{X} + \eta^{1/2}_1\|\hat{\beta}_m\|_{\Gamma}\Big\}	=O_p\Big(	\frac{ R^{1/2}  \eta_1^{-1/(4\iota)} }{M^{1/2}N^{1/2} } +K^{-\tau}+ \eta_1^{1/2} K^{(d- \nu)} \Big).
	$$
	The above can be optimized for $K$ at $K\asymp \eta_1^{-\frac{1}{2(\tau+d-\nu)}} $. Plugging in this value, it becomes
	$$
	\frac{1}{M}\sum_{m=1}^{M} \Big\{\|\hat{\beta}_m-\beta_0\|_{X} + \eta^{1/2}_1\|\hat{\beta}_m\|_{\Gamma}\Big\}	=O_p\Big(	\frac{ R^{1/2}  \eta_1^{-1/(4\iota)} }{M^{1/2}N^{1/2} } + \eta_1^{\frac{\tau}{2(\tau+d-\nu)}}\Big).
	$$
	By setting $\eta_1 \asymp (MN/R)^{-2\iota(\tau+d-\nu)/(\tau+d-\nu+2\iota\tau)}$, we get the rate of convergence as 
	$$
	\frac{1}{M}\sum_{m=1}^{M} \Big\{\|\hat{\beta}_m-\beta_0\|_{X} + \eta^{1/2}_1\|\hat{\beta}_m\|_{\Gamma}\Big\}	 =O_p\Big(		(MN/R)^{-\iota\tau/(\tau+d-\nu +2\iota\tau)}\Big).
	$$
	Correspondingly, the spline degrees of freedom $K$ satisfies $K\asymp (MN/R)^{\iota/(\tau+d-\nu+2\iota\tau)} $ in this case.
\end{proof}

\section{Technical Proofs for Section~\ref{sec:graphrate}}

\subsection{Proof of Lemma~\ref{lemma:graph:Qbound}}

\begin{proof} 
For the first conclusion~\eqref{eqn:aux:graphpenalty}, it is obvious that
\begin{equation} \label{eqn:traceToDeriv}
	\eta_1 \tr\big(\mB\trans\mGamma\mB\big) = \eta_1 \sum_{m=1}^M
\int_{\sT} \big\{\beta_m^{(d)}(t)  \big\}^2 \intd t .
\end{equation}
Meanwhile, notice  that
	\begin{align}
		&\tr\big(\mB\mOmega\mB\trans\widehat{\mSigma}\big) +\eta_1
		\tr\big(\mB\trans\mGamma\mB\mOmega\big) \nonumber \\
		=&   \sum_{v,v'=1}^M
		\frac{w_{vv'}}{NM} \sum_{n,m}
		\langle x_{nm},\beta_v -  \beta_{v'} \rangle^2   +
		\eta_1  \sum_{v,v'=1}^M
		w_{vv'}\int \big\{\beta_v^{(d)}(t) -  \beta_{v'}^{(d)}(t) \big\}^2 \intd t  \nonumber \\
		=&  \sum_{v,v'=1}^M w_{vv'}
		\Big[ \frac{1}{NM} \sum_{n,m}
		\langle x_{nm},\beta_v -  \beta_{v'} \rangle^2   +
		\eta_1 \int \big\{\beta_v^{(d)}(t) -  \beta_{v'}^{(d)}(t) \big\}^2 \intd t\Big]. \label{eqn:traceToEmNorm}
\end{align}
	Applying  Proposition~\ref{proposition:empiricalnorm}, we known each summand in the square bracket of~\eqref{eqn:traceToEmNorm} 
	converges in probability to
		$\| \beta_v -  \beta_{v'} \|_{X}^2  +\eta_1 \| \beta_v -  \beta_{v'} \|_{\Gamma}^2$. Therefore, the limit of~\eqref{eqn:traceToEmNorm} is
	\begin{align}
		 \sum_{v,v'=1}^M w_{vv'}
		\Big( \| \beta_v -  \beta_{v'} \|_X^2  +\eta_1 \| \beta_v -  \beta_{v'} \|_{\Gamma}^2\Big)  
		&=  \sum_{v,v'=1}^M w_{vv'}
		\Big( \| \vb_v -  \vb_{v'} \|_2^2  +\eta_1 \| \vb_v -  \vb_{v'} \|_{\Gamma}^2\Big)  \nonumber  \\
		&= \tr\big(\mB\mOmega\mB\trans\big) +\eta_1
		\tr\big(\mB\trans\mGamma\mB\mOmega\big).
	\end{align}
 This, together with~\eqref{eqn:traceToDeriv}, establishes the first conclusion~\eqref{eqn:aux:graphpenalty}. 
 
 Based on the above discussion, the second conclusion~\eqref{eqn:graph:normbound} is directly established by noticing that
		\begin{align} \label{eqn:graph:QboundExpand}
		\sQ_{\veta}^2(\mB) = \|\mB\|_F^2 +
		\eta_1 \tr\big(\mB\trans\mGamma\mB\big) +
		\eta_2 \tr\big(\mB\mOmega\mB\trans\widehat{\mSigma}\big) + \eta_1\eta_2
		\tr\big(\mB\trans\mGamma\mB\mOmega\big).
	\end{align}
The proof is completed.
\end{proof}

\subsection{Proof of  Lemma~\ref{lemma:graph:bias}}

The proof of  Lemma~\ref{lemma:graph:bias} depend on  some results  of \cite{trillos2020error}.  
In their Appendix, they consider the quantity 
$$
p_v=\sum_{v'=1}^M 
\frac{1}{ h^{\mu}M} G\Big(-\frac{\|\vs_v-\vs_{v'}\|_2}{h}\Big),
$$
for $v=1,\cdots, M$. We can interpret  $p_v$ as a kernel density estimate  of the true sampling density $p(\vs_{v})$  at $\vs_{v}$. Recall from Section~\ref{sec:graphrate:review} the regularity conditions on $p(\vs_{v})$.
The density $p(\vs_{v})$ is Lipschitz continuous with Lipschitz constant $L_p$, and it is bounded from below and above ($C_p > p(\vs) > 1/C_p$  for some constant $C_p>0$ and for all $\vs\in\sS$).  Equation (A.1) of \cite{trillos2020error} shows the following holds 
\begin{equation} \label{eqn:graph:bias:densityConverge}
\max_{v=1,\cdots,M} 
\big| p_v - p(\vs_{v})\big| \le C L_p h +  CC_p G(0) \frac{\epsilon}{h}
+ CC_p h^2
\end{equation}
for some constant $C$ depending on the  curvature and intrinsic dimension of the manifold $\sS$. 	
In the above, $\epsilon$ is the $\infty$-optimal transport distance between $P_n$ and $P$, where
$P_n$ is the empirical measure of $\vs_1,\cdots, \vs_M$ and $P$ is the measure whose density with respect to the volume form $\intd V_{\vs}$ is $p(\vs)$.

Theorem~2 of \cite{trillos2020error} indicates that the $\infty$-optimal transport distance
satisfies 
$\epsilon = O_p \big(\log(M)^{\zeta_\mu} /M^{1/\mu} \big)$,
where $\zeta_\mu=3/4$ if $\mu=2$ and $\zeta_\mu=1/\mu$ if $\mu\ge 3$. Because $p(\vs)$ is bounded by $C_p$ and due to~\eqref{eqn:graph:bias:densityConverge},  we have 
\begin{equation} \label{eqn:graph:bias:densityBound}
\max_{v=1,\cdots,M}  p_v = O_p(1),
\end{equation}
when $h\to 0$ and $h M^{1/\mu}/\log(M)^{\zeta_\mu} \to\infty$.


\vspace{10pt}

\noindent\textbf{Proof of Lemma~\ref{lemma:graph:bias}}.	 	
	At each $\vs_{v}$, we can construct an approximation $\tilde{\beta}_{0v}(t) = Q\beta_{0}(t,\vs_{v})$ in the spline space $\mathbb{S}_{K}$, where  the linear mapping $Q$ is defined in~(6.40) of~\cite{schumaker2007spline}.	The mapping $Q$ can control the approximation error to the optimal order. Similar to the proof of Lemma~\ref{lemma:reduced:approxerror}, the approximation error $\sE(\mathbb{S}_K)$ can be quantified by considering the constructed spline approximation $\tilde{\beta}_{0v}$. According to Lemma~\ref{lemma:graph:Qbound}, we have $\sE(\mathbb{S}_K) = O_p(E)$ where
	\begin{align} 
		E= &\sum_{m=1}^{M} \| \tilde{\beta}_{0m} - \beta_{0m}\|_{X}^2 + 
		\eta_1 \sum_{m=1}^M
		\int_{\sT} \big\{\tilde{\beta}_m^{(d)}(t)  \big\}^2 \intd t \nonumber \\
		&	 +	\eta_2 \underbrace{\sum_{v,v'=1}^M w_{vv'}
		\Expect_{x_{nm}}\langle x_{nm},\tilde{\beta}_{v} - \tilde{\beta}_{v'} \rangle^2 }_{E_1} 
		+ \eta_1 \eta_2 \underbrace{\sum_{v,v'=1}^M w_{vv'}\int \big\{\tilde{\beta}_{v}^{(d)}(t) -  \tilde{\beta}_{v'}^{(d)}(t) \big\}^2 \intd t}_{E_2}. 
		\label{eqn:penaltylimitE}
	\end{align}
	For the first two terms in the right hand side of~\eqref{eqn:penaltylimitE}, we have
	\begin{equation}
		\sum_{m=1}^{M} \| \tilde{\beta}_{0m} - \beta_{0m}\|_{X}^2 + 
		\eta_1 \sum_{m=1}^M
		\int_{\sT} \big\{\tilde{\beta}_m^{(d)}(t)  \big\}^2 \intd t  \asymp M\big\{K^{-2\tau} +
		\eta_1  K^{2(d-\nu)_{+}} \big\} 	\label{eqn:penaltylimitE_p1}
	\end{equation}
	according to Proposition~\ref{proposition:general:splineerror}. We continue to control the last two terms $E_1$ and $E_2$ of~\eqref{eqn:penaltylimitE} in the following.

	Suppose $\vs_{v'}$ is in the $2h$ geodesic neighbor of $\vs_{v}$. Denote $\gamma(t)$ as the geodesic curve connecting them with $\gamma(0) = \vs_{v'}$ and $\gamma(1) = \vs_{v}$, and $\| \cdot \|_{\vs}  = \langle\cdot, \cdot \rangle_{\vs}^{1/2}$ is the metric induced norm. Then, for the true slope function, we have
	\begin{align*}
		\big|\beta_0(t,\vs_v)-	\beta_0(t,\vs_{v'})\big| 
		&= \Big|\int_0^1  \frac{\intd}{\intd t} \beta_0(t,\gamma(t))
		\intd t \Big| \\
		&= \Big|\int_0^1 \big\langle \nabla_{\vs} \beta_0(t,\gamma(t)),\,
		\dot\gamma(t)
		\big\rangle_{\vs} \intd t \Big|\\
		&\le \int_0^1 \| \nabla_{\vs} \beta_0(t,\gamma(t))\|_{\vs}
		\cdot \|\dot\gamma(t)\|_{\vs} \intd t \\
		&\stackrel{(i)}{\le}  C_0 \int_0^1 \|\dot\gamma(t)\|_{\vs}\intd t 
		= C_0 d_{\sS}(\vs_{v'}, \vs_{v}) \le 2 C_0 h,
	\end{align*}
	where $d_{\sS}(\vs_{v'}, \vs_{v})$ represents the geodesic distance between $\vs_{v'}$ and $\vs_{v}$.
	In the above, (i) uses $\nabla_{\vs} \beta_0(t,\vs )$ is continuous over the compact set $\sT\times \sS$ in Condition~\ref{ass:betagraph}, hence its norm can be bounded by some constant $C_0$.
		It follows that, for all $t\in\sT$ and $\vs_{v},\vs_{v'}\in\sS$, we have the following uniform upper bound
	\begin{align}
		\Big| \frac{\beta_0(t,\vs_v)-	\beta_0(t,\vs_{v'})}{h}\Big| \le 2C_0.
	\end{align}

The difference of the constructed spline approximation at two distinct $\vs_{v}, \vs_{v'}$
can be expressed as
\begin{align*}
	\{\tilde{\beta}_{0v}(t) 
	-\tilde{\beta}_{0v'}(t)\}/h =Q \Big[ \frac{\beta_0(\cdot,\vs_v)-	\beta_0(\cdot,\vs_{v'})}{h} \Big].
\end{align*}
Applying Theorem~6.22 of \cite{schumaker2007spline}, we have
\begin{align} \label{eqn:graph:bias:betaDiff_1}
	\Big\| 	\frac{\{\tilde{\beta}_{0v}(t) -\tilde{\beta}_{0v'}(t) \}}{h}  \Big\|_{L_2}
	&=\Big\|	Q \Big[ \frac{\beta_0(\cdot,\vs_v)-	\beta_0(\cdot,\vs_{v'})}{h} \Big]
	\Big\|_{L_2} 
		\le  C_1:= 2C_0 \{2(\mathfrak{o}+1)\}^{\mathfrak{o}+1}.
\end{align}
Similarly, by Theorem~6.25 of \cite{schumaker2007spline}, we can show
\begin{equation}  \label{eqn:graph:bias:betaDiff_2}
	\Big\|	D^d Q \Big[ \frac{\beta_0(\cdot,\vs_v)-	\beta_0(\cdot,\vs_{v'})}{h} \Big] \Big\|_{L_2} \lesssim  K^{(d-\nu)_{+}}.
\end{equation}

Recall $\lambda_{01}$ is the largest eigenvalue of the covariance function $\mathcal{C}$ of $X_m$.
Based on the weight value $w_{vv'}$ of~\eqref{eqn:edgeweights}, 
we can control $E_1$ of~\eqref{eqn:penaltylimitE}  by
	\begin{align}
	E_1 &=	\sum_{v,v'=1}^M w_{vv'}
		\Expect_{x_{nm}}\langle x_{nm},\tilde{\beta}_{0v}
		-\tilde{\beta}_{0v'}\rangle^2  \nonumber   \\
		&\le	\sum_{v,v'=1}^M 
		\frac{2 \lambda_{01}}{\sigma_G h^{\mu}M} G\Big(-\frac{\|\vs_v-\vs_{v'}\|_2}{h}\Big)
		\cdot \| \tilde{\beta}_{0v}
		-\tilde{\beta}_{0v'}\|_{L_2}^2/h^2  \nonumber \\
		&\le  C_1\lambda_{01}	\sum_{v=1}^M  \sum_{v'=1}^M 
		\frac{2}{\sigma_G h^{\mu}M} G\Big(-\frac{\|\vs_v-\vs_{v'}\|_2}{h}\Big) \nonumber \\ 
		&=O_p(M). 	\label{eqn:penaltylimitE_p2}
	\end{align}
In the above, the second last line is due to~\eqref{eqn:graph:bias:betaDiff_1}, and the last equality is due to~\eqref{eqn:graph:bias:densityBound}.

Similarly,  we can show for $E_2$ of~\eqref{eqn:penaltylimitE}  that
\begin{align} 
	E_2 &=	 \sum_{v,v'=1}^M
		w_{vv'}\int \big\{\tilde{\beta}_v^{(d)}(t) -  \tilde{\beta}_{v'}^{(d)}(t) \big\}^2 dt \nonumber \\
		&= 	\sum_{v,v'=1}^M 
		\frac{2}{\sigma_G h^{\mu}M} K\Big(-\frac{\|\vs_v-\vs_{v'}\|_2}{h}\Big)
		\cdot \Big\|	D^d Q \Big[ \frac{\beta_0(\cdot,\vs_v)-	\beta_0(\cdot,\vs_{v'})}{h} \Big] \Big\|_{L_2}^2  \nonumber \\
		& = O_p\big( M K^{(d-\nu)_{+}} \big). 	\label{eqn:penaltylimitE_p3}
	\end{align}
	The conclusion of Lemma~\ref{lemma:graph:bias} follows by combining~\eqref{eqn:penaltylimitE}, \eqref{eqn:penaltylimitE_p1},  \eqref{eqn:penaltylimitE_p2}, and  \eqref{eqn:penaltylimitE_p3}. $\hfill \blacksquare$
%

\subsection{Proof of Lemma~\ref{lemma:graph:gamma}}
\begin{proof} Under the event~\eqref{eqn:graph:normbound}, it hold that
	$	\mathbb{N}(\bmB,\delta) $ is a subset of 	$\sN_1(\bmB,\delta) $
	\begin{align}  
		\sN_1(\bmB,\delta) = \big\{ \mDelta\in\bbR^{K\times M}:\ 	\| (\mI + \eta_1 \mGamma)^{1/2} \mDelta (\mI + \eta_2 \mOmega)^{1/2}\|_F\le2 \delta \big\}.
	\end{align}
	We can equivalently consider to bound the complexity for $\sN_1(\bmB,\delta) $. The idea is similar to the proof of~\eqref{eqn:complexity0:keybound}.
	We can recognize $\sN_1(\bmB,\delta)$ as an ellipsoid in $\bbR^{K\times M}$ by diagonalizing $\mOmega$. 	Suppose $\mOmega = \widetilde{\mX} \widetilde{\mOmega} \widetilde{\mX}\trans$	is the eigen-decomposition, and $\widetilde{\mOmega} =\diag(\tilde{w}_1,\cdots, \tilde{w}_M)$ is a diagonal matrix with $\tilde{\omega}_m = \lambda_m^{\uparrow}(\mOmega)$. Set $\widetilde{\mDelta} = \mDelta\widetilde{\mX}$, then
	$$
	\| (\mI + \eta_1 \mGamma)^{1/2} \mDelta (\mI + \eta_2 \mOmega)^{1/2}\|_F = 
	\| (\mI + \eta_1 \mGamma)^{1/2} \widetilde{\mDelta} (\mI + \eta_2 \widetilde{\mOmega} )^{1/2}\|_F .
	$$
	We get an equivalent representation of $\sN_1(\bmB,\delta)$ as
	\begin{equation}
		\sN_1(\bmB,\delta) = \Big\{ \widetilde{\mDelta} = \big(\tilde{a}_{km}\big)_{K\times M}:\; 
		\sum_{k=1}^{K} 	\sum_{m=1}^{M}  (1+\eta_1 \gamma_{k})(1+\eta_2 \tilde{\omega}_{m})  \tilde{a}_{km}^2 \le 4\delta^2\Big\}. \label{eqn:complexity2:ellipsoid}
	\end{equation}
	The half lengths of the principal axes are $2\delta/\sqrt{(1+\eta_1 \gamma_{k})(1+\eta_2 \tilde{\omega}_{m}) }$. 	Under Condition~\ref{ass:laplacianEigenvalue}, we have 
	$\tilde{\omega}_m=\lambda_m^{\uparrow}(\mOmega) \gtrsim  m^{2/\mu}$. 	Meanwhile, as in Proposition~\ref{prop:splineStruct}, $\mGamma=\mathrm{diag}(\gamma_1,\cdots, \gamma_K)$ is a diagonal matrix with $\gamma_k \gtrsim k^{(2q + 2d)}$. Therefore, based on Equation (2.115) and Theorem 4.1.11 of \citep{talagrand2014upper}, we have that
	\begin{align*}
		\gamma_{2}(E,d) &= O_p\big(	\sN_1(\bmB,\delta) \big) \\
		 &= O_p\bigg(	  \delta \Big\{ 	\sum_{k=1}^{K} 	\sum_{m=1}^{M} \frac{1}{(1+\eta_1 \gamma_{k})(1+\eta_2 \tilde{\omega}_{m})} \Big\}^{1/2}  \bigg) \\
		&= O_p\Big(\delta  \Big\{ \sum_{k=1}^{K} 	\frac{1}{1+\eta_1 \gamma_{k}} \Big\}^{1/2} \cdot  \Big\{ 		\sum_{m=1}^{M} \frac{1}{1+\eta_2 \tilde{\omega}_{m}} \Big\}^{1/2}  \Big)\\
        &= O_p\Big(	 \delta (M^{1/2}\wedge \eta_2^{-\mu/4})\cdot \big\{ K^{1/2}\wedge \eta_1^{-1/(4q+4d)} \big\}\Big).
	\end{align*}
The last inequality uses result~\eqref{eqn:complexity0:keybound} and the similar bound $\sum_{m=1}^M \frac{1}{1+ \eta_2 \tilde{\omega}_m}  \lesssim M \wedge \eta_2^{-\mu/2}$.
\end{proof}

\subsection{Proof of Corollary~\ref{corollary:graph:rate1}}

Consider the case of weak graph regularization $\eta_2 \lesssim M^{-2/\mu}$.
Suppose $d \le \nu$, $\eta_1 \rightarrow 0$, and $M,N\to \infty$. The upper bound of the convergence rate can be written as 
$$
(1/M) \Big\{\sum_{m=1}^{M} \|\hat{\beta}_m-\beta_{0m}\|_{X} +  \sP_{\veta}^{1/2}(\hat{\vbeta})\Big\} =O_p\Big(  	\frac{K^{1/2} \wedge \eta_1^{-1/(4\iota)}}{N^{1/2}} + K^{-\tau} + \eta_1^{1/2} + \eta_2^{1/2}\Big).
$$
Specifically, we will consider three subcases, which are determined by the relative magnitudes of $K$, $\eta_1^{-\frac{1}{2(\iota \vee \tau)}}$, and $\eta_1^{-\frac{1}{2(\iota \wedge \tau)}}$.

\noindent \textbf{Subcase (i)}.  Consider the case when $K \lesssim \eta_1^{-\frac{1}{2(\iota \vee \tau)}} \le \eta_1^{-\frac{1}{2(\iota \wedge \tau)}}$. The
upper bound of the convergence rate in this case simplifies to 
$$
(1/M) \Big\{\sum_{m=1}^{M} \|\hat{\beta}_m-\beta_{0m}\|_{X} +  \sP_{\veta}^{1/2}(\hat{\vbeta})\Big\} =O_p\Big( \frac{K^{1/2}}{N^{1/2}} + K^{-\tau} + \eta_2^{1/2}\Big).$$
After optimizing the right hand side by $K$, we find the minimum is achieved at $K \asymp N^{1/(2\tau + 1)}$. Therefore, we have
$$
(1/M) \Big\{\sum_{m=1}^{M} \|\hat{\beta}_m-\beta_{0m}\|_{X} +  \sP_{\veta}^{1/2}(\hat{\vbeta})\Big\} =O_p\Big(  N^{-\tau/(2\tau + 1)}\Big),$$
when $K \asymp N^{1/(2\tau + 1)}$, $\eta_1 \lesssim N^{-2(\iota\vee \tau)/(2\tau + 1)}$ and $\eta_2 \lesssim N^{-2\tau/(2\tau + 1)}$. \\

\noindent \textbf{Subcase (ii)}. Consider the case when $\eta_1^{-\frac{1}{2(\iota \vee \tau)}}  \le  \eta_1^{-\frac{1}{2(\iota \wedge \tau)}}\lesssim K$. 
The rate of convergence in this case is upper bounded by
$$
(1/M) \Big\{\sum_{m=1}^{M} \|\hat{\beta}_m-\beta_{0m}\|_{X} +  \sP_{\veta}^{1/2}(\hat{\vbeta})\Big\} =O_p\Big(   \frac{\eta_1^{-1/(4\iota)}}{N^{1/2}} + \eta_1^{1/2} + \eta_2^{1/2}\Big).
$$
After optimizing the right hand side with respect to $\eta_1$, we also find the optimal value is achieved at
$\eta_1 \asymp N^{-2\iota/(2\iota+1)}$. Therefore, we obtain the rate of convergence
$$
(1/M) \Big\{\sum_{m=1}^{M} \|\hat{\beta}_m-\beta_{0m}\|_{X} +  \sP_{\veta}^{1/2}(\hat{\vbeta})\Big\} =O_p\Big(  N^{-\iota/(2\iota+1)}\Big),
$$
when $\eta_1 \asymp N^{-2\iota/(2\iota+1)}$, $\eta_2 \lesssim N^{-2\iota/(2\iota+1)} $, and $K \gtrsim N^{\frac{\iota}{(\iota\wedge\tau)(2\iota+1)}}$. \\

\noindent \textbf{Subcase (iii)}. Consider the case when $\eta_1^{-\frac{1}{2(\iota \vee \tau)}} \lesssim K \lesssim  \eta_1^{-\frac{1}{2(\iota \wedge \tau)}}$. In the following, we will show  Subcase (iii) either corresponds to Subcase~(i) or  Subcase~(ii), depending on the relative size of $\tau$ and $\iota$.

\begin{enumerate}
	\item  When $\tau \le \iota $, the upper bound of the convergence rate is then of order 
	$$
	(1/M) \Big\{\sum_{m=1}^{M} \|\hat{\beta}_m-\beta_{0m}\|_{X} +  \sP_{\veta}^{1/2}(\hat{\vbeta})\Big\} =O_p\Big(  \frac{\eta_1^{-1/(4\iota)}}{N^{1/2}} + K^{-\tau} + \eta_2^{1/2}\Big).
	$$
	Optimizing the right hand side with respect to $K$ in the range $\eta_1^{-\frac{1}{2(\iota \vee \tau)}} \lesssim K \lesssim  \eta_1^{-\frac{1}{2(\iota \wedge \tau)}}$, we find the optimal value is achieved when $K$ has the order of its upper bound in this range, i.e.,
	$K\asymp  \eta_1^{-\frac{1}{2(\iota \wedge \tau)}}\asymp  \eta_1^{-\frac{1}{2\tau}}$. Plugging in this value to the upper bound of the convergence rate, we obtain
	\begin{equation} \label{proof:graphrate1:31}
		(1/M) \Big\{\sum_{m=1}^{M} \|\hat{\beta}_m-\beta_{0m}\|_{X} +  \sP_{\veta}^{1/2}(\hat{\vbeta})\Big\} =O_p\Big(  \frac{\eta_1^{-1/(4\iota)}}{N^{1/2}} + \eta_1^{1/2} + \eta_2^{1/2}\Big).
	\end{equation}
	
	\item When $\tau > \iota $, the upper bound of the convergence rate is then of order 
	$$
	(1/M) \Big\{\sum_{m=1}^{M} \|\hat{\beta}_m-\beta_{0m}\|_{X} +  \sP_{\veta}^{1/2}(\hat{\vbeta})\Big\} =O_p\Big(  \frac{K^{1/2}}{N^{1/2}} + \eta_1^{1/2} + \eta_2^{1/2}\Big).
	$$
	Optimizing the right hand side with respect to $K$ in the range $\eta_1^{-\frac{1}{2(\iota \vee \tau)}} \lesssim K \lesssim  \eta_1^{-\frac{1}{2(\iota \wedge \tau)}}$, we find that the optimal value is achieved when $K$ has the order of its lower bound in this range, i.e.,
	$K\asymp  \eta_1^{-\frac{1}{2(\iota \vee \tau)}}\asymp  \eta_1^{-\frac{1}{2\tau}}$. Plugging in this value to the upper bound of the convergence rate, we obtain
	\begin{equation} \label{proof:graphrate1:32}
		(1/M) \Big\{\sum_{m=1}^{M} \|\hat{\beta}_m-\beta_{0m}\|_{X} +  \sP_{\veta}^{1/2}(\hat{\vbeta})\Big\} =O_p\Big(  \frac{\eta_1^{-{1}/{(4 \tau)}}}{N^{1/2}} + \eta_1^{1/2} + \eta_2^{1/2}\Big).
	\end{equation}
\end{enumerate}
The above two rates of convergence~\eqref{proof:graphrate1:31} and~\eqref{proof:graphrate1:32}
can be summarized as
$$
(1/M) \Big\{\sum_{m=1}^{M} \|\hat{\beta}_m-\beta_{0m}\|_{X} +  \sP_{\veta}^{1/2}(\hat{\vbeta})\Big\} =O_p\Big(  \frac{\eta_1^{-{1}/{(4\iota \vee 4\tau)}}}{N^{1/2}} + \eta_1^{1/2} + \eta_2^{1/2}\Big).
$$
In the above, the optimal $\eta_1$ is of order
$\eta_1 \asymp N^{-\frac{2(\iota \vee \tau)}{2(\iota \vee \tau) +1}}$, which then leads to
$$
(1/M) \Big\{\sum_{m=1}^{M} \|\hat{\beta}_m-\beta_{0m}\|_{X} +  \sP_{\veta}^{1/2}(\hat{\vbeta})\Big\} =O_p\Big(  N^{-\frac{(\iota \vee \tau)}{2(\iota \vee \tau) +1}}\Big).
$$
In summary, the above rate of convergence for Subcase~(iii) is obtained when $\eta_1 \asymp N^{-\frac{2(\iota \vee \tau)}{2(\iota \vee \tau) +1}}$, $\eta_2 \lesssim N^{-\frac{2(\iota \vee \tau)}{2(\iota \vee \tau) +1}}$, and $K \asymp N^{\frac{(\iota \vee \tau)/\tau}{2(\iota \vee \tau) +1}}$. 
It can be seen that Subcase~(iii) exactly corresponds to Subcases~(i) and (ii), when $\tau\ge \iota$ and $\tau<\iota$, respectively. 
$\hfill \blacksquare$

\subsection{Proof of Corollary~\ref{corollary:graph:rate2}}
In the case of strong graph regularization where $\eta_2 \gtrsim M^{-2/\mu}$, $N,M \rightarrow \infty$, and $\eta_1\to 0$, the upper bound of the convergence rate can be written as
$$
(1/M) \Big\{\sum_{m=1}^{M} \|\hat{\beta}_m-\beta_{0m}\|_{X} +  \sP_{\veta}^{1/2}(\hat{\vbeta})\Big\} 
=O_p\Big( 	\frac{ \eta_2^{-\mu/4} \{K^{1/2} \wedge \eta_1^{-1/(4\iota)}\} }{ M^{1/2} N^{1/2}} + K^{-\tau} + \eta_1^{1/2} + \eta_2^{1/2} \Big).
$$
Specifically, we will consider three subcases, which are determined by the relative magnitudes of $K$, $\eta_1^{-\frac{1}{2(\iota \vee \tau)}}$, and $\eta_1^{-\frac{1}{2(\iota \wedge \tau)}}$.

\noindent \textbf{Subcase (i)}.  Consider the case when $K \lesssim \eta_1^{-\frac{1}{2(\iota \vee \tau)}} \le \eta_1^{-\frac{1}{2(\iota \wedge \tau)}}$. The
upper bound of the convergence rate in this case simplifies to 
$$
(1/M) \Big\{\sum_{m=1}^{M} \|\hat{\beta}_m-\beta_{0m}\|_{X} +  \sP_{\veta}^{1/2}(\hat{\vbeta})\Big\} 
=O_p\Big( 		\frac{ \eta_2^{-\mu/4} K^{1/2}  }{ M^{1/2} N^{1/2}} + K^{-\tau}  + \eta_2^{1/2} \Big).
$$
Fixing $K$ and optimizing the right hand side with respect to the parameter $\eta_2$, the minimum is achieved when 
$$
\eta_2 \asymp K^{\frac{2}{\mu+2}} (MN)^{-\frac{2}{\mu+2}}.
$$
Plugging in this value into the upper bound of the convergence rate, we obtain
$$
(1/M) \Big\{\sum_{m=1}^{M} \|\hat{\beta}_m-\beta_{0m}\|_{X} +  \sP_{\veta}^{1/2}(\hat{\vbeta})\Big\} 
=O_p\Big( 	 K^{-\tau}  + K^{\frac{1}{\mu+2}} (MN)^{-\frac{1}{\mu+2}}\Big).
$$	
Taking further optimization of the right hand side with respect to $K$, it can be directly seen that the minimum is obtained at
$$
K\asymp (MN)^{\frac{1}{\tau (\mu+2)+1}}.
$$
Therefore, we have
$$
	(1/M) \Big\{\sum_{m=1}^{M} \|\hat{\beta}_m-\beta_{0m}\|_{X} +  \sP_{\veta}^{1/2}(\hat{\vbeta})\Big\} 
	=O_p\Big(  (NM)^{-\frac{\tau}{\tau(2+\mu) + 1}}\Big).
$$
In summary, this rate of convergence for Subcase~(i) is achieved when $\eta_1 \lesssim (MN)^{-\frac{2(\iota\vee\tau)}{\tau(2+\mu) + 1}}$, $\eta_2 \asymp (MN)^{-\frac{2\tau}{\tau(2+\mu) + 1}}$, and $K \asymp (MN)^{\frac{1}{\tau(2+\mu) + 1}}$.\\

\noindent \textbf{Subcase (ii)}. Consider the case when $\eta_1^{-\frac{1}{2(\iota \vee \tau)}}  \le  \eta_1^{-\frac{1}{2(\iota \wedge \tau)}}\lesssim K$.  The upper bound of the convergence rate in this case can be organized as 
$$
(1/M) \Big\{\sum_{m=1}^{M} \|\hat{\beta}_m-\beta_{0m}\|_{X} +  \sP_{\veta}^{1/2}(\hat{\vbeta})\Big\} 
=O_p\Big( 		\frac{ \eta_2^{-\mu/4} \eta_1^{-1/4\iota}  }{ M^{1/2} N^{1/2}} + \eta_1^{1/2}  + \eta_2^{1/2}\Big).
$$
We then optimize the right hand side of the above with respect to $\eta_1$ and $\eta_2$.
For a fixed $\eta_1$, it is directly seen that the optimal $\eta_2$ satisfies
$$
\eta_2 \asymp \eta_{1}^{-\frac{1}{\iota (\mu+2)}} (MN)^{-\frac{2}{\mu+2}},
$$
and plug in this value into the rate of convergence to obtain
$$
(1/M) \Big\{\sum_{m=1}^{M} \|\hat{\beta}_m-\beta_{0m}\|_{X} +  \sP_{\veta}^{1/2}(\hat{\vbeta})\Big\} 
=O_p\Big(  \eta_1^{1/2}  +  \eta_{1}^{-\frac{1}{2\iota (\mu+2)}} (MN)^{-\frac{1}{\mu+2}} \Big).
$$
Taking further optimization with respect to $\eta_1$, we get
$$
(1/M) \Big\{\sum_{m=1}^{M} \|\hat{\beta}_m-\beta_{0m}\|_{X} +  \sP_{\veta}^{1/2}(\hat{\vbeta})\Big\} 
\asymp  (MN)^{-\frac{\iota}{\iota(2+\mu) + 1}},
$$
where the tuning parameters are configured as $\eta_1 \asymp (MN)^{-\frac{2\iota}{\iota(2+\mu) + 1}}$, $\eta_2 \asymp (MN)^{-\frac{2\iota}{\iota(2+\mu) + 1}}$, and $K \gtrsim (MN)^{\frac{\iota}{\{\iota(2+\mu) + 1\} (\iota \wedge \tau)}}$. \\

\noindent \textbf{Subcase (iii)}. Consider the case when $\eta_1^{-\frac{1}{2(\iota \vee \tau)}} \lesssim K \lesssim  \eta_1^{-\frac{1}{2(\iota \wedge \tau)}}$. In the following, we will show  Subcase (iii) either corresponds to Subcase~(i) or Subcase~(ii), depending on the relative size of $\tau$ and $\iota$. 

\begin{enumerate}
	\item If $\tau\le \iota$, the rate of convergence is then upper bounded by 
	$$
	(1/M) \Big\{\sum_{m=1}^{M} \|\hat{\beta}_m-\beta_{0m}\|_{X} +  \sP_{\veta}^{1/2}(\hat{\vbeta})\Big\} 
	=O_p\Big(  \frac{ \eta_2^{-\mu/4} \eta_1^{-1/4\iota}  }{ M^{1/2} N^{1/2}} + K^{-\tau}  + \eta_2^{1/2} \Big).
	$$	
	In the range $\eta_1^{-\frac{1}{2(\iota \vee \tau)}} \lesssim K \lesssim  \eta_1^{-\frac{1}{2(\iota \wedge \tau)}}$, the optimal $K$ for the right hand side of the above reaches at $K\asymp  \eta_1^{-\frac{1}{2 \tau}}$. This implies the upper bound of the convergence rate becomes 
	\begin{equation} \label{proof:graphrate2:31}
		(1/M) \Big\{\sum_{m=1}^{M} \|\hat{\beta}_m-\beta_{0m}\|_{X} +  \sP_{\veta}^{1/2}(\hat{\vbeta})\Big\} 
		=O_p\Big(  \frac{ \eta_2^{-\mu/4} \eta_1^{-1/4\iota}  }{ M^{1/2} N^{1/2}} + \eta_1^{1/2}  + \eta_2^{1/2} \Big).
	\end{equation}
	It shows that the above has the same form as in Subcase~(ii).  
	
	
	\item If $\tau> \iota$, the rate of convergence is then upper bounded by 
	$$
	(1/M) \Big\{\sum_{m=1}^{M} \|\hat{\beta}_m-\beta_{0m}\|_{X} +  \sP_{\veta}^{1/2}(\hat{\vbeta})\Big\} 
	=O_p\Big(  \frac{ \eta_2^{-\mu/4} K^{1/2}  }{ M^{1/2} N^{1/2}} + \eta_1^{1/2}  + \eta_2^{1/2}\Big).
	$$	
	In the range $\eta_1^{-\frac{1}{2(\iota \vee \tau)}} \lesssim K \lesssim  \eta_1^{-\frac{1}{2(\iota \wedge \tau)}}$, the optimal $K$ for the right hand side of the above reaches at $K\asymp  \eta_1^{-\frac{1}{2 \tau}}$. This implies the upper bound of the convergence rate becomes 
	\begin{equation}  \label{proof:graphrate2:32}
		(1/M) \Big\{\sum_{m=1}^{M} \|\hat{\beta}_m-\beta_{0m}\|_{X} +  \sP_{\veta}^{1/2}(\hat{\vbeta})\Big\} 
	=O_p\Big(  \frac{ \eta_2^{-\mu/4} K^{1/2}  }{ M^{1/2} N^{1/2}} + K^{-\tau}  + \eta_2^{1/2}\Big).
	\end{equation}
	It shows that the above has the same form as in Subcase~(i). 
%
\end{enumerate}
The results of this corollary follow by combining Subcases (i)--(iii).  $\hfill \blacksquare$

\vskip 0.2in

\end{document}